\documentclass[12pt]{amsart}
\usepackage{amsmath,amsthm,amsfonts,epsfig,amssymb}
\usepackage{dsfont}

\usepackage[colorlinks]{hyperref}
\usepackage[alphabetic,initials]{amsrefs}

\theoremstyle{plain}
\newtheorem{thm}{Theorem}[section]
\newtheorem{prop}[thm]{Proposition}

\newtheorem{cor}[thm]{Corollary}
\newtheorem{lemma}[thm]{Lemma}
\newtheorem{assumption}[thm]{Assumption}

\theoremstyle{definition}
\newtheorem{example}[thm]{Example}

\newtheorem{defn}[thm]{Definition}

\theoremstyle{remark}
\newtheorem{remark}[thm]{Remark}

\newcommand{\dist}{\operatorname{dist}}
\newcommand{\End}{\operatorname{End}}
\newcommand{\ev}{\operatorname{ev}}

\newcommand{\Hom}{\operatorname{Hom}}

\newcommand{\ind}{\operatorname{ind}}
\newcommand{\inj}{\operatorname{injrad}}

\newcommand{\interior}{\operatorname{int}}

\newcommand{\loc}{\operatorname{loc}}
\newcommand{\lk}{\operatorname{lk}}
\newcommand{\muCZ}{\mu_{\text{CZ}}}
\newcommand{\piTu}{\pi Tu}

\newcommand{\sign}{\operatorname{sgn}}

\newcommand{\wind}{\operatorname{wind}}
\newcommand{\windpi}{\operatorname{wind}_\pi}

\newcommand{\SL}{\operatorname{SL}}

\newcommand{\CC}{{\mathbb C}}
\newcommand{\DD}{{\mathbb D}}

\newcommand{\HH}{{\mathbb H}}
\newcommand{\NN}{{\mathbb N}}
\newcommand{\QQ}{{\mathbb Q}}
\newcommand{\RR}{{\mathbb R}}
\newcommand{\ZZ}{{\mathbb Z}}

\newcommand{\dD}{{\mathcal D}}

\newcommand{\fF}{{\mathcal F}}
\newcommand{\hH}{{\mathcal H}}

\newcommand{\mM}{{\mathcal M}}

\newcommand{\pP}{{\mathcal P}}
\newcommand{\tT}{{\mathcal T}}
\newcommand{\uU}{{\mathcal U}}

\newcommand{\1}{\mathds{1}}
\newcommand{\p}{\partial}
\newcommand{\Cinftyloc}{C^\infty_{\loc}}

\numberwithin{equation}{section}

\title[Foliations on Overtwisted Contact Manifolds]{Finite Energy Foliations
on Overtwisted Contact Manifolds}
\author{Chris Wendl}
\address{Departement Mathematik, HG G38.1 \\ 
R\"amistrasse 101 \\
8092 Z\"urich \\ 
Switzerland}
\email{wendl@math.ethz.ch}
\urladdr{http://www.math.ethz.ch/~wendl/}
\thanks{Research partially supported by NSF grants DMS-0102298
and DMS-0603500}

\begin{document}

\begin{abstract}
We develop a method for preserving pseudoholomorphic curves in contact
$3$--manifolds under surgery along transverse links.  This makes use of a
geometrically natural boundary value problem for holomorphic curves in
a $3$--manifold with stable Hamiltonian structure, where the boundary
conditions are defined by $1$--parameter families of 
totally real surfaces.  The technique
is applied here to construct a finite energy
foliation for every closed overtwisted contact $3$--manifold.
\end{abstract}

\maketitle

\tableofcontents

\newpage

\section{Introduction}
\label{sec:intro}

\subsection{Definitions and main result}
\label{subsec:mainresult}

Finite energy foliations of contact $3$--manifolds were introduced by 
Hofer, Wysocki and Zehnder in
\cite{HWZ:foliations}, where they were shown to exist for generic tight
three-spheres, giving powerful consequences for the Reeb dynamics.  The
present work is a step toward generalizing such existence results:
we prove that for every closed overtwisted contact $3$--manifold, one
can choose a contact form (of Morse-Bott type) and compatible complex 
multiplication such that a stable finite energy foliation exists.  

Fix a closed oriented $3$--manifold $M$ with a cooriented, positive
contact structure
$\xi$: this is by definition the kernel of a smooth $1$--form $\lambda$
which satisfies $\lambda \wedge d\lambda > 0$.  The choice of contact
form defines the \emph{Reeb vector field} $X$ by the conditions
$$
d\lambda(X,\cdot) \equiv 0
\quad \text{ and } \quad
\lambda(X) \equiv 1.
$$
Then the flow of $X$ preserves $\xi$, along with the symplectic
vector bundle structure on $\xi \to M$ defined by $d\lambda$.

Recall that a contact structure $\xi$ is called \emph{overtwisted} if
there exists an \emph{overtwisted disk}: an embedded disk $\dD \subset M$ 
such that for all
$m \in \p\dD$, $T_m(\p\dD) \subset \xi_m$ but $T_m\dD \ne \xi_m$.
By Eliashberg's classification result \cite{Eliashberg:overtwisted},
contactomorphism classes of overtwisted contact structures on $M$ are in
one-to-one correspondence with homotopy classes of cooriented 
$2$--plane distributions.

The following is the main result of this paper.  

\begin{thm}
\label{thm:mainresult}
Suppose $(M,\xi)$ is a closed oriented $3$--manifold with a positive overtwisted
contact structure.  Then there exists a contact form $\lambda$ and
admissible complex multiplication $J$ such that $(M,\lambda,J)$ admits a
stable finite energy foliation of Morse-Bott type.  The foliation has
precisely one nondegenerate asymptotic orbit and one or more Morse-Bott
tori of asymptotic orbits, and every leaf is either an orbit cylinder or an 
index~$2$ finite energy sphere with distinct
simply covered asymptotic orbits, all positive.
\end{thm}

We will spend the rest of \S\ref{subsec:mainresult} explaining the 
definitions needed to understand this statement.

Denote the time-$t$ flow of $X$ by $\varphi_X^t$, and
recall that a closed orbit $x : \RR \to M$ of $X$ with period 
$T > 0$ is called \emph{nondegenerate} if the linearized return map
$d\varphi_X^T(x(0))|_{\xi_{x(0)}}$ does not have~$1$ in its spectrum.
More generally, a \emph{Morse-Bott manifold of $T$--periodic orbits}
is a submanifold $N \subset M$ tangent to $X$ such that
$\varphi_X^T|_N$ is the identity, and for every $m \in N$,
$$
T_m N = \ker (d\varphi_X^T(m) - \1).
$$
In this paper we shall deal exclusively with situations where $N$ is a
circle (i.e. a nondegenerate orbit) or a two-dimensional torus.  For
the latter case, one can show (see \cite{Wendl:BP1}) 
that all orbits in $N$ have the same
minimal period $\tau > 0$, and $N$ is a Morse-Bott family of $k\tau$--periodic
orbits for all $k \in \NN$.  Thus we will call such submanifolds
\emph{Morse-Bott tori} without reference to the period, and a single
closed orbit will be called simply \emph{Morse-Bott} if it either is
nondegenerate or belongs to a Morse-Bott torus.

The \emph{symplectization} of $M$ is the open $4$--manifold $\RR\times M$
with symplectic structure $d(e^a\lambda)$, where $a$ denotes the coordinate
on the $\RR$--component.  We consider a natural class of
$\RR$--invariant almost complex structures compatible with this symplectic
form, defined as follows.  Note the choice of contact form $\lambda$ defines
a splitting
$$
T(\RR\times M) = (\RR \oplus \RR X) \oplus \xi,
$$
where the first factor also comes with a natural trivialization.
An \emph{admissible complex multiplication} is a choice of
complex structure $J$ for the bundle $\xi \to M$, compatible with the
symplectic structure, i.e.~we require that $d\lambda(\cdot, J\cdot)$ define
a bundle metric.  Given any such choice, we define an almost complex
structure $\tilde{J}$ on $\RR\times M$ in terms of the above splitting
and trivialization by
$$
\tilde{J} = i \oplus J,
$$
where $i$ is understood as the natural complex structure acting on
$\CC = \RR^2$.  We will call $\tilde{J}$ the \emph{almost complex 
structure associated to $\lambda$ and $J$}.  

Given such a structure, we consider $\tilde{J}$--holomorphic curves
$$
\tilde{u} = (a,u) : (\dot{\Sigma},j) \to (\RR\times M, \tilde{J}),
$$
where the domain $\dot{\Sigma} = \Sigma \setminus \Gamma$ is a Riemann 
surface $(\Sigma,j)$ with a discrete set of points 
$\Gamma\subset \Sigma$ removed.  The \emph{energy} of such a curve is
defined as
\begin{equation}
\label{eqn:energy}
E(\tilde{u}) = \sup_{\varphi\in\tT} \int_{\dot{\Sigma}} 
 \tilde{u}^* d(\varphi\lambda),
\end{equation}
where $\tT := \{ \varphi\in C^\infty(\RR,[0,1])\ |\ \varphi' \ge 0 \}$.
An easy computation shows that the integrand is nonnegative whenever
$\tilde{u}$ is $\tilde{J}$--holomorphic, and such a curve is constant if
and only if $E(\tilde{u}) = 0$.  When $\Sigma$ is closed,
$\tilde{J}$--holomorphic curves $\tilde{u} : \dot{\Sigma} \to \RR\times M$
with $E(\tilde{u}) < \infty$ are called \emph{finite energy surfaces}.
By results in \cites{Hofer:weinstein,HWZ:props1,HWZ:props4}, these have
nicely controlled asymptotic behavior near the punctures, which can be
described as follows.  Denote by $\DD \subset \CC$ the
closed unit disk with its natural complex structure, and let
$\DD_r \subset \CC$ be the closed disk of radius $r$ for any $r > 0$.

\begin{prop}
\label{prop:asymptotics}
Suppose $\tilde{u} = (a,u) : \dot{\DD} = \DD\setminus\{0\} \to \RR\times M$
is a $\tilde{J}$--holomorphic map with 
$0 < E(\tilde{u}) < \infty$.  If $\tilde{u}$ is bounded, then $\tilde{u}$
extends to a $\tilde{J}$--holomorphic map $\DD\to\RR\times M$.  Otherwise,
$\tilde{u}$ is a proper map, and
for every sequence $s_k \to \infty$ there is a subsequence
such that the loops $t\mapsto u(e^{-2\pi (s_k+it)})$ converge
in $C^\infty(S^1,M)$ to a loop $t\mapsto x(Qt)$.  Here 
$x : \RR\to M$
is a periodic orbit of $X$ with period $T = |Q|$, where
\begin{equation}
\label{eqn:charge}
Q = - \lim_{\epsilon \to 0} \int_{\p \DD_\epsilon} u^*\lambda \ne 0.
\end{equation}
Moreover, $t\mapsto a(e^{-2\pi (s_k+it)})/s_k$ converges
in $C^\infty(S^1,\RR)$ to the constant map $t\mapsto Q$.

If the orbit $x$ is Morse-Bott, then in fact the maps
$t \mapsto u(e^{-2\pi(s+it)})$ and $t \mapsto a(e^{-2\pi (s+it)})/s$
converge in $C^\infty(S^1)$ as $s \to \infty$.
\end{prop}
The number $Q \in \RR\setminus\{0\}$ appearing above is called the
\emph{charge} of the puncture, and we call the puncture positive/negative
in accordance with the sign of $Q$.  This defines a partition of the
set of punctures:
$$
\Gamma = \Gamma^+ \cup \Gamma^-,
$$
and one can use the maximum principle to show that finite energy
surfaces always have $\#\Gamma^+ \ge 1$, cf.~\cite{HWZ:props2}.

The simplest example of a finite energy surface is the so-called
\emph{orbit cylinder} or \emph{trivial cylinder} over a $T$--periodic
orbit $x : \RR\to M$.  Indeed, it's
easy to check that the map
$$
\tilde{u} : \RR\times S^1 \to \RR\times M : (s,t) \mapsto (Ts, x(Tt))
$$
is $\tilde{J}$--holomorphic and has finite energy; after reparametrization,
it is a sphere with one positive puncture and one negative.  
Prop.~\ref{prop:asymptotics} above is a precise way of saying that any
finite energy surface looks approximately like an orbit cylinder near
each puncture.

\begin{defn}
\label{defn:fef}
A \emph{finite energy foliation} for $(M,\lambda,J)$ is a smooth 
two-dimensional foliation $\fF$ of $\RR\times M$ such that
\begin{enumerate}
\item Each leaf $F\in\fF$ can be presented as the image of an embedded
$\tilde{J}$--holomorphic finite energy 
surface, and there exists a constant that bounds the energy of every leaf
uniformly.
\item For every leaf $F\in\fF$, the set 
$\sigma + F := \{ (\sigma + a,m) \ | \ (a,m)\in F \}$ for $\sigma\in\RR$ 
is also a leaf 
of the foliation, and thus either disjoint from or identical to $F$.
\end{enumerate}
\end{defn}

We shall often abuse notation and write $\tilde{u} \in \fF$, meaning that
the finite energy surface $\tilde{u}$ parametrizes a leaf of $\fF$.
The $\RR$--invariance assumption says that $\tilde{u} = (a,u) \in \fF$ if 
and only if $\tilde{u}^\sigma := (a+\sigma,u) \in \fF$ for all $\sigma\in \RR$.
This has several consequences for the projection of $\fF$ to the underlying
contact manifold.

\begin{prop}
\label{prop:foliationProperties}
Let $\fF$ be a finite energy foliation.  Then
\begin{enumerate}
\item[(i)] If $P \subset M$ is a periodic orbit which is an asymptotic limit 
for some leaf $\tilde{u} \in \fF$, then the orbit cylinder
$\RR\times P$ is also a leaf of $\fF$.
\item[(ii)] For each leaf $\tilde{u} = (a,u) : \dot{\Sigma} \to \RR\times M$
 of $\fF$ that is not an orbit cylinder, 
 the map $u : \dot{\Sigma} \to M$ is embedded and does not
 intersect its asymptotic limits.
\item[(iii)] If $\tilde{u} = (a,u) : \dot{\Sigma} \to \RR\times M$ and
 $\tilde{v} = (b,v) : \dot{\Sigma}' \to \RR\times M$ are two leaves of $\fF$, 
 then $u(\dot{\Sigma})$ and $v(\dot{\Sigma}')$ are 
either disjoint or identical.
\end{enumerate}
\end{prop}
The proofs of these properties are mostly straightforward exercises using
positivity of intersections; we refer to \cite{Wendl:thesis} for details.
The only detail not covered there is the fact that the maps
$u : \dot{\Sigma} \to M$ are not just injective but also \emph{embedded}: for
this one uses intersection theory to show that a critical point of
$u$ at $z \in \dot{\Sigma}$ implies intersections between 
$(a,u)$ and $(a + \epsilon,u)$ near $z$ for small $\epsilon$;
cf.~\cite{Wendl:BP1}.
Denote by $\pP_\fF \subset M$
the union of all the closed orbits that occur as asymptotic limits for
leaves of $\fF$; equivalently, this is the projection down to $M$ of all
the orbit cylinders in $\fF$.  Then Prop.~\ref{prop:foliationProperties}
can be rephrased as follows.

\begin{cor}
If $\fF$ is a finite energy foliation, then the projections of its leaves 
from $\RR\times M$ to $M$ form a smooth foliation of $M\setminus \pP_\fF$.
\end{cor}

To explain the stability condition, we need to introduce some more
technical details.
In the following, write periodic orbits of $X$ in the notation
$P := (x(\RR),T)$ where $x : \RR \to M$ is a $T$--periodic solution of
the Reeb flow equation; note that $T$ need not be the \emph{minimal} period,
and the parametrization $x$ can always be changed by a time shift.  We
shall sometimes abuse notation and regard $P$ as a subset of $M$, keeping
in mind that $T$ is also part of the data.  Recall that if $P$ is 
nondegenerate, we can choose a unitary trivialization $\Phi$ of $\xi$ along
$P$ and define the \emph{Conley-Zehnder index}
$\muCZ^\Phi(P)$ as in \cite{HWZ:props2}.  Then a finite energy surface
$\tilde{u} : \dot{\Sigma} \to \RR\times M$ with
only nondegenerate asymptotic orbits $\{ P_z \}_{z \in \Gamma}$ is assigned
the Conley-Zehnder index
$$
\muCZ(\tilde{u}) = \sum_{z\in\Gamma^+} \muCZ^\Phi(P_z) -
\sum_{z\in\Gamma^-} \muCZ^\Phi(P_z),
$$
where the trivializations $\Phi$ are chosen so as to admit an extension to
a global complex trivialization of $\xi$ along $\tilde{u}$; 
then $\muCZ(\tilde{u})$ doesn't depend on this choice.

This index can be extended to the Morse-Bott case in the following
straightforward manner.  Given $P$, define an \emph{admissible parametrization}
of $P$ to be any map $\mathbf{x} : S^1 \to P \subset M$ such that
$\lambda(\dot{\mathbf{x}}) \equiv T$.
This defines the so-called \emph{asymptotic operator}
\begin{equation}
\label{eqn:asympOp}
\mathbf{A}_{\mathbf{x}} : \Gamma(\mathbf{x}^*\xi) \to
\Gamma(\mathbf{x}^*\xi) : v \mapsto -J (\nabla_t v - T \nabla_v X),
\end{equation}
where $\nabla$ is any symmetric connection on $M$; one can check that
this expression gives a well defined section of $\mathbf{x}^*\xi$, not
depending on $\nabla$.  As an unbounded operator on
$L^2(\mathbf{x}^*\xi)$ with domain $H^1(\mathbf{x}^*\xi)$,
$\mathbf{A}_{\mathbf{x}}$ is self-adjoint, with spectrum consisting of
discrete real eigenvalues of finite multiplicity, accumulating only at
infinity.  The equation $\mathbf{A}_{\mathbf{x}} v = 0$ then defines
the linearized Reeb flow restricted to $\xi$ along $P$, 
and $P$ is nondegenerate if and only
if $\ker\mathbf{A}_{\mathbf{x}} = \{0\}$.  When this is the case, one can
define the linearized Reeb flow purely in terms of the equation
$\mathbf{A}_{\mathbf{x}} v = 0$ and sensibly denote the Conley-Zehnder index by
$\muCZ^\Phi(P) = \muCZ^\Phi(\mathbf{A}_{\mathbf{x}})$.  
The key observation now is
that for any $c \in \RR$, the equation $(\mathbf{A}_{\mathbf{x}} - c) v = 0$
also defines a linear Hamiltonian flow, and thus yields a well defined
Conley-Zehnder index if $c$ is not an eigenvalue of $\mathbf{A}_{\mathbf{x}}$.
Then if $P$ belongs to a 
Morse-Bott manifold $N \subset M$, we can pick any sufficiently
small number $\epsilon > 0$ and define the \emph{perturbed} Conley-Zehnder
indices
\begin{equation}
\label{eqn:perturbedCZ}
\muCZ^{\Phi\pm}(P) = \muCZ^\Phi(\mathbf{A}_{\mathbf{x}} \pm \epsilon).
\end{equation}
This doesn't depend on $\epsilon$ if the latter is sufficiently small,
but does depend on the sign choice whenever $\ker\mathbf{A}_{\mathbf{x}}$
is nontrivial, i.e.~when $P$ is degenerate.
For $\tilde{u}$ with Morse-Bott asymptotic orbits, we now define its
Conley-Zehnder index by
$$
\muCZ(\tilde{u}) = \sum_{z\in\Gamma^+} \muCZ^{\Phi-}(P_z) -
\sum_{z\in\Gamma^-} \muCZ^{\Phi+}(P_z)
$$
which is equal to the previous definition if all $P_z$ are nondegenerate.

The \emph{moduli space of finite energy surfaces} $\mM_{\tilde{J}}$ is the 
set of equivalence classes $C = [(\Sigma,j,\Gamma,\tilde{u})]$, where 
$\tilde{u} : (\Sigma\setminus\Gamma,j) \to (\RR\times M,\tilde{J})$ 
is a finite energy surface, $\Gamma$ is assigned an \emph{ordering}, and
$(\Sigma,j,\Gamma,\tilde{u}) \sim (\Sigma',j',\Gamma',\tilde{u}')$ if and
only if there is a biholomorphic map $\varphi : (\Sigma,j) \to (\Sigma',j')$
that takes $\Gamma$ to $\Gamma'$ with ordering preserved and satisfies
$\tilde{u} = \tilde{u}' \circ \varphi$.  We shall sometimes abuse notation
and write $\tilde{u} \in \mM_{\tilde{J}}$ when there is no confusion.

To define a topology on $\mM_{\tilde{J}}$,
first note that the punctured Riemann surface
$(\dot{\Sigma},j)$ can be regarded as a surface with cylindrical ends, which
then admits a natural compactification.  Indeed, for each $z \in \Gamma^\pm$,
pick a closed disk-like neighborhood $\dD_z$ of $z$ in $\Sigma$ and a 
biholomorphic map $\dot{\dD}_z := \dD\setminus\{z\} \to Z_\pm$, where
\begin{equation}
\label{eqn:halfCylinders}
Z_+ = [0,\infty)\times S^1, \qquad
Z_- = (-\infty,0]\times S^1,
\end{equation}
both with the standard complex structure $i$.  This decomposes $\dot{\Sigma}$
in the form
$$
\dot{\Sigma} \cong \Sigma_0 \cup 
\left( \bigcup_{z\in\Gamma^\pm} Z_\pm \right),
$$
where $\Sigma_0$ is a compact surface with boundary.  We now define the
compactified surface $\overline{\Sigma}$
by adding ``circles at infinity,'' which means replacing
each $Z_\pm$ with $\overline{Z}_\pm$, where
$$
\overline{Z}_+ = [0,\infty]\times S^1, \qquad
\overline{Z}_- = [-\infty,0]\times S^1.
$$
Denote the components of $\p\overline{\Sigma}$ by 
$\delta_z \cong \{\pm\infty\}\times S^1$ for each $z \in \Gamma^\pm$.
We shall not place a smooth structure on $\overline{\Sigma}$.
It is naturally a compact \emph{topological} manifold with boundary,
where the interior $\dot{\Sigma} \subset \overline{\Sigma}$ 
and the boundary components $\delta_z$ all have natural smooth structures;
in fact the latter have natural identifications with $S^1$ up to translation,
and one can show that none of this structure depends on the choices.

The symplectization $W := \RR\times M$ also has a natural compactification
$\overline{W} := [-\infty,\infty]\times M$, which we again regard as a
topological manifold with boundary, on which the interior and the boundary
separately have natural smooth structures.  It is then convenient to
observe that any
finite energy surface $\tilde{u} : \dot{\Sigma} \to W$ extends naturally
to a continuous map $\bar{u} : \overline{\Sigma} \to \overline{W}$,
whose restriction to each $\delta_z$ gives an admissible parametrization
of the corresponding orbit $P_z \subset \{\pm\infty\}\times M$.

We say that a sequence $C_k \in \mM_{\tilde{J}}$ converges to 
$C \in \mM_{\tilde{J}}$ if for
sufficiently large $k$ there exist
representatives $(\Sigma,j_k,\Gamma,\tilde{u}_k) \in C_k$ and
$(\Sigma,j,\Gamma,\tilde{u}) \in C$ such that
\begin{enumerate}
\setlength{\itemsep}{0in}
\item $j_k \to j$ in $C^\infty(\Sigma)$
\item $\tilde{u}_k \to \tilde{u}$ in $\Cinftyloc(\dot{\Sigma},W)$,
\item $\bar{u}_k \to \bar{u}$ in $C^0(\overline{\Sigma},\overline{W})$.
\end{enumerate}
This defines the topology on the moduli space $\mM_{\tilde{J}}$.

For any finite energy surface $\tilde{u} : \dot{\Sigma} \to \RR\times M$
with Morse-Bott asymptotic orbits, define the \emph{index} of $\tilde{u}$ by
\begin{equation}
\label{eqn:FredholmIndexBasic}
\ind(\tilde{u}) = \muCZ(\tilde{u}) - \chi(\dot{\Sigma}).
\end{equation}
This is the Fredholm index of the \emph{linearized normal Cauchy-Riemann
operator} $\mathbf{L}_{\tilde{u}}$, which is explained in
\cite{HWZ:props3} for the nondegenerate case and \cite{Wendl:BP1} in
general.  We call $\tilde{u}$ \emph{regular} if the operator
$\mathbf{L}_{\tilde{u}}$ is surjective; in this case the implicit function
theorem allows us to describe a neighborhood of $\tilde{u}$ in 
$\mM_{\tilde{J}}$ as
a smooth manifold of dimension 
$\ind(\tilde{u})$.  In the general Morse-Bott case there is a stronger
notion of regularity: suppose $\tilde{u}$ has a puncture $z \in \Gamma$
at which the asymptotic orbit belongs to a Morse-Bott torus $N \subset M$,
and let $\uU_{\tilde{u}} \subset \mM_{\tilde{J}}$ denote a connected open
neighborhood of $\tilde{u}$.
We can assume without loss of generality that all curves in $\uU_{\tilde{u}}$ 
are parametrized on the same domain $\Sigma$ with the puncture 
$z \in \Gamma \subset \Sigma$
in a fixed position.  The Reeb flow along $N$
defines an $S^1$--action so that $N / S^1$ is a circle, and there is then
a natural evaluation map
$$
\ev_z : \uU_{\tilde{u}} \to N / S^1,
$$
defined by assigning to any curve in $\uU_{\tilde{u}}$ its asymptotic orbit 
$P_z \subset N$.
We say that $\tilde{u}$ is \emph{strongly regular} if it is regular and
for every $z \in \Gamma$ where $P_z$ is degenerate,
$\ev_z$ has a surjective linearization at $\tilde{u}$.

\begin{defn}
\label{defn:stability}
A finite energy foliation $\fF$ is called \emph{stable} if $\pP_\fF$ is
a finite union of nondegenerate Reeb orbits, and every leaf $F\in \fF$
is parametrized by a \emph{regular} finite energy surface 
$\tilde{u}_0 \in \mM_{\tilde{J}}$
such that all other curves $\tilde{u} \in \mM_{\tilde{J}}$ 
near $\tilde{u}_0$ also parametrize leaves of $\fF$.

We say that $\fF$ is a \emph{stable} foliation of \emph{Morse-Bott type}
if $\pP_\fF$ is a finite union of nondegenerate Reeb orbits and Morse-Bott
tori, and each leaf is a \emph{strongly regular} finite energy surface
whose neighbors in $\mM_{\tilde{J}}$ are also leaves of $\fF$.
\end{defn}

Observe that leaves of stable finite energy foliations can
only have index~$0$, $1$ or~$2$.  The index~$0$ leaves are precisely the
orbit cylinders, while index~$1$ leaves are called \emph{rigid surfaces},
because they project to isolated leaves in the foliation of
$M \setminus \pP_\fF$.  Index~$2$ leaves come in $\RR$--invariant
$2$--parameter families, which project to $1$--parameter families in
$M \setminus \pP_\fF$.  In the Morse-Bott case, orbit cylinders can also
have index~$1$, projecting to $M$ as $1$--parameter families moving along
Morse-Bott tori.

\begin{figure}
\begin{center}
\includegraphics{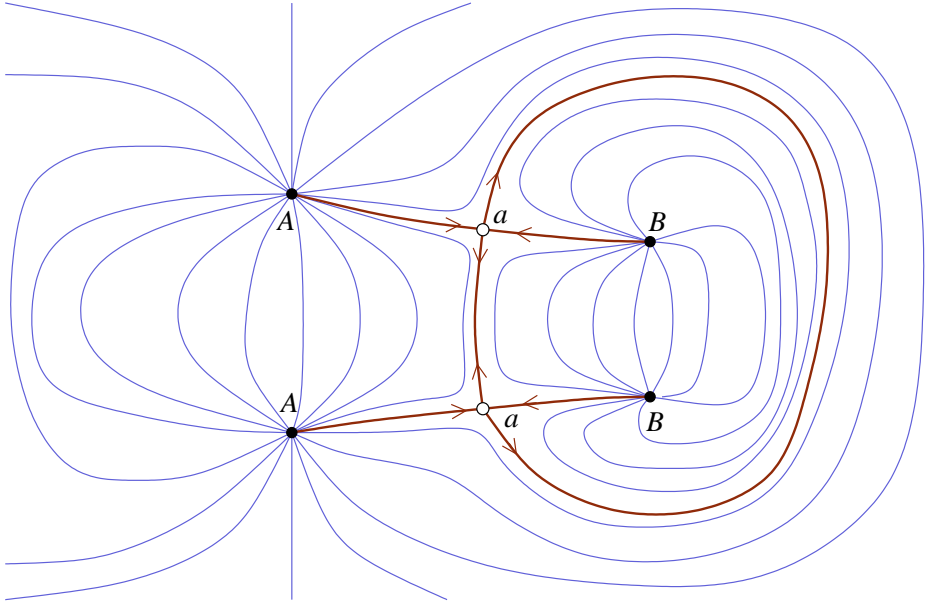}
\end{center}
\caption[Stable finite energy foliation of 
 $S^3$]{\label{fig:fol_S3} A cross section of a stable finite energy
foliation on $S^3 = \RR^3 \cup \{\infty\}$, with three asymptotic orbits
cutting transversely through the page.  The hyperbolic orbit $a$
is the limit of two rigid planes, and is connected to two elliptic
orbits $A$ and $B$ by rigid cylinders.  All other leaves are
index~2 planes asymptotic to $A$ or $B$.
Arrows represent the signs of the punctures at $a$:
a puncture is positive/negative if the arrow points away from/toward
the orbit.}
\end{figure}

\subsection{Outline of the proof}

The surgery construction involves two main technical ingredients.
The first is the Fredholm and intersection theory for the mixed boundary
value problem considered in \cite{Wendl:BP1}, which we review 
in~\S\ref{sec:MBVP}.  The crucial
point is to observe that embedded index~$2$ curves
with certain properties are always regular, and give rise to
non-intersecting $2$--parameter families of embedded curves, which
project to $1$--parameter families of embeddings in $M$.  

\begin{figure}
\begin{center}
\includegraphics{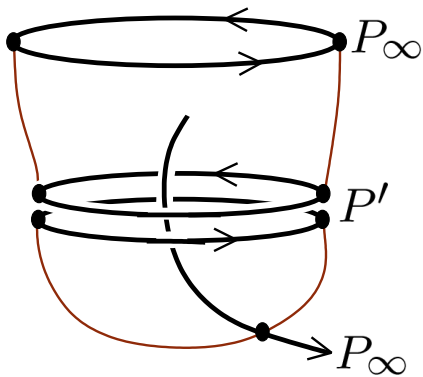}
\end{center}
\caption[A topological obstructions to noncompactness]{
\label{fig:bubbling} A finite energy plane bubbles off and produces an
illegal intersection with the asymptotic limit $P_\infty$.  This gives
a topological obstruction to noncompactness.}
\end{figure}

The second
main ingredient is a compactness argument: this is explained in
\S\ref{sec:compactness} and constitutes the bulk of the technical work
in this paper.  While we borrow certain ideas from the compactness
theorems of Symplectic Field Theory \cite{SFTcompactness}, those
results cannot be applied and are in fact not true in our setup, because
we make only very weak nondegeneracy assumptions on our data.
This is necessary in order to accommodate nontrivial homotopies of the data,
but it allows potentially quite strange asymptotic behavior for holomorphic
curves.  Thus in our situation, the moduli space generally has no natural
compactification---yet we'll 
find that the particular spaces of interest encounter topological obstructions
to noncompactness, which are peculiar to the
low-dimensional setting.  In this way, our arguments are quite
different from those in \cite{SFTcompactness}.  This is illustrated by
the example in Figure~\ref{fig:bubbling}.  Here we consider the
degeneration of a sequence of finite energy planes $\tilde{u}_k$
in $\RR\times S^3$, all asymptotic to the same simply covered orbit 
$P_\infty \subset S^3$ and not intersecting it.
Without assuming that the Reeb vector field is nondegenerate, it is
sometimes possible to show that any other closed orbit $P \subset S^3$
must be nontrivially linked with $P_\infty$.  Then if a plane bubbles off
as in the picture, its asymptotic limit $P'$ must be linked with $P_\infty$,
implying that the new plane intersects $P_\infty$.  But then $P_\infty$
must also intersect $\tilde{u}_k$ for sufficiently large~$k$, giving a
contradiction.  Some more elaborate variations on this argument will be used 
repeatedly in \S\ref{sec:compactness}.

The main result is then proved in \S\ref{sec:mainConstruction}: starting
with a stable finite energy foliation on the tight three-sphere
(an open book decomposition with one nondegenerate binding orbit), 
we perform a combination of rational Dehn surgery and Lutz twists along a 
transverse link $K \subset S^3$ and show that the resulting contact manifold 
also admits a stable foliation, now of Morse-Bott type.  The topological
preliminaries on surgery and Lutz twists are explained in
\S\ref{subsec:DehnLutz}.  In \S\ref{subsec:local}, we tackle the easiest
step in the foliation construction, finding families of holomorphic curves
to fill the solid tori that are glued in by surgery.  
This is done by explicitly solving
the nonlinear Cauchy-Riemann equation on $S^1 \times B^2 \subset S^1 \times
\RR^2$, with rotationally symmetric contact forms and complex structures.

\begin{figure}
\begin{center}
\includegraphics{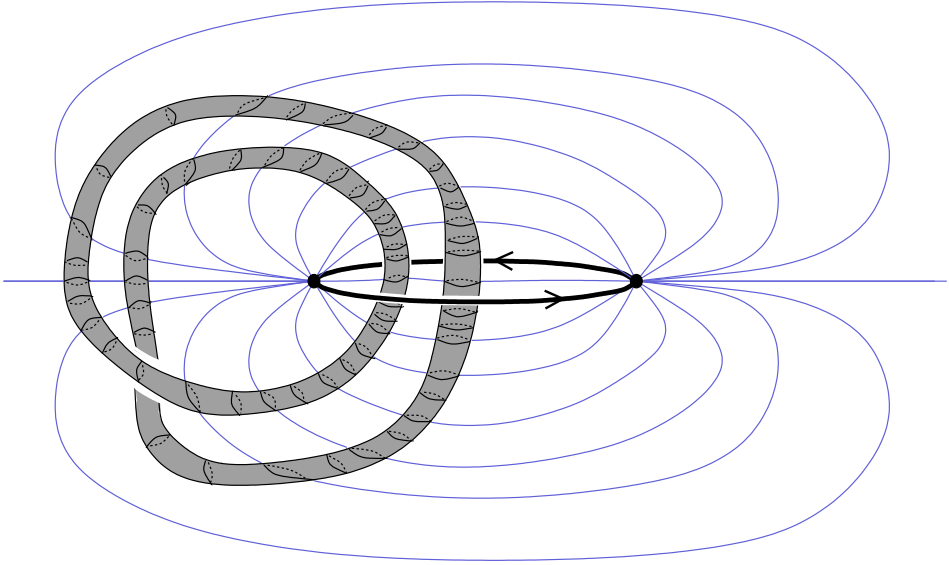}
\end{center}
\caption[Stable Morse-Bott foliation after surgery]{\label{fig:MorseBottFol} 
The stable Morse-Bott finite energy foliation obtained from an open book
decomposition of $S^3$ after surgery along a transverse knot $K$, linked twice
with the binding orbit.  Each leaf in the region outside the surgery has
three punctures: one at the original binding orbit, and two at Morse-Bott
orbits along the torus around $K$.}
\end{figure}

The technical background of \S\ref{sec:MBVP} and \S\ref{sec:compactness}
is then applied in \S\ref{subsec:surgery} to change a given open book
decomposition of $S^3$ (which can also be constructed explicitly) into
a stable Morse-Bott foliation in the complement of a transverse link
neighborhood.  The key is to 
remove a collection of disks from each page of the open book, 
obtaining holomorphic curves with boundary, which satisfy a problem of the type
considered in \S\ref{sec:MBVP}, with images avoiding the region inside 
a set of small tori.
Thus we are now free to perform surgery and Lutz twists inside 
these tori without
killing the holomorphic curves outside; the region inside can afterwards be
filled in by the explicitly constructed curves from \S\ref{subsec:local}.

In principle, the curves with boundary filling the region outside the tori can
be homotoped as the contact form is twisted, so that in the
limit, all boundary components degenerate to punctures, giving rise to
finite energy surfaces without boundary.  The actual argument is not quite
so simple, because the curves obtained by cutting out disks 
generally satisfy a boundary 
condition that is \emph{totally real} but not \emph{Lagrangian}, thus
lacking \textit{a priori} energy bounds.  This problem does not appear
to be solvable for all curves at once, but we can deal with a single curve
in the case where each component of $K$ is only singly linked with the
binding orbit: this makes it possible to construct the totally real surface
in $\RR\times M$ so that it is equivalent via a diffeomorphism to a
Lagrangian surface in the symplectization of $M$ with a \emph{stable
Hamiltonian structure}, i.e.~the generalization of a contact form
described in \cite{SFTcompactness}.  We can then homotop this back to
contact data but keep the Lagrangian boundary condition, and use the
implicit function theorem to extend the single curve again to a whole
foliation.  The only catch is that this trick requires the restrictive
assumption that components of $K$ link only once with the binding orbit.
We fix this in \S\ref{subsec:braids} by viewing the general case as a
branched cover: the foliation can then be lifted to the cover using
intersection theory.

\subsection{Discussion}
\label{subsec:discussion}

The result proved here is one step in a program proposed by Hofer
to study Reeb dynamics on arbitrary closed contact $3$--manifolds via
finite energy foliations; this is joint work in progress by the author
with Hofer, R.~Siefring and J.~Fish.  As was shown in \cite{HWZ:foliations}, 
the existence of finite energy foliations on $(S^3,\xi_0)$ implies that generic
Reeb vector fields in that setting admit either two or infinitely many
periodic orbits.  The present work does not imply such a result for overtwisted
contact manifolds, because we fix a very specific contact form.  The
next step would therefore be a homotopy argument in which one shows that
a foliation of $(M,\lambda,J)$ gives rise to a foliation (or something
similar) of
$(M,f\lambda,J')$ for generic positive smooth functions $f$ and complex
multiplications $J'$.  One can then try to extend this to tight
contact manifolds by the following trick: any $(M,\xi)$ can be made 
overtwisted by taking a connected sum of $(M,\xi)$ with 
$(S^3,\xi_{\text{ot}})$ for some
$\xi_{\text{ot}}$ overtwisted.  One would then like to understand what
happens to a sequence of foliations on the overtwisted object as one
pinches off the overtwisted part.  It is known that this program \emph{cannot}
in general
lead to a stable finite energy foliation for generic $(M,\lambda,J)$, as
there are examples of tight contact manifolds where stable foliations
don't exist (cf.~\cite{Wendl:BP1}).  Nonetheless, the limits obtained in
such spaces from sequences of foliations should be interesting objects,
with possible dynamical implications.

A related program is of a more topological
nature: the author proposed in \cite{Wendl:thesis}*{Chapter~6} an equivalence
relation for stable finite energy foliations, called \emph{concordance},
which is defined by the existence of
stable holomorphic foliations on cylindrical symplectic cobordisms.
The goal would then be to classify all foliations for a given $(M,\xi)$ up to
concordance.  Conjecturally, two concordance classes can be
distinguished by an invariant $HC_*(\fF)$, which is a version of
\emph{contact homology} (or more generally, symplectic field theory
\cite{SFT}) that counts the orbits and
rigid surfaces in the foliation.  In this framework, the constructions
of the present paper show that every overtwisted $(M,\xi)$ admits a foliation
$\fF$ for which $HC_*(\fF)$ is trivial.  As suggested however by an example
in \cite{Wendl:thesis}, this is not true for \emph{all} 
foliations on overtwisted contact manifolds.

We should mention the related work of 
C.~Abbas \cite{Abbas:openbook}, which uses the open book decompositions
of Giroux \cites{Giroux:ICM,Giroux:openbook} to produce 
(in the planar case)
open book decompositions with pseudoholomorphic pages.  Due to Etnyre's
result that all overtwisted contact structures are planar
\cite{Etnyre:planar}, this also produces a finite energy foliation for
all overtwisted contact manifolds.  The two constructions are however
quite different, e.g.~ours is not an open book decomposition, and the
foliations of Abbas appear to have nontrivial contact homology in the sense
described above, suggesting that they are not equivalent to ours
via concordance.

\subsubsection*{Acknowledgements}
Much of the material covered here appeared originally---in more detailed 
and sometimes more complicated forms---in my Ph.D. 
thesis \cite{Wendl:thesis}, and I thank Helmut Hofer for many valuable
suggestions and insights.
I'm grateful also to Denis Auroux, Kai Cieliebak, Yasha Eliashberg, 
Brett Parker, Richard Siefring and Katrin Wehrheim for helpful conversations.

\section[Boundary value problem]{A mixed boundary value problem}
\label{sec:MBVP}

In this section we review the basic facts about the boundary value problem
considered in \cite{Wendl:BP1}, referring to that paper for all proofs.

Let $M$ be a closed, oriented $3$--manifold.  A \emph{stable Hamiltonian
structure} on $M$ is a tuple $\hH = (\xi,X,\omega,J)$ where
\begin{itemize}
\setlength{\itemsep}{0in}
\item $\xi$ is a smooth cooriented $2$--plane distribution on $M$
\item $\omega$ is a smooth closed $2$--form on $M$ which restricts to a 
 symplectic structure on the vector bundle $\xi \to M$
\item $X$ is a smooth vector field which is transverse to $\xi$, satisfies
 $\omega(X,\cdot) \equiv 0$, and whose flow preserves $\xi$
\item $J$ is a smooth complex structure on the bundle $\xi \to M$, compatible
 with $\omega$ in the sense that $\omega(\cdot, J\cdot)$ defines a
 bundle metric
\end{itemize}
It follows from these definitions that the flow of $X$ also preserves the
symplectic structure defined by $\omega$ on $\xi$, and the special
$1$--form $\lambda$ associated to $\xi$ and $X$ by the conditions
$$
\lambda(X) \equiv 1,
\qquad
\ker\lambda \equiv \xi,
$$
satisfies $d\lambda(X,\cdot) \equiv 0$.

\begin{example}
\label{ex:contact}
Given a contact form $\lambda$ on $M$ with contact structure 
$\xi = \ker\lambda$, Reeb vector field $X$ and admissible complex 
multiplication $J$, the data $\hH = (\xi,X,d\lambda,J)$ define
a stable Hamiltonian structure.
\end{example}

The stable Hamiltonian structure of Example~\ref{ex:contact} is referred
to as \emph{the contact case}.  We will define and use a non-contact
example in \S\ref{subsec:surgery}.

An $\RR$--invariant almost complex
structure $\tilde{J}$ on $\RR\times M$ is associated to
any $\hH = (\xi,X,\omega,J)$ by defining
$\tilde{J}\p_a = X$ and $\tilde{J}v = Jv$
for $v \in \xi$, where again $a$ denotes the coordinate on the $\RR$--factor
and $\p_a$ is the unit vector in the $\RR$--direction.  Since $d\lambda|_\xi$
may now be degenerate, we generalize the definition of
energy for pseudoholomorphic curves $\tilde{u} = (a,u) : (\Sigma,j) \to
(\RR\times M,\tilde{J})$ by
$$
E(\tilde{u}) = E_\omega(\tilde{u}) + E_\lambda(\tilde{u}),
$$
where
$$
E_\omega(\tilde{u}) = \int_{\Sigma} u^*\omega
$$
is the so-called \emph{$\omega$--energy}, and
$$
E_\lambda(\tilde{u}) = \sup_{\varphi\in\tT} \int_{\Sigma} 
\tilde{u}^*(d\varphi \wedge \lambda),
$$
with $\tT$ defined as in \eqref{eqn:energy}.  In the contact case this
definition is equivalent to \eqref{eqn:energy}, in the sense that uniform
bounds on one imply uniform bounds on the other.  Punctured 
$\tilde{J}$--holomorphic curves with finite energy in this generalized sense
also have asymptotically cylindrical behavior near punctures, the same as
in Prop.~\ref{prop:asymptotics}.  The next result, which follows from
arguments in \cites{HWZ:props1,HWZ:props4,Siefring:thesis}, 
gives a more precise and
useful statement.  Recall from \eqref{eqn:halfCylinders} the definition
of the positive and negative half-cylinders~$Z_\pm$.

\begin{prop}
\label{prop:asympFormula}
Let $\hH = (\xi,X,\omega,J)$ be a stable Hamiltonian structure on $M$ with
associated almost complex structure $\tilde{J}$, and choose a metric
on $M$.  Suppose $\tilde{u} = (a,u) : Z_\pm \to \RR\times M$ 
is a finite energy 
$\tilde{J}$--holomorphic map asymptotic (with sign corresponding to the
choice of $Z_+$ or $Z_-$) to a Morse-Bott orbit $P \subset M$
with admissible parametrization $\mathbf{x} : S^1 \to M$.  There is then
a smooth 
map $h : Z_\pm \to \mathbf{x}^*\xi$ with $h(s,t) \in \xi_{\mathbf{x}(t)}$
such that, up to translation by constants in $s$ and $t$, 
$u(s,t) = \exp_{\mathbf{x}(t)}h(s,t)$ for $|s|$ sufficiently large.
Moreover, either $h(s,t) \equiv 0$ or it satisfies the formula
\begin{equation}
\label{eqn:asympFormula}
h(s,t) = e^{\mu s} (\eta(t) + r(s,t))
\end{equation}
where $\eta \in \Gamma(\mathbf{x}^*\xi)$ is an eigenfunction of the
asymptotic operator $\mathbf{A}_{\mathbf{x}}$ in \eqref{eqn:asympOp},
$\mu \ne 0$ is the corresponding eigenvalue, whose sign is opposite
that of the puncture, and $r(s,t) \to 0$ uniformly in~$t$ 
as $s \to \pm\infty$.
\end{prop}

\begin{defn}
The section $e \in \Gamma(\mathbf{x}^*\xi)$ appearing in
\eqref{eqn:asympFormula} is called the \emph{asymptotic eigenfunction}
of $\tilde{u}$ at the puncture; 
it is well defined up to a positive multiple.  Note
that $e(t)$ is never zero, so given a trivialization $\Phi$ of 
$\mathbf{x}^*\xi$, there is a well defined winding number
$\wind^\Phi(e) \in \ZZ$.
\end{defn}

For some integer $m \ge 0$, let $L_1,\ldots,L_m \subset M$ be a collection 
of smoothly embedded surfaces
which are everywhere tangent to $X$.  Choosing smooth functions
$g_j : L_j \to \RR$, the graphs
\begin{equation}
\label{eqn:graphs}
\tilde{L}_j := \{ (g_j(x),x) \in \RR\times M\ |\ x \in L_j \}
\end{equation}
are then totally real submanifolds of $(\RR\times M,\tilde{J})$, and so
are their $\RR$--translations
$$
\tilde{L}_j^\sigma := \{ (g_j(x) + \sigma,x) \in \RR\times M\ |\ x \in L_j \}
$$
for every $\sigma \in \RR$.  Denote 
$\Lambda = (\tilde{L}_1,\ldots,\tilde{L}_m)$.
Then we define the moduli space
$$
\mM_{\hH,\Lambda}
$$
to consist of equivalence classes $[(\Sigma,j,\Gamma,\tilde{u})]$ where
$(\Sigma,j)$ is a compact Riemann surface with an ordered set of $m$ 
boundary components
$$
\p\Sigma = \gamma_1 \cup \ldots \cup \gamma_m,
$$
$\Gamma \subset \interior{\Sigma}$ is a finite ordered set of \emph{interior}
points giving rise to the punctured Riemann surface
$\dot{\Sigma} = \Sigma\setminus\Gamma$ with boundary
$\p\dot{\Sigma} = \p\Sigma$,
and $\tilde{u} : (\dot{\Sigma},j) \to (\RR\times M,\tilde{J})$
is a pseudoholomorphic map with $E(\tilde{u}) < \infty$ and satisfying
the following boundary condition:
$$
\begin{minipage}{4.5in}
\textit{For each component $\gamma_j \subset \p\Sigma$, there exists a number
$\sigma_j \in \RR$ 
such that $\tilde{u}(\gamma_j) \subset \tilde{L}_j^{\sigma_j}$.}
\end{minipage}
$$
As before, we will abuse notation and write $\tilde{u} \in \mM_{\hH,\Lambda}$
whenever convenient.
Equivalence in $\mM_{\hH,\Lambda}$ is defined by biholomorphic maps
that preserve the ordering of both the punctures and the boundary
components, and the definition of convergence used in \S\ref{subsec:mainresult}
to topologize $\mM_{\tilde{J}}$ also naturally gives a topology on 
$\mM_{\hH,\Lambda}$.  Note that in the case where $\hH$ is contact and
$m = 0$ (i.e.~$\p\Sigma = \emptyset$), $\mM_{\hH,\Lambda}$ is simply
$\mM_{\tilde{J}}$, the space of $\tilde{J}$--holomorphic finite energy surfaces.
Observe also that $\mM_{\hH,\Lambda}$ admits a natural $\RR$--action
$$
\RR\times \mM_{\hH,\Lambda} \to \mM_{\hH,\Lambda} :
(\sigma, (a,u)) \mapsto \tilde{u}^\sigma := (a + \sigma,u).
$$

The Conley-Zehnder index generalizes to 
$\tilde{u} = (a,u) \in \mM_{\hH,\Lambda}$
as follows.  For each puncture $z \in \Gamma$, let $P_z$ be the corresponding
orbit of $X$ approached by $\tilde{u}$.  The bundle $\xi$ along $P_z$ has
a symplectic structure $\omega$, which is preserved by the linearized
flow of $X$, thus we can choose a unitary trivialization $\Phi$ and
define $\muCZ^{\Phi\pm}(P_z)$ as in \eqref{eqn:perturbedCZ}.
Now at each component $\gamma_j \subset \p\Sigma$,
$u^*\xi$ has a totally real subbundle
$\ell \to \gamma_j$ defined for $z \in \gamma_j$ by
$$
\ell_z = \xi_{u(z)} \cap T_{u(z)}L_j;
$$
here we exploit the fact that $X$ is tangent to $L_j$, hence $\xi$ and
$L_j$ are transverse.  Then if $m + \#\Gamma > 0$, we can choose
trivializations $\Phi$ at each
orbit and each boundary component so that these extend to a global
complex trivialization of $u^*\xi \to \dot{\Sigma}$, and define
the \emph{Maslov index of $\tilde{u}$} as
$$
\mu(\tilde{u}) = \sum_{z \in \Gamma^+} \muCZ^{\Phi-}(P_z) -
\sum_{z \in \Gamma^-} \muCZ^{\Phi+}(P_z) +
\sum_{j=1}^m \mu^\Phi(u^*\xi|_{\gamma_j},\ell),
$$
where the last term is a sum of Maslov indices for the loops of totally
real subspaces $\ell$ along $\gamma_j$ with respect to the
complex trivialization $\Phi$.  If $m = \#\Gamma = 0$, so $\dot{\Sigma}$
is closed, we instead define
$$
\mu(\tilde{u}) = 2 c_1(u^*\xi).
$$
The \emph{Fredholm index of $\tilde{u}$} is then the integer
\begin{equation}
\label{eqn:FredholmIndex}
\ind(\tilde{u}) = \mu(\tilde{u}) - \chi(\dot{\Sigma}) + m.
\end{equation}

As with finite energy surfaces, a neighborhood of $\tilde{u} \in
\mM_{\hH,\Lambda}$ can be described via a linearized normal 
Cauchy-Riemann operator $\mathbf{L}_{\tilde{u}}$,
and the previous definitions of \emph{regular} and \emph{strongly regular}
carry over immediately.  

We now collect some useful results about the moduli space $\mM_{\hH,\Lambda}$.
The first important observation is that $\tilde{u} = (a,u) 
\in \mM_{\hH,\Lambda}$
is never equivalent to its $\RR$--translations $\tilde{u}^\sigma = 
(a + \sigma,u)$ for small
$\sigma$ unless it is tangent everywhere to $\RR \oplus \RR X \subset
T(\RR\times M)$, which implies $E_\omega(\tilde{u}) = 0$.  Thus
when $E_\omega(\tilde{u}) > 0$, $\tilde{u}$ and $\tilde{u}^\sigma$ have
only isolated, positive intersections.  The infinitesimal version of this
statement deals with the zeros of the section
$$
\piTu : \dot{\Sigma} \to \Hom_\CC(T\dot{\Sigma},u^*\xi),
$$
which composes the tangent map $Tu$ with the projection $\pi : TM \to \xi$
along $X$.  Counting these zeros algebraically defines the nonnegative 
homotopy invariant $\windpi(\tilde{u})$, which is a half-integer in general 
because zeros at the boundary must be counted with half weight.  
The map $u : \dot{\Sigma} \to M$ is immersed if and only if
$\windpi(\tilde{u}) = 0$, and in this case the Cauchy-Riemann equation
implies it is also transverse to $X$.  Due to
Prop.~\ref{prop:asympFormula} and the relations proved in \cite{HWZ:props2}
between the spectrum of the asymptotic operator at an orbit and the 
orbit's Conley-Zehnder index, $\windpi(\tilde{u})$ is also bounded
from above:

\begin{thm}
\label{thm:windpi}
For any $\tilde{u} \in \mM_{\hH,\Lambda}$ 
with $E_\omega(\tilde{u}) > 0$,
$$
0 \le 2\windpi(\tilde{u}) \le \ind(\tilde{u}) - 2 + 2g + \#\Gamma_0,
$$
where $g$ denotes the genus of $\Sigma$ and
$\Gamma_0$ is the set of $z \in \Gamma^\pm$ for which
$\muCZ^{\Phi\mp}(P_z)$ is even (for any trivialization $\Phi$).
\end{thm}

Nondegenerate \emph{elliptic} orbits always have odd Conley-Zehnder index,
thus punctures at such orbits never belong to $\Gamma_0$.  The next lemma
gives a useful criterion in the degenerate Morse-Bott case.  It follows,
for instance, that $z \not\in\Gamma_0$ if $P_z$ belongs to a Morse-Bott
torus which never intersects the image of $u$.

\begin{lemma}
\label{lemma:MBparity}
Suppose $\tilde{u} = (a,u) \in \mM_{\hH,\Lambda}$ is asymptotic at 
$z \in \Gamma^\pm$
to an orbit $P_z$ belonging to a Morse-Bott torus $N_z \subset M$.
Let $\Phi_0$ be the natural trivialization of $\xi$ along
$P_z$ determined by the directions tangent to $N_z$, and suppose that the
asymptotic eigenfunction $e$ of $\tilde{u}$ at $z$
satisfies $\wind^{\Phi_0}(e) = 0$.  Then $\muCZ^{\Phi_0\mp}(P_z) = \pm 1$.
\end{lemma}

Analogously to $\piTu$, sections in $\ker\mathbf{L}_{\tilde{u}}$ have
only positive zeros and the count of these satisfies a similar upper bound.
When this bound is zero in particular, we find
$\dim\ker\mathbf{L}_{\tilde{u}} \le 2$, leading to
the conclusion in the embedded index~$2$ case that $\tilde{u}$ is
regular without need for any genericity assumption.  Moreover,
nearby curves in the moduli space can be identified with sections in
$\ker\mathbf{L}_{\tilde{u}}$ which have no zeros, and therefore the 
nearby curves have no
intersections except possibly near infinity.  The latter can
be excluded \textit{a priori} in the situation of interest to us here, 
and we obtain the following special case of a result in \cite{Wendl:BP1}:

\begin{thm}
\label{thm:IFT}
Suppose $\tilde{u} = (a,u) : \dot{\Sigma} \to \RR\times M$ 
represents an element of $\mM_{\hH,\Lambda}$ and has the following
properties:
\begin{enumerate}
\setlength{\itemsep}{0in}
\item $\tilde{u}$ is embedded.
\item $u$ is injective or $\p\Sigma = \emptyset$.
\item All asymptotic orbits $P_z$ for $z \in \Gamma$ are Morse-Bott, 
simply covered and geometrically distinct from one another.
\item $\ind(\tilde{u}) = 2$.
\item $\Sigma$ has genus~$0$.
\item For each $z \in \Gamma^\pm$, $\muCZ^{\Phi\mp}(P_z)$ is odd
(for any trivialization $\Phi$).
\end{enumerate}
Then $\tilde{u}$ is strongly regular and a neighborhood of $\tilde{u}$
in $\mM_{\hH,\Lambda}$ has naturally the structure of a smooth
$2$--dimensional manifold.  In fact, there exists an embedding
\begin{equation*}
\begin{split}
\RR\times (-1,1) \times\dot{\Sigma} 
  &\stackrel{\tilde{F}}\longrightarrow \RR\times M \\
(\sigma,\tau,z) &\longmapsto (a_\tau(z) + \sigma, u_\tau(z))
\end{split}
\end{equation*}
such that
\begin{enumerate}
\setlength{\itemsep}{0in}
\item[(i)] For $\sigma\in\RR$ and $\tau\in (-1,1)$, the embeddings
$\tilde{u}_{(\sigma,\tau)} = \tilde{F}(\sigma,\tau,\cdot) :
\dot{\Sigma} \to \RR\times M$ parametrize elements of $\mM_{\hH,\Lambda}$,
and $\tilde{u}_{(0,0)} = \tilde{u}$.
\item[(ii)] The map $F(\tau,z) = u_\tau(z)$ is an embedding
$(-1,1)\times\dot{\Sigma} \hookrightarrow M$, and its
image never intersects any of the orbits $P_z$ or Morse-Bott tori
$N_z$.  In particular the maps 
$u_\tau : \dot{\Sigma}\to M$ are embedded for each
$\tau\in(-1,1)$, with mutually disjoint images which
do not intersect their asymptotic limits.
\item[(iii)] If $P_z$ belongs to a Morse-Bott torus $N_z$, then
$u_\tau$ and $u_{\tau'}$ are asymptotic at $z$ to distinct orbits in 
$N_z$ whenever $\tau \ne \tau'$.
\item[(iv)] Every $\tilde{u}'$ sufficiently close to $\tilde{u}$ in
$\mM_{\hH,\Lambda}$ is parametrized by $\tilde{u}_{(\sigma,\tau)}$ for some
unique $\sigma\in\RR$, $\tau \in (-1,1)$.
\end{enumerate}
\end{thm}

From this and the smoothness of the nonlinear normal Cauchy-Riemann
operator defined in \cite{Wendl:BP1}, we obtain the following
deformation result.

\begin{thm}
\label{thm:parametrizedIFT}
Suppose $\tilde{u} \in \mM_{\hH,\Lambda}$ is as in 
Theorem~\ref{thm:IFT}, and
$$
\hH_{r} = (\xi_{r},X_{r},\omega_{r},J_{r}),
\qquad r \in (-1,1)
$$ 
is a smooth $1$--parameter family of stable Hamiltonian structures 
such that $\hH_0 = \hH$, and for each $r \in (-1,1)$, all of the orbits $P_z$ 
and Morse-Bott tori $N_z$ remain Morse-Bott orbits/tori with respect to 
$X_{r}$, while the surfaces $L_1,\ldots,L_m$ remain tangent to $X_r$.

Then there exists a number $\epsilon \in (0,1]$ and a smooth 
$1$--parameter family of maps
$$
\tilde{F}_r : \RR \times (-1,1) \times \dot{\Sigma} \to
\RR\times M,
\qquad
r \in (-\epsilon,\epsilon)
$$
such that $\tilde{F}_0(0,0,\cdot) = \tilde{u}$ and each map
$\tilde{F}_r$ has the properties of $\tilde{F}$ in Theorem~\ref{thm:IFT} 
with respect to the moduli space $\mM_{\hH_r,\Lambda}$.
\end{thm}

These local perturbation theorems start from the assumption that
$\tilde{u} = (a,u)$ is embedded, and usually also that $u$ is injective.
To study foliations, we need to know that such conditions are preserved under
convergence to limits.  One needs therefore a version of positivity of
intersections for holomorphic curves with boundary: the crucial observation
is that such a result holds whenever one can guarantee that all boundary
intersections are ``one-sided'' (cf.~\cite{Ye:filling}).
This is easy to show under the
conditions of Theorem~\ref{thm:IFT}, where the assumptions guarantee that
$\windpi(\tilde{u}) = 0$, implying $u$ is transverse to $X$ and
thus $\tilde{u}$ is transverse to the hypersurfaces $\RR\times L_j$.

\begin{thm}
\label{thm:injective}
Suppose $\hH_k = (\xi_k,X_k,\omega_k,J_k)$ is a sequence of stable Hamiltonian
structures converging in $C^\infty(M)$ to $\hH = (\xi,X,\omega,J)$,
such that $X_k$ is tangent to the
surfaces $L_1,\ldots,L_m$ for all $k$, and $\tilde{u}_k = (a_k,u_k) :
\dot{\Sigma} \to \RR\times M$ is a sequence of
curves in $\mM_{\hH_k,\Lambda}$ such that for all $k$, $u_k$ is embedded and
intersects each $L_j$ only at $\p\Sigma$.  Then if
$\tilde{u}_k$ converges in $\Cinftyloc$ to a somewhere injective curve
$\tilde{u} = (a,u) \in \mM_{\hH,\Lambda}$ with
$E_\omega(\tilde{u}) > 0$ and $\windpi(\tilde{u}) = 0$,
$u : \dot{\Sigma} \to M$ is embedded.

Moreover, suppose $\tilde{u}_k = (a_k,u_k)$, $\tilde{v}_k = (b_k,v_k) \in
\mM_{\hH_k,\Lambda}$ are sequences such that $u_k$ and $v_k$ both satisfy
the conditions above, and never intersect each other.
Then if $\tilde{u}_k \to \tilde{u} = (a,u)$ and $\tilde{v}_k \to 
\tilde{v} = (b,v)$ in $\Cinftyloc$ with $E_\omega(\tilde{u})$ and
$E_\omega(\tilde{v})$ both positive and $\windpi(\tilde{u}) =
\windpi(\tilde{v}) = 0$,
the images of $u$ and $v$ are either disjoint or identical.
\end{thm}

\section{Compactness}
\label{sec:compactness}

\subsection{The setup}
\label{subsec:setup}

For any pair of oriented knots $\gamma$ and $\gamma' \subset S^3$, denote
their linking number by $\lk(\gamma,\gamma') \in \ZZ$.  
Let $P_\infty \subset S^3$
be an oriented knot, and $K = K_1 \cup \ldots \cup K_m \subset S^3
\setminus P_\infty$ an oriented link whose components satisfy
$\lk(P_\infty,K_j) > 0$.  Each knot $K_j$ is the center of some solid torus
$N_j \subset S^3$; we assume that these solid tori are pairwise disjoint and
that $N := N_1 \cup \ldots \cup N_m \subset S^3$ is disjoint from $P_\infty$.
Denote $\p N_j = L_j$ and $M = S^3 \setminus
(\interior{N})$, so the oriented boundary of $M$ is $\p M = - \bigcup_j L_j$.
Let $\hH_k = (\xi_k,X_k,\omega_k,J_k)$ be a sequence of stable Hamiltonian
structures on $M$ which converge in $C^\infty(M)$ to a stable Hamiltonian 
structure $\hH_\infty = (\xi_\infty,X_\infty,\omega_\infty,J_\infty)$
and have the following properties for all $k \le \infty$:
\begin{enumerate}
\item $P_\infty$ is a nondegenerate periodic orbit of $X_k$.
\item Any other periodic orbit $P \subset M \setminus P_\infty$ of
  $X_k$ satisfies $\lk(P,P_\infty) \ne 0$.
\item $X_k$ is tangent to each torus $L_j$
\item There are trivializations $\Phi_k$ of $\xi_k|_M$
such that $\muCZ^{\Phi_k}(P_\infty) = 3$ and,
 if $\gamma \subset L_j$ is a positively oriented meridian,
 $\wind_\gamma^{\Phi_k}(TL_j \cap \xi_k) = 1$.
\end{enumerate}

The stable Hamiltonian structures $\hH_k$ define $\RR$--invariant
almost complex structures $\tilde{J}_k$ on $\RR\times M$, for which
the surfaces
$$
\tilde{L}_j := \{0\}\times L_j
$$
are totally real submanifolds.  We then consider a sequence of
embedded $\tilde{J}_k$--holomorphic curves 
$$
\tilde{u}_k = (a_k,u_k) : 
(\dot{\Sigma},j_k) \to (\RR\times M,\tilde{J}_k),
$$
where 
$$
\Sigma = S^2 \setminus \bigcup_{j=1}^m \dD_j
$$
for some set of pairwise disjoint open disks $\dD_j \subset \CC$,
$j_k$ is an arbitrary sequence of smooth complex structures on
$\Sigma$ and $\dot{\Sigma} = \Sigma\setminus\{\infty\}$.
Assume each $\tilde{u}_k$ has finite energy with respect to
$\hH_k$, is positively asymptotic at $\infty$ to the orbit $P_\infty$,
and satisfies the boundary condition
$$
\tilde{u}_k(\gamma_j) \subset \{\text{const}\} \times L_j
$$
for $\gamma_j = \p\dD_j$, where the constant in the $\RR$--factor is
arbitrary and independent for each component.  Thus
$\tilde{u}_k \in \mM_{\hH_k,\Lambda}$, where
$\Lambda = (\tilde{L}_1,\ldots,\tilde{L}_m)$.
We assume also that the maps $u_k : \dot{\Sigma} \to M$ are all injective
and have the following topological behavior at the boundary:
$$
\begin{minipage}{4.2in}
\textit{For each component $\gamma_j \subset \p\Sigma$,
the oriented loop $u_k(\gamma_j)$ is homotopic along $L_j$ to a
negatively oriented meridian, i.e.~$\lk(u_k(\gamma_j),P_\infty) = 0$
and $\lk(u_k(\gamma_j),K_j) = -1$.}
\end{minipage}
$$

\begin{remark}
\label{remark:repeatedKnots}
We've implicitly assumed thus far that $L_i$ and $L_j$ are disjoint when
$i \ne j$, but it's convenient also to allow the possibility that $K$ has
fewer than $m$ components.  As a notational convenience, we will continue
to list these components as $K_1 \cup \ldots \cup K_m$, with the understanding
that some of the $K_j$ (and the corresponding $N_j$ and $L_j$) may be
identical.  In this case different components of $\p\Sigma$ may satisfy
the same boundary condition.
\end{remark}

\begin{lemma}
\label{lemma:indexuk}
$\ind(\tilde{u}_k) = 2$.
\end{lemma}
\begin{proof}
Using the global trivialization $\Phi_k$, the boundary Maslov index
at $\gamma_j$ is twice the winding number of $TL_j \cap \xi_k$ along
a negatively oriented meridian; this gives $-2$.  We then have
$$
\mu(\tilde{u}_k) = \muCZ^{\Phi_k}(P_\infty) + \sum_{j=1}^m
\mu^{\Phi_k}_{u_k(\gamma_j)}(u^*\xi_k,TL_j \cap \xi_k)
= 3 - 2m,
$$
hence by \eqref{eqn:FredholmIndex},
$$
\ind(\tilde{u}_k) = \mu(\tilde{u}_k) - \chi(\dot{\Sigma}) + m = 
(3 - 2m) - (1 - m) + m = 2.
$$
\end{proof}

Applying Theorem~\ref{thm:IFT} to $\tilde{u}_k$ now yields:
\begin{cor}
\label{cor:uk}
For each $\tilde{u}_k = (a_k,u_k)$, $u_k : \dot{\Sigma} \to M$ is
embedded, transverse to $X_k$ and does not intersect $P_\infty$.  
In particular then, $u_k$ is transverse to $\p M$.
\end{cor}

Our main goal in the next few subsections will be to prove the following.

\begin{thm}[Compactness]
\label{thm:compactness}
There exists a sequence of numbers $c_k \in \RR$ and diffeomorphisms
$\varphi_k : \Sigma \to \Sigma$ which fix $\infty$ and preserve each
component of $\p\Sigma$, such that a subsequence of
$(a_k + c_k,u_k) \circ \varphi_k$ converges in $\Cinftyloc(\dot{\Sigma})$
to some $\tilde{u} = (a,u) \in \mM_{\hH_\infty,\Lambda}$, and the
continuous extensions of these maps over $\overline{\Sigma}$ converge
in~$C^0$.  Moreover, $\tilde{u}$ is embedded and $u$ is injective.
\end{thm}

A closely related result involves the degeneration of such a sequence
as $X_k$ is twisted to the point where the meridians
on $L_j$ become periodic orbits.
The following will be crucial in passing from holomorphic foliations
with boundary to stable finite energy foliations of Morse-Bott type:

\begin{thm}[Degeneration]
\label{thm:degeneration}
Assume $\tilde{u}_k$ and $\hH_k$ are as described above for all $k < \infty$,
but with the following change for $k=\infty$:
\begin{equation*}
\begin{minipage}{4.2in}
\textit{Any periodic orbit $P$ of $X_\infty$ in 
$(\interior{M})\setminus P_\infty$ satisfies $\lk(P,P_\infty) \ne 0$, and
each $L_j \subset \p M$ is a Morse-Bott torus with respect to $X_\infty$,
with orbits $P$ satisfying $\lk(P,P_\infty) = 0$ and $\lk(P,K_j) = -1$.}
\end{minipage}
\end{equation*}
Then there is a finite set $\Gamma' \subset \CC$ with $\#\Gamma' = m$,
a sequence of numbers $c_k \in \RR$ and diffeomorphisms
$\varphi_k : S^2 \setminus \Gamma' \to \interior{\Sigma}$ that fix $\infty$,
such that after passing to a subsequence, 
$\varphi_k^*j_k \to i$ in $\Cinftyloc(S^2\setminus \Gamma')$
and $(a_k + c_k, u_k) \circ \varphi_k$ converges in 
$\Cinftyloc(\CC\setminus \Gamma',\RR\times M)$ to a
$\tilde{J}_\infty$--holomorphic finite energy surface
$$
\tilde{u}_\infty : S^2 \setminus (\{\infty\} \cup \Gamma') \to \RR\times M.
$$
The map $\tilde{u}_\infty = (a_\infty,u_\infty)$ is embedded, with
$u_\infty : \CC\setminus\Gamma' \to M$ injective and disjoint from $\p M$,
all the punctures are positive, the asymptotic limit
at $\infty \in S^2$ is $P_\infty$, and for each component $\gamma_j \subset
\p\Sigma$, there is a corresponding puncture $z_j \in \Gamma'$ such that
the asymptotic limit at $z_j$ is a simply covered orbit on $L_j$.
Moreover for any sequence $\zeta_k \in \CC\setminus\Gamma'$ approaching
a puncture $z_j \in \Gamma'$, we have $u_k\circ\varphi_k(\zeta_k) \to L_j$.
\end{thm}

\subsection{Deligne-Mumford theory with boundary}
\label{subsec:DM}

\begin{figure}
\begin{center}
\includegraphics{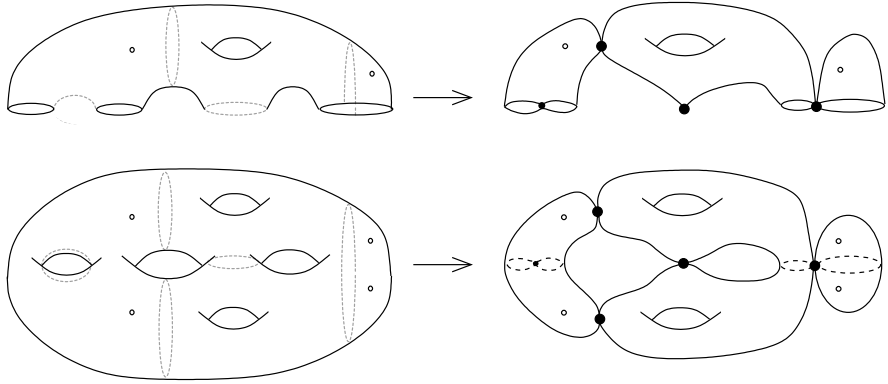}
\end{center}
\caption[Degeneration of a stable Riemann surface with
 boundary]{\label{fig:DM} Degeneration of a stable Riemann surface
$(\Sigma,j,\Gamma)$ with genus~$1$, four boundary components and two
interior marked points, together with its symmetric doubled surface.
The lightly shaded curves on the left are the geodesic loops and arcs that
shrink to zero length in the limit.  The right side shows the corresponding
singular surface $\widehat{\mathbf{\Sigma}}$ and its double after
degeneration; $\widehat{\mathbf{\Sigma}}$ has one interior double point,
two boundary double points and one unpaired node.}
\end{figure}

Before proving the main compactness results, we briefly review the
Deligne-Mumford compactification of the moduli space of Riemann surfaces
with boundary and interior marked points.  This space can always be
identified with a space of \emph{symmetric} surfaces without boundary,
``symmetric'' meaning there exists an antiholomorphic involution which
is respected by all the data (see Figure~\ref{fig:DM}).  A more detailed
discussion of this correspondence may be found in \cite{Wendl:thesis},
and for proofs of the compactness theorem we refer to
\cites{SeppalaSorvali,Hummel}.

Let $(\Sigma,j)$ be a compact connected 
Riemann surface, possibly with boundary, and let
$\Gamma \subset \interior{\Sigma}$ be a finite \emph{ordered} subset.  
As usual, denote the corresponding punctured surface by
$\dot{\Sigma} = \Sigma\setminus \Gamma$.
If the Euler
characteristic $\chi(\dot{\Sigma}) < 0$, then we call the
triple $(\Sigma,j,\Gamma)$ a \emph{stable Riemann surface with boundary
and interior marked points}.  The stability condition means
$$
2g + m + \#\Gamma > 2,
$$
where $g$ is the genus of $\Sigma$ and $m$ is the number of boundary
components.  Equivalently, one can say that $(\Sigma,j,\Gamma)$ is
stable if the Riemann surface obtained by doubling
$(\dot{\Sigma},j)$ along its boundary 
has negative Euler characteristic; this definition
has the advantage of also being correct when $\Sigma$ has marked points
on the boundary.

It is a standard fact that every Riemann surface without boundary and
with negative Euler characteristic admits a unique complete metric that
is compatible with the conformal structure and has constant curvature
$-1$: this is called the \emph{Poincar\'e metric}.  For a stable Riemann 
surface
$(\Sigma,j,\Gamma)$ with boundary, we define the Poincar\'e metric 
as the restriction to $\dot{\Sigma}$ of the Poincar\'e metric on the
symmetric double of $(\dot{\Sigma},j)$.  In this way, each component
of $\p\Sigma$ becomes a geodesic.

Denote by $\mM_{g,m,p}$ the moduli space of equivalence classes of
compact connected Riemann surfaces
$(\Sigma,j,\Gamma)$ with genus $g$, $m \ge 0$ boundary components and 
$p = \#\Gamma$ 
interior marked points $\Gamma\subset\interior{\Sigma}$.  Recall that
the points of $\Gamma$ come with an ordering.  Equivalence
$(\Sigma,j,\Gamma) \sim (\Sigma',j',\Gamma')$ means that 
there exists a biholomorphic
map $\varphi : (\Sigma,j) \to (\Sigma',j')$ that takes $\Gamma$ to
$\Gamma'$, preserving the ordering.
The topology on $\mM_{g,m,p}$ is defined by saying that
$[(\Sigma_k,j_k,\Gamma_k)] \to [(\Sigma,j,\Gamma)]$ if for sufficiently
large $k$ there exist diffeomorphisms $\varphi_k : \Sigma \to \Sigma_k$
mapping $\Gamma \to \Gamma_k$ (with the right ordering) and such that
$\varphi_k^* j_k \to j$ in $C^\infty$.

A \emph{nodal Riemann surface with boundary and interior marked points} 
$\mathbf{\Sigma} = (\Sigma, j, \Gamma, \Delta, N)$
consists of a Riemann surface $(\Sigma,j)$ with finitely many connected
components $\Sigma = \Sigma_1 \cup \ldots \cup \Sigma_q$, each of which
is a compact surface, possibly with boundary.  The marked point set
$\Gamma$ is a finite ordered set of interior points in $\Sigma$, 
and $\Delta$ is a set of unordered pairs of points in $\Sigma$,
$$
\Delta = \{ \{z_1,z_1'\}, \ldots, \{z_d,z_d'\} \},
$$
called \emph{double points}.  By assumption, the points
$z_1,z_1',\ldots,z_d,z_d'$ are all distinct,
and $z_j \in\p\Sigma$ if and only if $z_j'\in\p\Sigma$.  We will sometimes
abuse notation and regard $\Delta$ as a subset of $\Sigma$, rather than
a set of pairs.
There is also a finite unordered set
$N$ of interior points, which we'll call \emph{unpaired nodes}.
We assume the sets $\Gamma$, $\Delta$ and $N$ are all disjoint.
Intuitively, one thinks of $\mathbf{\Sigma}$ as the topological space obtained
from $\Sigma$ by identifying each pair of double points:
$$
\widehat{\mathbf{\Sigma}} = \Sigma / \{ z_j \sim z_j'
 \text{ for each pair $\{ z_j, z_j' \} \in \Delta$} \}.
$$
The point in $\widehat{\mathbf{\Sigma}}$ determined by a given pair of double
points $\{z_j,z_j'\}\in\Delta$ is called a \emph{node}.
We say that $\mathbf{\Sigma}$ is \emph{connected} whenever 
$\widehat{\mathbf{\Sigma}}$ is connected.  The connected components of 
$\mathbf{\Sigma}$ may be regarded as Riemann surfaces with boundary and
marked points
$(\Sigma_j, j|_{\Sigma_j}, (\Gamma \cup \Delta \cup N) \cap \Sigma_j)$,
which give rise to punctured surfaces
$$
\dot{\Sigma}_j =\Sigma_j \setminus ((\Gamma \cup \Delta \cup N) \cap \Sigma_j),
$$
having potentially both interior and boundary punctures.  We then say that
$\mathbf{\Sigma}$ is \emph{stable} if its connected components are all
stable; this means each $\dot{\Sigma}_j$ has negative Euler characteristic
after doubling.

The punctured components $\dot{\Sigma}_j$ can be compactified naturally
as follows: for an interior puncture $z$,
choose holomorphic coordinates identifying $z$ with $0$ in the standard
unit disk, use the map $(s,t) \mapsto e^{-2\pi(s+it)}$ to identify this
biholomorphically with the half-cylinder $[0,\infty) \times S^1$,
and compactify by adding the ``circle at infinity''
$\delta_z \cong \{\infty\} \times S^1$.
For punctures $z \in \p\Sigma$, one instead uses the map 
$(s,t) \mapsto e^{-\pi(s+it)}$ to identify a punctured upper half-disk
with the half strip $[0,\infty) \times [0,1]$, and then adds
the ``arc at infinity'' $\delta_z \cong \{\infty\} \times [0,1]$.  
Doing this for all punctures yields a compact surface
with piecewise smooth boundary.  Denote by $\overline{\Sigma}_j$ the
compactification of $\dot{\Sigma}_j \cup (\Gamma \cap \Sigma_j)$
obtained by adding such circles for each interior double point and
unpaired node, and arcs for each boundary double point in $\Sigma_j$.

Given a nodal surface
$(\Sigma,j,\Gamma,\Delta,N)$, an \emph{asymptotic marker}
at $z \in \Delta \cap \interior{\Sigma}$ is a choice 
of direction $\mu \in (T_z\Sigma \setminus \{0\}) /
\RR_+$, where $\RR_+$ is the group of positive real numbers, acting by
scalar multiplication.  A choice of asymptotic markers 
$r = \{ \{\mu_1,\mu_1'\}, \ldots, \{ \mu_d, \mu_d'\} \}$ 
for every pair of interior double points
is called a \emph{decoration}, and we then
call $(\mathbf{\Sigma},r) = (\Sigma,j,\Gamma,\Delta,N,r)$ 
a \emph{decorated nodal Riemann surface}.  For each pair 
$\{z,z'\} \in \Delta$ with asymptotic
markers $\{\mu,\mu'\}$, the conformal structure $j$ determines a natural
choice of orientation reversing map
$$
r_z : (T_{z}\Sigma \setminus \{0\}) / \RR_+ \to
  (T_{z'}\Sigma \setminus \{0\}) / \RR_+
$$
such that $r_z(\mu) = \mu'$, and hence also an orientation reversing
diffeomorphism $\bar{r}_z : \delta_z \to \delta_{z'}$.
For boundary pairs $\{z,z'\}$, the boundary determines natural
asymptotic markers and thus diffeomorphisms $\bar{r}_z$ between
the arcs $\delta_z$ and $\delta_{z'}$.  Then define
$$
\overline{\mathbf{\Sigma}}_r =
(\overline{\Sigma}_1 \sqcup \ldots \sqcup \overline{\Sigma}_q) /
\{ w \sim \bar{r}_z(w) \text{ for all $w\in \delta_{z}$, 
  $z \in \Delta \cup N$} \}.
$$
This is a compact surface with smooth boundary, 
and is connected if and only if 
$\widehat{\mathbf{\Sigma}}$ is connected.  In that case, we define the 
\emph{arithmetic signature} of $\mathbf{\Sigma}$ to be the pair
$(g,m)$ where $g$ is the genus of
$\overline{\mathbf{\Sigma}}_r$ (the \emph{arithmetic genus of
$\mathbf{\Sigma}$}) and $m$ is
the number of connected components of $\p\overline{\mathbf{\Sigma}}_r$.
We shall denote the union of the
special circles and arcs $\delta_z$ for $z \in \Delta \cup N$ by 
$\Theta_{\Delta,N}
\subset \overline{\mathbf{\Sigma}}_r$.  
The conformal structure $j$ on $\Sigma$
defines a singular conformal structure $j_{\mathbf{\Sigma}}$
on $\overline{\mathbf{\Sigma}}_r$,
which degenerates at $\Theta_{\Delta,N}$.
If $\mathbf{\Sigma}$ is stable, then there is similarly a ``singular
Poincar\'e metric'' $h_{\mathbf{\Sigma}}$ 
on $\overline{\mathbf{\Sigma}}_r \setminus \Gamma$, defined
by choosing the Poincar\'e metric on each of the punctured components
$\dot{\Sigma}_j$.  This metric
degenerates at $\Theta_{\Delta,N}$ as well as at $\Gamma$; 
in particular the distance from
a marked point $z_0\in \Gamma$ or a circle or arc 
$\delta_z \subset \Theta_{\Delta,N}$ to any other point
in $\overline{\mathbf{\Sigma}}_r$ is infinite, and the circles and arcs 
$\delta_z$ have length~$0$.  Observe that in the stable
case, $\chi(\overline{\mathbf{\Sigma}}_r \setminus\Gamma) < 0$,
i.e.~a stable nodal surface $(\Sigma,j,\Gamma,\Delta,N)$ with arithmetic
signature $(g,m)$ satisfies $2g + m + \#\Gamma > 2$.

\begin{figure}
\begin{center}
\includegraphics{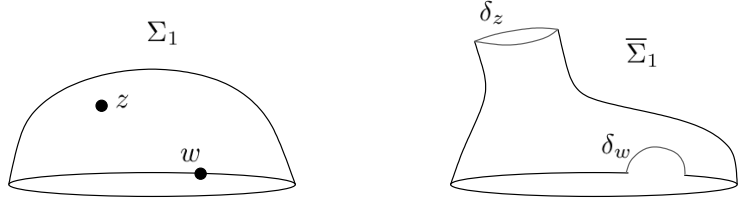}
\end{center}
\caption[Compactification with circles and arcs at infinity]
{\label{fig:DMblowup} A component $\Sigma_1$ with its compactification
$\overline{\Sigma}_1$.  Here there's one interior double point
$z \in \Delta \cap \interior{\Sigma}_1$ and one boundary double point
$w \in \Delta \cap \p\Sigma_1$.}
\end{figure}

Assume $2g + m + p > 2$ and let $\overline{\mM}_{g,m,p}$ denote the moduli 
space of equivalence classes of stable nodal Riemann surfaces
$\mathbf{\Sigma} = (\Sigma, j, \Gamma, \Delta,N)$ with arithmetic 
signature~$(g,m)$
and~$p = \#\Gamma$ interior marked points.  
We say $(\Sigma,j,\Gamma,\Delta,N) \sim
(\Sigma',j',\Gamma',\Delta',N')$ if there is a biholomorphic map
$\varphi : (\Sigma,j) \to (\Sigma',j')$ taking $\Gamma$ to $\Gamma'$
with the proper ordering, and such that $\varphi(N) = N'$ and 
$\{\varphi(z_1),\varphi(z_2)\} \in \Delta'$
if and only if $\{z_1,z_2\} \in \Delta$.
There is a natural inclusion $\mM_{g,m,p} \hookrightarrow 
\overline{\mM}_{g,m,p}$ defined by assigning to 
$[(\Sigma,j,\Gamma)]$ an empty set of double points and unpaired nodes.

The topology of $\overline{\mM}_{g,m,p}$ is determined by the following
notion of convergence.

\begin{defn}
\label{defn:DMconvergence}
A sequence $[\mathbf{\Sigma}_k] = 
[(\Sigma_k,j_k, \Gamma_k,\Delta_k,N_k)] \in \overline{\mM}_{g,m,p}$ 
converges to $[\mathbf{\Sigma}] = [(\Sigma,j, \Gamma,\Delta,N)] \in 
\overline{\mM}_{g,m,p} $ if there are decorations $r_k$ for
$\mathbf{\Sigma}_k$ and $r$ for $\mathbf{\Sigma}$, and diffeomorphisms
$\varphi_k : \overline{\mathbf{\Sigma}}_r \to 
(\overline{\mathbf{\Sigma}}_k)_{r_k}$, 
for sufficiently large $k$, with the following properties:
\begin{enumerate}
\item $\varphi_k$ sends $\Gamma$ to $\Gamma_k$, preserving the ordering.
\item $\varphi_k^* j_{\mathbf{\Sigma}_k} \to j_{\mathbf{\Sigma}}$ in 
 $\Cinftyloc(\overline{\mathbf{\Sigma}}_r 
 \setminus \Theta_{\Delta,N})$.
\item $\varphi_k^{-1}(\Theta_{\Delta_k,N_k}) \subset 
 \Theta_{\Delta,N}$, and all circles in
 $\varphi_k(\Theta_{\Delta,N}) \setminus \Theta_{\Delta_k,N_k}$
 are closed geodesics for the Poincar\'e metric $h_{\mathbf{\Sigma}_k}$ on
 $(\overline{\mathbf{\Sigma}}_k)_{r_k}$; similarly all arcs in
 $\varphi_k(\Theta_{\Delta,N}) \setminus \Theta_{\Delta_k,N_k}$
 are geodesic arcs with endpoints on 
$\p ( (\overline{\mathbf{\Sigma}}_k)_{r_k} )$.
\end{enumerate}
\end{defn}

\begin{thm}
\label{DM}
$\overline{\mM}_{g,m,p}$ is compact.  In particular, any sequence
of stable Riemann surfaces $(\Sigma_k,j_k,\Gamma_k)$ with boundary and
interior marked points, having fixed topological type and number of
marked points, has a subsequence convergent (in the sense of
Definition~\ref{defn:DMconvergence}) to a stable
nodal Riemann surface $(\Sigma,j,\Gamma,\Delta, N)$ with boundary and
interior marked points.
\end{thm}

\subsection{Preparation and removal of singularities}
\label{subsec:preparation}

In this section, fix a closed $3$--manifold $M$ with stable Hamiltonian
structure $\hH = (\xi,X,\omega,J)$ and an embedded surface $L \subset M$
tangent to $X$.  Let $\HH \subset \CC$ denote the closed upper half-plane,
and $\DD^+ := \DD \cap \HH$.
We now collect some lemmas that will be useful in the
compactness arguments to come.

\begin{lemma}[\cite{HoferZehnder}, Sec.~6.4, Lemma~5]
\label{lemma:Hofer}
Let $(X,d)$ be a complete metric space and $f : X\to [0,\infty)$ a continuous
function.  Then given any $x_0\in X$ and $\epsilon_0 > 0$, there exist
$x\in \overline{B_{2\epsilon_0}(x_0)}$ 
and $\epsilon \in (0,\epsilon_0]$ such that
\begin{equation*}
f(x)\epsilon \ge f(x_0)\epsilon_0 \quad \text{ and } \quad
f(y) \le 2f(x) \text{ for all $y \in \overline{B_\epsilon(x)}$}.
\end{equation*}
\end{lemma}

The next statement follows from the fact that when $\tilde{u} = (a,u)
: \dot{\Sigma} \to \RR\times M$ is a finite energy surface, the map $a$
is both proper and subharmonic (cf.~\cite{HWZ:props2}).
Note however that in non-contact cases, there can generally exist
nonconstant \emph{closed} $\tilde{J}$--holomorphic curves in $\RR\times M$.

\begin{lemma}
\label{lemma:positivePuncture}
Every finite energy surface with nonremovable punctures and no
boundary has at least one positive puncture.
\end{lemma}

\begin{lemma}
\label{lemma:zeroArea}
Suppose $(\Sigma,j)$ is a simply connected Riemann surface and
$\tilde{u} = (a,u) : (\Sigma,j) \to (\RR\times M,\tilde{J})$ is
pseudoholomorphic with $E_\omega(\tilde{u}) = 0$.  Then $\tilde{u}$ has
the form
$$
\tilde{u}(z) = (a(z),x(f(z)))
$$
where $x : \RR \to M$ is an orbit of $X$ (not necessarily periodic) and
$a + if : \Sigma \to \CC$ is a holomorphic function.
\end{lemma}
\begin{proof}
The integrand of $E_\omega(\tilde{u}) = \int_{\Sigma} u^*\omega$ is
everywhere nonnegative, and vanishes at $z \in \Sigma$ if and only if
the image of $du(z)$ is contained in $\RR X(u(z))$.  By assumption this
is true everywhere, thus $u(\Sigma)$ is contained in the image of some orbit
$x : \RR \to M$ of $X$.  If $x$ is not periodic, then it's a diffeomorphism
onto its image and can be inverted, allowing us to find a function
$f : \Sigma \to \RR$ such that $u = x \circ f$.  If $x$ is periodic, then
it can be viewed as a covering map onto its image, and the simple
connectedness of $\Sigma$ means that $u$ can be lifted to the cover,
producing again the map $f$.  An easy
computation now shows that $\tilde{u} = (a,u)$ is $\tilde{J}$--holomorphic
if and only if $a + if$ satisfies the standard Cauchy-Riemann equations.
\end{proof}

\begin{prop}
\label{prop:constants}
Suppose $\tilde{u} : \Sigma \to \RR\times M$ is $\tilde{J}$--holomorphic with
finite energy $E(\tilde{u}) < \infty$, $E_\omega(\tilde{u}) = 0$ and
$\tilde{u}(\p\Sigma) \subset \{0\}\times L$.
\begin{enumerate}
\setlength{\itemsep}{0in}
\item If $\Sigma$ is $S^2$, $\DD$, $\CC$ or $\HH$, then 
$\tilde{u}$ is constant.
\item If $\Sigma$ is $\RR\times S^1$ or $(-\infty,0]\times S^1$, there
exists a number $Q \in \RR$ and an orbit $x : \RR \to M$ of $X$ such that
up to $\RR$--translation,
$\tilde{u}(s,t) = (Qs, x(Qt))$.
\end{enumerate}
\end{prop}
\begin{proof}
Applying Lemma~\ref{lemma:zeroArea}, the function $a$ is harmonic,
and constant on the boundary.  For the compact cases $\Sigma = S^2$ and
$\DD$ this is enough: we conclude that $a$ is globally constant, and
so therefore is $f$, its harmonic conjugate.

For $\Sigma = \CC$, an argument from \cite{Hofer:weinstein}*{Lemma~28}
uses the finiteness of $E_\lambda(\tilde{u})$ to show that 
$\Phi := a + if$ is constant.  In brief, the energy can be rewritten as
$$
E_\lambda(\tilde{u}) = \sup_{\varphi\in\tT}
 \int_{\CC} \tilde{u}^*(d\varphi\wedge \lambda)
= \sup_{\varphi\in\tT} \int_{\CC} \Phi^* d(\varphi(s)\ dt),
$$
and a bubbling off argument shows that $\Phi$ must be constant if
the latter is finite.  We use this argument also for the case
$\Sigma = \HH$, after observing that the boundary condition implies
$\Phi(\RR) \subset i\RR$, so $\Phi$ extends by the Schwartz reflection
principle to an entire function on $\CC$.

The cylinder cases also follow from arguments in \cite{Hofer:weinstein}:
a simple bubbling off argument using the finite energy shows that 
$|d\tilde{u}|$ is globally bounded, and we then lift the domain to
$\CC$ or $\HH$ and apply Lemma~\ref{lemma:zeroArea}.  The holomorphic
function $\Phi = a + if$ need not be constant in this case 
but is affine due to the
gradient bound, so periodicity in $t$ implies $\Phi(s,t) = Qs + c + iQt$
for some constants $Q, c \in \RR$, and $x$ is $|Q|$--periodic unless $Q=0$.
\end{proof}

We will often use the fact that interior punctures of holomorphic curves with
finite energy are either asymptotic to periodic orbits or are removable;
the latter is the case whenever the image is contained within a compact
subset of $\RR\times M$.  This follows from the standard theorem on
removal of singularities (cf.~\cite{McDuffSalamon:Jhol}), together with
the following observation: if $\varphi : \RR \to (0,\epsilon)$ is a smooth
increasing function and $\epsilon$ is sufficiently small, then
\begin{equation}
\label{eqn:sympForm}
d(\varphi\lambda) + \omega
\end{equation}
is a symplectic form on $\RR\times M$, and any 
$\tilde{J}$--holomorphic map $\tilde{u}$ with
$E(\tilde{u}) < \infty$ also has finite symplectic area with respect to
this form.

We will need a corresponding statement for boundary punctures.  Since our
usual boundary condition on curves $\tilde{u} : \dot{\Sigma} \to \RR\times M$
constrains $\tilde{u}(\p\Sigma)$ to lie in a compact submanifold of
$\RR\times M$, one expects boundary punctures to removable.  One only needs
to show that such maps cannot become unbounded on the interior in a
neighborhood of a puncture; we will show that
this is always the case when the boundary condition has the form
$\tilde{L} = \{\text{const}\} \times L$.\footnote{As shown in 
\cite{Wendl:thesis}, a somewhat more general statement is true, but it's
not clear whether these arguments can be generalized to accomodate
arbitrary graphs $\tilde{L} = \{ (g(x),x) \in \RR\times L \}$.}

\begin{thm}[Removal of boundary singularities]
\label{thm:removeSing}
Suppose
$\tilde{u} = (a,u) : \dot{\DD}^+ = \DD^+\setminus\{0\} \to \RR\times M$
is a $\tilde{J}$--holomorphic map with 
$E(\tilde{u}) < \infty$ and $\tilde{u}(\dot{\DD}^+ \cap \RR) \subset
\{0\} \times L$.  Then $\tilde{u}$ extends to
a $\tilde{J}$--holomorphic half-disk $\DD^+ \to\RR\times M$ with
$\tilde{u}(\DD^+ \cap \RR) \subset \{0\} \times L$.
\end{thm}
\begin{proof}
By the remarks above, the result will follow from the standard removal of
singularities theorem after showing
that $\tilde{u}(\dot{\DD}^+)$ is contained in
a compact subset of $\RR\times M$.  To see that this is the case, 
compose $\tilde{u}$ with
the biholomorphic map $\psi : [0,\infty)\times [0,1] \to \dot{\DD}^+:
(s,t)\mapsto e^{-\pi (s+it)}$ and consider the pseudoholomorphic half-strip
$$
\tilde{v} = (b,v) = \tilde{u}\circ \psi : 
  [0,\infty)\times [0,1] \to \RR\times M.
$$
We claim that $|d \tilde{v}|$ is bounded on $[0,\infty)\times [0,1]$,
where the norm is defined with respect to the Euclidean metric on
$[0,\infty)\times [0,1] \subset \CC$ and any fixed $\RR$--invariant metric
on $\RR\times M$.  It will follow from this that $\tilde{v}$
(and hence $\tilde{u}$) is bounded, as
$\tilde{v}([0,\infty)\times \{0\})$ and
$\tilde{v}([0,\infty)\times \{1\})$ are contained in the compact set
$\{0\}\times L$.

If $|d \tilde{v}|$ is not bounded, there is a sequence $z_k =
s_k + i t_k \in [0,\infty) \times [0,1] \subset\CC$ such that
$R_k := |d\tilde{v}(z_k)| \to \infty$.  We may assume $s_k \to \infty$.
Choose a sequence of positive
numbers $\epsilon_k \to 0$ such that $\epsilon_k R_k \to \infty$; by
Lemma~\ref{lemma:Hofer} we can assume without loss of generality
that $|d \tilde{v}(z)| \le
2|d\tilde{v}(z_k)|$ whenever $|z - z_k| \le \epsilon_k$.  We will
define a sequence of rescaled maps which converge to either a plane or
a half-plane, depending on whether and how fast $z_k$ approaches the
boundary of $[0,\infty) \times [0,1]$.  We consider three cases.

\textit{Case~1: assume $t_k R_k$ and $(1 - t_k) R_k$ are both unbounded:}
then we can pass to a subsequence so that both approach $\infty$.
Let 
$$
r_k := \min \{ \epsilon_k R_k, t_k R_k, (1 - t_k) R_k \},
$$
so $r_k \to \infty$ and there are embeddings 
$$
\psi_k : \DD_{r_k} \hookrightarrow [0,\infty) \times [0,1] :
z \mapsto z_k + \frac{z}{R_k}.
$$
Use these to define rescaled maps
$$
\tilde{v}_k = (b \circ \psi_k - b(z_k), v \circ \psi_k) :
\DD_{r_k} \to \RR\times M,
$$
which satisfy a uniform $C^1$--bound and have a subsequence convergent 
in $\Cinftyloc$ to a non-constant $\tilde{J}$--holomorphic plane
$\tilde{v}_\infty : \CC\to \RR\times M$.  This map has finite energy
$E(\tilde{v}_\infty) \le E(\tilde{v}) < \infty$, but also
$E_\omega(\tilde{v}_\infty) = 0$, giving a contradiction to
Prop~\ref{prop:constants}.

\textit{Case~2: assume $t_k R_k$ is bounded.}  This means $z_k$ is approaching
the half-line $[0,\infty) \times \{0\}$.  Let
$$
\psi_k : \DD_{\epsilon_k R_k}^+ \hookrightarrow [0,\infty) \times [0,1] :
z \mapsto s_k + \frac{z}{R_k},
$$
and define a sequence of rescaled maps $\tilde{v}_k : \DD_{\epsilon_k R_k}^+
\to \RR\times M$ by $\tilde{v}_k = \tilde{v} \circ \psi_k$.  These maps
satisfy the boundary condition $\tilde{v}_k(\DD_{\epsilon_k R_k} \cap \RR)
\subset \{0\}\times L$.  Moreover, the points $\tilde{v}_k(0) = \tilde{v}(s_k)$
are contained in the compact set $\{0\}\times L$, 
and there is a uniform gradient
bound $|d\tilde{v}_k(z)| \le 2$ for all $z \in \DD_{\epsilon_k R_k}^+$.  Thus
a subsequence converges in $\Cinftyloc(\HH,\RR\times M)$ to a
$\tilde{J}$--holomorphic half-plane $\tilde{v}_\infty : \HH \to \RR\times M$
with $E(\tilde{v}_\infty) < \infty$ and $E_\omega(\tilde{v}_\infty) = 0$.
We claim however that $\tilde{v}_\infty$ is not constant.  Indeed, 
$| \psi_k^{-1}(z_k) | = R_k | z_k - s_k | = t_k R_k$ is bounded,
thus passing to a subsequence, $\psi_k^{-1}(z_k) \to \zeta \in \HH$,
and $|d \tilde{v}_\infty(\zeta)| = \lim_k 
 |d \tilde{v}_k(\psi_k^{-1}(z_k))| = 1$.  Thus the existence of
$\tilde{v}_\infty$ again contradicts Prop~\ref{prop:constants}.

\textit{Case~3: assume $(1 - t_k) R_k$ is bounded.}  This is very similar
to the previous case; this time $z_k$ is approaching the half-line
$[0,\infty)\times\{1\}$, so we rescale using the embeddings
$$
\psi_k : \DD_{\epsilon_k R_k}^+ \hookrightarrow [0,\infty) \to [0,1] :
z \mapsto (s_k + i) - \frac{z}{R_k}.
$$
Then by the same arguments used above, $\tilde{v}_k = \tilde{v}\circ \psi_k$ 
has a subsequence convergent to a non-constant finite energy half-plane
$\tilde{v}_\infty : \HH\to\RR\times M$ with boundary condition
$\tilde{v}_\infty(\RR) \subset \{0\}\times L$ and 
$E_\omega(\tilde{v}_\infty) = 0$, giving another contradiction.
\end{proof}

\subsection{Taming forms and energy bounds}
\label{subsec:energy}

We now proceed toward the proofs of Theorems~\ref{thm:compactness}
and~\ref{thm:degeneration}, so let $M$, $K$, $L_j$, $\hH_k$
and $\tilde{u}_k$ be as defined in \S\ref{subsec:setup}.  Observe that
since $X_k$ is always tangent to the tori $L_j$, the $2$--forms
$d\lambda_k$ and $\omega_k$ arising from $\hH_k$ vanish on $L_j$.

\begin{lemma}
\label{lemma:exactTaming}
The taming forms $\omega_k$ for $k \le \infty$ are exact on $M$.
\end{lemma}
\begin{proof}
Applying a Mayer-Vietoris sequence to $S^3 = M \cup N$, $H_2(M)$ is generated
by the inclusions of the fundamental classes $[L_j] \in H_2(L_j)$ for 
$j = 1,\ldots,m$.  Then $\int_{L_j} \omega_k = 0$ implies that $\omega_k$
vanishes on $H_2(M; \RR)$, hence $[\omega_k] = 0 \in H^2(M; \RR)$.
\end{proof}

Denote by $E_k(\tilde{u}_k)$ the energy of
$\tilde{u}_k$, computed with respect to $\omega_k$ and $\lambda_k$.
Surfaces of the form $\{\text{const}\}\times L_j$ are not only totally
real in $(\RR\times M,\tilde{J})$ but also \emph{Lagrangian} with
respect to natural symplectic forms as in \eqref{eqn:sympForm}.
This will permit a uniform energy bound for the sequence $\tilde{u}_k$.

\begin{lemma}
\label{lemma:contactArea}
$\int_{\dot{\Sigma}} u_k^*\omega_k$ is uniformly bounded for all $k$. 
\end{lemma}
\begin{proof}
Since $\omega_k$ is exact and vanishes on each torus $L_j$,
$\int_{\dot{\Sigma}} u_k^*\omega_k$ depends only on the asymptotic behavior
and the homotopy class of $u_k|_{\p\Sigma} : \p\Sigma \to \p M$.
Thus we can pick any smooth map $u : \dot{\Sigma} \to M$ with the proper
behavior and write
$\int_{\dot{\Sigma}} u_k^*\omega_k = \int_{\dot{\Sigma}} u^*\omega_k
\to \int_{\dot{\Sigma}} u^*\omega_\infty$.
\end{proof}

\begin{lemma}
\label{lemma:omegaMetric}
There exists a constant $C > 0$ independent of $k$ such that
$| u^*d\lambda_k | \le C \cdot u^*\omega_k$ for every 
$\tilde{J}_k$--holomorphic curve $\tilde{u} = (a,u)$ in $\RR\times M$.
\end{lemma}
\begin{proof}
Let $C_k > 0$ be the $C^0$--norm of the bilinear form
$$
\xi_k \times \xi_k \to \RR : (v,w) \mapsto d\lambda_k(v,J_k w)
$$
on the bundle $\xi_k \to M$ with respect to the bundle metric
$|v|_k^2 := \omega_k(v,J_k v)$; this is finite since $M$ is compact.
Since $\omega_k$ and $\lambda_k$ each converge in $C^\infty$ as $k\to \infty$,
the sequence $C_k$ is also \emph{bounded}, $C_k \le C$.  Denote by
$\pi_k : TM \to \xi_k$ the projection along $X_k$, and note that both
$d\lambda_k$ and $\omega_k$ annihilate $X_k$.  Then in any local
holomorphic coordinate system $(s,t)$ on the domain,
\begin{multline*}
| u^* d\lambda_k(\p_s,\p_t) | = | d\lambda_k(\pi_k u_s, \pi_k u_t) |
= | d\lambda_k(\pi_k u_s, J_k \pi_k u_s) | \le C_k | \pi_k u_s |_k^2 \\
 \le C \cdot \omega_k(\pi_k u_s, J_k \pi_k u_s)
 = C \cdot \omega_k(\pi_k u_s, \pi_k u_t)
 = C \cdot u_k^*\omega_k(\p_s,\p_t).
\end{multline*}
\end{proof}

\begin{prop}
\label{prop:energyBound}
There exists a constant $C > 0$ such that $E_k(\tilde{u}_k) < C$.
\end{prop}
\begin{proof}
Writing $E_k(\tilde{u}_k) = E_{\omega_k}(\tilde{u}_k) + 
E_{\lambda_k}(\tilde{u}_k)$, the first term is uniformly bounded due to
Lemma~\ref{lemma:contactArea}.  The second is
$\sup_{\varphi \in \tT} E_{\lambda_k}^\varphi(\tilde{u}_k)$, where
$$
E_{\lambda_k}^\varphi(\tilde{u}_k) := 
\int_{\dot{\Sigma}} \tilde{u}_k^* (d\varphi \wedge \lambda_k).
$$
Writing $d(\varphi \lambda_k) = 
\varphi\ d\lambda_k + d\varphi \wedge \lambda_k$, then applying
Stokes' theorem, Lemma~\ref{lemma:omegaMetric} and the fact that
$a_k$ is locally constant at $\p\Sigma$, we have
\begin{equation*}
\begin{split}
\left| E_{\lambda_k}^\varphi(\tilde{u}_k) \right| 
 &\le \left| \int_{\dot{\Sigma}} 
 \tilde{u}_k^*(\varphi\ d\lambda_k) \right| +
 \left| \int_{\dot{\Sigma}} \tilde{u}_k^*d(\varphi\lambda_k) \right| \\
&\le \int_{\dot{\Sigma}} | \varphi\circ a_k| \cdot |u_k^* d\lambda_k |
 + \int_{P_\infty} \lambda_k + 
 \left| \int_{\p\Sigma} \tilde{u}_k^*(\varphi\lambda_k) \right| \\
&\le C \int_{\dot{\Sigma}} u_k^*\omega_k + T_k + 
 \left| \int_{\p\Sigma} u_k^*\lambda_k \right|,
\end{split}
\end{equation*}
where $T_k$ is the period of $P_\infty$ as an orbit of $X_k$.  This is clearly
bounded as $k \to \infty$, and so is the first term, by 
Lemma~\ref{lemma:contactArea}; it remains only to bound
$\int_{\p\Sigma} u_k^*\lambda_k$.  Here we use the fact that
$d\lambda_k(X_k,\cdot) \equiv 0$, hence $d\lambda_k$ annihilates each
$L_j$, and this integral therefore only depends on $\lambda_k$ and
$[u_k|_{\p\Sigma}] = \ell \in H_1(L_1 \cup \ldots \cup L_m)$.  In particular,
it approaches $\int_{\ell} \lambda_\infty$ as $k \to \infty$.
\end{proof}

\subsection{Bubbling}
\label{subsec:gradBounds}

In this section we establish uniform bounds on the first derivatives of
the maps $\tilde{u}_k$.  In the \emph{non-stable} ($m < 2$) case,
the arguments of this section suffice to prove $\Cinftyloc$--convergence.
For $m \ge 2$, we also need to ensure that the sequence of conformal
structures induced by $j_k$ is compact; this issue will be dealt with in
\S\ref{subsec:conformal}.
The fundamental argument is that any gradient blow up causes
the bubbling off of a holomorphic plane or disk, which for topological
reasons cannot exist.

We focus first on the stable case, thus assume $m \ge 2$,
$\chi(\dot{\Sigma}) < 0$.  Then each of the stable Riemann surfaces
$(\Sigma,j_k,\{\infty\})$ determines a Poincar\'e metric $h_k$, which is
the restriction of a complete metric $h_k^D$ of constant curvature~$-1$
on the Riemann surface 
$$
(\dot{\Sigma}^D,j^D_k),
$$
obtained by doubling $(\dot{\Sigma},j_k)$ along the
boundary.  Denote the injectivity radius of $h_k^D$ at any point
$z \in \dot{\Sigma}$ by $\inj_k(z)$.

Fix any metric on $M$ and extend it to an $\RR$--invariant metric on
$\RR\times M$ in the natural way.  In the following, we will always use
the Euclidean metric on subsets of $\CC$ or $\RR\times S^1$, and
one of the Poincar\'e metrics $h_k$ on $\dot{\Sigma}$, with dependence on
$k$ reflected in the notation.  So for instance, given
$\tilde{u} : \dot{\Sigma} \to \RR\times M$, $| d\tilde{u}(z) |_k$ is the
norm of the linear map $d\tilde{u}(z) : T_z\dot{\Sigma} \to 
T_{\tilde{u}(z)}(\RR\times M)$ with respect to $h_k$ and the fixed metric
on $\RR\times M$.  For $\varphi : \DD \to \dot{\Sigma}$, define
$|d\varphi(z)|_k$ with respect to the Euclidean metric on the
domain and $h_k$ on the target.

The following technical lemma provides good coordinates in a neighborhood of
any point $z_0 \in \dot{\Sigma}$.  They are constructed by lifting 
$\dot{\Sigma}$ to the hyperbolic open disk so that~$0$ covers $z_0$,
then projecting the embedding $\DD \hookrightarrow \interior{\DD} : 
z \mapsto rz$
for sufficiently small $r > 0$ down to $\dot{\Sigma}$.  It follows from
the hyperbolic geometry of the disk that the resulting embedding has
the desired properties; see \cite{Wendl:thesis} for details.  Denote by
$\DD_\rho$ the standard closed disk of radius $\rho > 0$ in $\CC$.

\begin{lemma}
\label{lemma:coords}
Let $(\dot{\Sigma},j)$ be a stable punctured Riemann surface without boundary,
with Poincar\'e metric $h$, whose injectivity radius at $z \in \dot{\Sigma}$
we denote by $\inj(z)$.  There are positive constants 
$c_i$ and $C_i$ depending only on the topological type of 
$\dot{\Sigma}$ (i.e.~not on $j$), such that the following is true: 
for any $z_0 \in \dot{\Sigma}$ and any geodesic $\gamma$ passing
through $z_0$,
there is a holomorphic embedding $\varphi : \DD \hookrightarrow \dot{\Sigma}$ 
such that $\varphi(0) = z_0$, $\varphi$ maps $\RR\cap \DD$ to $\gamma$
preserving orientation, and
\begin{equation}
\label{eqn:chartGrad}
c_1 \cdot \inj(z_0) \le | d\varphi(z) |_h \le C_1 \cdot \inj(z_0) \quad
 \text{ for all $z\in \DD$}.
\end{equation}
For any $\rho\in[0,1]$, the image $\varphi(\DD_\rho)$
is then a closed ball of radius $d(\rho)$ in $(\dot{\Sigma},h)$, where
\begin{equation}
\label{eqn:bigBall}
c_2 \rho \cdot \inj(z_0) \le d(\rho) \le C_2 \rho \cdot \inj(z_0),
\end{equation}
and the injectivity radius at any point $\varphi(w)$ for $w \in \DD$ with
$|w| = \rho$ can be estimated by
\begin{equation}
\label{eqn:injradBound}
(c_3 - c_4 \rho) \cdot \inj(z_0) \le \inj(\varphi(w)) 
  \le (1 + C_3\rho) \cdot \inj(z_0)
\end{equation}
\end{lemma}

\begin{remark}
\label{remark:halfDisk}
Lemma~\ref{lemma:coords} extends to surfaces $\dot{\Sigma}$ with 
nonempty boundary as
follows: for any $z_0 \in \p\Sigma$, the component
$\gamma \subset \p\Sigma$ containing $z_0$ is a closed geodesic in the doubled
surface $(\dot{\Sigma}^D,h^D)$.  Thus the lemma gives an
embedding $\varphi : \DD^+ \to \dot{\Sigma}$ of the upper half disk,
sending $0$ to $z_0$ and $\RR\cap \DD^+$ into $\p\Sigma$.
\end{remark}

The first main objective in this section is the following result, which
says in effect that there is no bubbling off in the sequence $\tilde{u}_k$.

\begin{prop}
\label{prop:gradBounds}
If $\chi(\dot{\Sigma}) < 0$, then there is a constant $C > 0$ such that
\begin{equation}
\label{eqn:gradBound}
| d \tilde{u}_k(z) |_k \le \frac{C}{\inj_{k}(z)}
\end{equation}
for all $z \in \dot{\Sigma}$ and all $k$.
\end{prop}
\begin{proof}
Assume there exists a sequence $z_k \in \dot{\Sigma}$ such that
$\inj_{k}(z_k) \cdot | d \tilde{u}(z_k) |_k \to \infty$.
Using Lemma~\ref{lemma:coords}, choose a sequence of holomorphic 
embeddings
$$
\varphi_k : \DD \hookrightarrow \dot{\Sigma}^D
$$
such that $|d\varphi_k|_k$, the radii of the disks $\varphi_k(\DD)$
and the injectivity radius satisfy the bounds specified in the lemma.
Let 
$$
\rho_k = \min \big\{ |\zeta|\ \big|\ \zeta \in \varphi_k^{-1}(\p\Sigma) \big\},
$$
or $\rho_k = \infty$ if $\varphi_k(\DD) \cap \p\Sigma = \emptyset$.
The sequence $\rho_k$ determines whether or not we can restrict the embeddings
$\varphi_k$ in a uniform way so that their images are in $\dot{\Sigma}$.

\textit{Case~1: assume there is a number $\rho \in (0,1]$ and a subsequence
for which $\rho_k \ge \rho$.}  Then the restrictions of $\varphi_k$ to 
$\DD_\rho$
are embeddings into $\dot{\Sigma}$, and we can define a sequence of
pseudoholomorphic disks
$$
\tilde{v}_k = (b_k,v_k) = \tilde{u}_k \circ \varphi_k : 
 \DD_\rho \to \RR\times S^3,
$$
which satisfy a uniform energy bound
$$
E_k(\tilde{v}_k) \le E_k(\tilde{u}_k) \le C.
$$
Denoting the Euclidean metric on $\DD$ by $\eta$, the fact that
$\varphi_k : (\DD_\rho,\eta) \to (\dot{\Sigma},h_k)$ is conformal implies that
the norms of $d\varphi_k(z)$ and its inverse are reciprocals.  Then a
simple computation shows
$$
| d\tilde{v}_k(0) | = | d\tilde{u}_k(z_k) |_k \cdot |d\varphi_k(0)|_k
\ge c_1 |d\tilde{u}_k(z_k) |_k \cdot \inj_{k}(z_k) \to \infty.
$$
By Lemma~\ref{lemma:Hofer}, we can choose a sequence $\zeta_k \in 
\DD_\rho$
and positive numbers $\epsilon_k \to 0$ such that 
$R_k := |d\tilde{v}(\zeta_k)| \to \infty$,
$\epsilon_k R_k \to \infty$ and $|d\tilde{v}(\zeta)| \le 2R_k$ for
all $\zeta \in 
\DD_\rho$ with $|\zeta - \zeta_k| \le \epsilon_k$.  Assume without
loss of generality that $\overline{B_{\epsilon_k}(\zeta_k)} \subset \DD_\rho$
and define 
$$
\psi_k : \DD_{\epsilon_k R_k} \to
\overline{B_{\epsilon_k}(\zeta_k)} : \zeta \to \zeta_k + \frac{\zeta}{R_k}.
$$
Then we can define a rescaled sequence of $\tilde{J}_k$--holomorphic maps 
$\tilde{w}_k = (\beta_k, w_k) : \DD_{\epsilon_k R_k} \to \RR\times S^3$ by
$$
(\beta_k(\zeta) , w_k(\zeta)) = (b_k \circ \psi_k(\zeta) - b_k(\zeta_k), 
  v_k \circ \psi_k(\zeta)).
$$
These satisfy the uniform gradient bound 
$|d \tilde{w}_k(\zeta)| \le 2$, and they all map $0$ into a compact
subset of $\RR\times S^3$, thus a subsequence converges in 
$\Cinftyloc$ to a $\tilde{J}_\infty$--holomorphic plane
$$
\tilde{w}_\infty = (\beta_\infty, w_\infty) : \CC \to \RR\times S^3.
$$
The bound on $E_k(\tilde{v}_k)$ gives also a bound on
$E_k(\tilde{w}_k)$ and thus implies $E_\infty(\tilde{w}_\infty) < \infty$,
so $\tilde{w}_\infty$ is a non-constant finite energy plane.  

If the puncture at $\infty$ is removable, $\tilde{w}_\infty$ extends to a
nonconstant $\tilde{J}_\infty$--holomorphic sphere in $\RR\times M$.
Recall now that $\xi_\infty$ admits a global trivialization over $M$, thus
$c_1(w_\infty^*\xi_\infty) = 0$ and we have
$$
\ind(\tilde{w}_\infty) = \mu(\tilde{w}_\infty) - \chi(S^2) = 
2 c_1(w_\infty^*\xi_\infty) - 2 = -2.
$$
Since $\tilde{w}_\infty$ is nonconstant, it cannot have 
$E_{\omega_\infty}(\tilde{w}_\infty) = 0$, by Prop.~\ref{prop:constants}.
Thus $E_{\omega_\infty}(\tilde{w}_\infty) > 0$ and Theorem~\ref{thm:windpi}
gives a contradiction, in the form
$$
0 \le 2\windpi(\tilde{w}_\infty) \le \ind(\tilde{w}) - 2 = -4.
$$

It remains to find a contradiction in the case of a non-removable puncture
at~$\infty$: then $\tilde{w}_\infty$ is asymptotic
to some periodic orbit $P$ of $X_{\infty}$. 
We now use a topological argument to show that this is impossible.
 
If $P$ is geometrically distinct from $P_\infty$, then $\lk(P,P_\infty) \ne 0$
by assumption.  For some large radius $R$, the image $w_\infty(\p \DD_R)$ is
uniformly close to $P$, and we may assume the same is true of 
$P' := w_k(\p \DD_R)$ for sufficiently large $k$, thus
$\lk(P',P_\infty) \ne 0$.  But since $w_k$ is a reparametrization of
$u_k : \dot{\Sigma} \to S^3$ over some disk, this means there is a disk
$\dD \subset \dot{\Sigma}$ such that $P' = u_k(\p\dD)$.  The linking condition
then implies that $u_k(\dD)$ intersects $P_\infty$, contradicting the result
of Corollary~\ref{cor:uk}.

Suppose now that $P$ is identical to $P_\infty$ or some cover thereof.
For any component $K_j \subset K$,
observe that $u_k(\dot{\Sigma})$ never intersects
$K_j$.  Then repeating the argument above, we find a disk $\dD\subset
\dot{\Sigma}$ such that for sufficiently large $k$, $u_k(\p\dD)$ is a
knot uniformly close to $P_\infty$.  This implies $\lk(P_\infty,K_j) = 0$,
another contradiction, thus proving that the
plane $\tilde{w}_\infty$ cannot exist.

\textit{Case~2: assume $\rho_k \to 0$.}
Here we will find that either a plane or a disk bubbles off, depending
on how fast $\rho_k$ approaches $0$.  Choose a sequence $\zeta_k'\in \DD$
such that $\zeta_k := \varphi_k(\zeta_k') \in \p\Sigma$ and 
$|\zeta_k'| = \rho_k$.
By Remark~\ref{remark:halfDisk}, we can find a sequence of holomorphic
embeddings
$$
\varphi^+_k : \DD^+ \hookrightarrow \dot{\Sigma}
$$
that map $0$ to $\zeta_k$ and $\DD^+\cap \RR$ into $\p\Sigma$, 
and satisfy the bounds specified in
Lemma~\ref{lemma:coords}.  We claim there is a sequence of radii
$r_k \to 0$ such that $z_k \in \varphi^+_k(\DD_{r_k}^+)$.
Indeed, from Lemma~\ref{lemma:coords}, we know that
$\varphi^+(\DD_{r_k}^+)$ contains all points $\zeta \in \dot{\Sigma}$ with
$\dist_{h_k}(\zeta,\zeta_k) \le d_k$, where
$$
d_k \ge c_2 r_k \cdot \inj_{k}(\zeta_k).
$$
We have also the estimates
\begin{equation*}
\begin{split}
\dist_{h_k}(z_k,\zeta_k) &\le C_2 \rho_k \cdot \inj_{k}(z_k),\\
\inj_{k}(\zeta_k) &\ge (c_3 - c_4 \rho_k) \cdot \inj_{k}(z_k).
\end{split}
\end{equation*}
Then when $\rho_k$ is sufficiently small we can set
$$
r_k = \frac{2 C_2}{c_2 (c_3 - c_4 \rho_k)} \rho_k
$$
and compute,
$$
\dist_{h_k}(\zeta_k,z_k) \le 
\frac{C_2}{c_3 - c_4 \rho_k} \rho_k \cdot \inj_{k}(\zeta_k)
= \frac{1}{2} c_2 r_k \inj_{k}(\zeta_k) < d_k.
$$
We can thus choose a sequence $z_k' \in \DD^+$ with $z_k' \to 0$ and
$\varphi_k^+(z_k') = z_k$.  Defining a sequence of $\tilde{J}_k$--holomorphic
half-disks 
$$
\tilde{v}_k = \tilde{u}_k \circ \varphi^+ : \DD^+ \to \RR\times S^3,
$$
we have
\begin{multline*}
R_k := 
| d\tilde{v}_k(z_k') | = | d\tilde{u}_k(z_k) |_k \cdot 
 |d\varphi^+(z_k')|_k
 \ge C | d\tilde{u}_k(z_k) |_k \cdot \inj_{k}(\zeta_k)\\
 \ge C (c_3 - c_4 \rho_k) | d\tilde{u}_k(z_k) |_k \cdot \inj_{k}(z_k)
 \to \infty.
\end{multline*}
Using Lemma~\ref{lemma:Hofer}, we may assume there is a sequence of
positive numbers $\epsilon_k \to 0$ such that $\epsilon_k R_k \to \infty$
and $| d\tilde{v}_k(z) | \le 2 R_k$ for all $z\in \DD^+$ with
$|z - z_k'| \le \epsilon_k$.  Writing $z_k' = s_k + it_k$, there are
two possibilities:

\textit{Case~2a: assume $t_k R_k$ is unbounded.}  Passing to a subsequence,
we may assume $t_k R_k \to \infty$, thus $r_k' := \min\{ \epsilon_k R_k,
t_k R_k\} \to \infty$.  Then for sufficiently large $k$ we can
define embeddings $\psi_k : \DD_{r_k'} \hookrightarrow \DD^+$ by
$$
\psi_k(z) = z_k + \frac{z}{R_k}.
$$
Arguing as in case~1, there is now a sequence of rescaled maps 
$$
\tilde{w}_k = (\beta_k,w_k) = \tilde{v}_k \circ
\psi_k : \DD_{r_k'} \to \RR\times S^3
$$
and constants $c_k\in\RR$ such that a subsequence of $(\beta_k + c_k, w_k)$
converges in $\Cinftyloc$ to a nonconstant
$\tilde{J}_\infty$--holomorphic finite energy plane 
$\tilde{w}_\infty = (\beta_\infty,w_\infty) : \CC\to \RR\times S^3$,
giving the same contradiction as before.

\textit{Case~2b: assume $t_k R_k$ is bounded.}  
Now define $\psi_k : \DD_{\epsilon_k R_k}^+ \hookrightarrow \DD^+$ by
$$
\psi_k(z) = s_k + \frac{z}{R_k},
$$
and let 
$$
\tilde{w}_k = (\beta_k,w_k) = 
\tilde{v}_k \circ \psi_k : \DD_{\epsilon_k R_k}^+ \to
\RR\times S^3.
$$
Then $|d\tilde{w}_k|$ is uniformly bounded.  Passing to a subsequence,
$(\beta_k - \beta_k(0), w_k)$ converges in $\Cinftyloc$ to a
$\tilde{J}_\infty$--holomorphic half-plane
$$
\tilde{w}_\infty = (\beta_\infty, w_\infty) : \HH \to \RR\times S^3,
$$
satisfying $E_\infty(\tilde{w}_\infty) < \infty$ and the boundary condition
$\tilde{w}_\infty(\RR) \subset \{0\} \times L_j$ for some 
$j \in \{1,\ldots,m\}$.  It is not constant, since
$|d \tilde{w}_k( it_k R_k)| = \frac{1}{R_k} 
|d \tilde{v}_k(s_k + it_k)| = 1$ and a subsequence of $i t_k R_k$
converges in $\DD^+$.  Now identifying $\HH$ conformally with
$\DD\setminus \{1\}$, we can regard $\tilde{w}_\infty$ as
a holomorphic disk with a puncture on the boundary, and
Theorem~\ref{thm:removeSing} tells us that the puncture is removable.
Thus extending over the puncture defines
a $\tilde{J}_\infty$--holomorphic disk
$$
\tilde{w} = (\beta, w) : \DD \to \RR\times S^3
$$
with $w(\p \DD) \subset L_j$.  By topological considerations, we can
severely restrict the homotopy class of the loop
$\gamma = w|_{\p \DD} : \p \DD \to L_j$.
Indeed, choose a radius $r$ slightly less than $1$ so that
$w|_{\p \DD_r} : \p \DD_r \to S^3$ is uniformly 
close to $\gamma$.  Returning
to the half-plane $\HH$, there is then a large simply connected region
$\Omega \subset \HH$ with smooth boundary such that for large $k$,
$w_k|_{\p \Omega} : \p\Omega \to S^3$ is also uniformly
close to $\gamma$.  Undoing the reparametrization one step further, there
is then an embedded disk $\dD \subset \dot{\Sigma}$ such that for some
large $k$,
$$
u_k|_{\p\dD} : \p\dD \to S^3
$$
is uniformly close to $\gamma$.  Since $u_k$ does not intersect either
$P_\infty$ or any of the knots $K_j\subset K$, this implies
$$
\lk(\gamma,P_\infty) = \lk(\gamma,K_1) = \ldots = \lk(\gamma,K_m) = 0.
$$
This is only possible if $\gamma$ is contractible on $L_j$.  But this
implies that the Maslov index $\mu(\tilde{w})$ is zero.  In this case
$\ind(\tilde{w}) = -\chi(\DD) + 1 = 0$, and Theorem~\ref{thm:windpi} gives
$$
0 \le 2\windpi(\tilde{w}) \le \ind(\tilde{w}) - 2 = -2,
$$
unless $E_{\omega_\infty}(\tilde{w}) = 0$.  The latter is also impossible
by Prop~\ref{prop:constants}, since $\tilde{w}$ is not constant.
\end{proof}

Having proved the gradient bound when $m \ge 2$, we now apply similar
arguments for $m \le 1$ and finish the proof of $\Cinftyloc$--covergence 
for this case.

\begin{prop}
\label{prop:stableCompactness}
The statement about $\Cinftyloc$--convergence in
Theorem~\ref{thm:compactness} holds if $\chi(\dot{\Sigma}) \ge 0$.
\end{prop}
\begin{proof}
This includes two cases: $\dot{\Sigma}$ is diffeomorphic to either a
plane or a singly punctured disk.
In both cases the space of conformal structures on the domain
is trivial, so we can assume $(\dot{\Sigma},j_k)$ is either
$(\CC,i)$ or $(\CC\setminus \dD, i)$ for all $k$, where $\dD = \interior{\DD}$.
We then have a sequence of maps
$\tilde{u}_k = (a_k,u_k) : \dot{\Sigma} \to \RR\times S^3$ satisfying
$T\tilde{u}_k \circ i = \tilde{J}_k \circ T\tilde{u}$,
all positively asymptotic at $\infty\in\Sigma$ to the simply
covered orbit $P_\infty$.

We address first the case $\dot{\Sigma} = \CC$; then there is no boundary
condition and $M = S^3$.  By the same rescaling
argument as in Prop.~\ref{prop:gradBounds},
we may assume after reparametrization that $|d\tilde{u}_k(z)|$
satisfies a global uniform bound and that $|d\tilde{u}_k(0)|$ is bounded
away from zero.  Then a subsequence of $(a_k - a_k(0),u_k)$ converges
in $\Cinftyloc$ to a nonconstant $\tilde{J}_\infty$--holomorphic finite 
energy plane $\tilde{u}_\infty = (a_\infty,u_\infty) : \CC \to \RR\times S^3$.
If the puncture is removable, then just as in Prop.~\ref{prop:gradBounds},
$\tilde{u}$ extends to a nonconstant sphere of index~$-2$, which violates
Theorem~\ref{thm:windpi}.  Thus $\tilde{u}_\infty$ is positively 
asymptotic to a periodic orbit $P$, which we claim must be $P_\infty$.
Indeed, if $P$ and $P_\infty$ are geometrically distinct, then the same
linking argument implies
$\lk(P,P_\infty) = 0$, a contradiction.  Suppose now that $P$ is an
$n$--fold cover of $P_\infty$ for some integer $n \ge 1$.
Writing $\omega_k = d\alpha_k$ by Lemma~\ref{lemma:exactTaming} and fixing
any smooth map $u : \CC \to S^3$ that approaches $P_\infty$ asymptotically,
$$
E_{\omega_k}(\tilde{u}_k) = \int_\CC u^*\omega_k = 
\int_{P_\infty} \alpha_k \to E_\infty := \int_\CC u^*\omega_\infty =
\int_{P_\infty} \alpha_\infty.
$$
But then
$$
n E_\infty = \int_{P} \alpha_\infty = \int_{\CC} u_\infty^*\omega_\infty
 \le \lim \int_{\CC} u_k^*\omega_k = E_\infty.
$$
The left hand side equals $E_\infty(\tilde{u}_\infty)$ and must therefore
be positive, so we conclude $n = 1$.

Next suppose $\dot{\Sigma} = \CC \setminus \dD$.  We claim that
$|d\tilde{u}_k|$ is uniformly bounded.  If not, then as in
Prop.~\ref{prop:gradBounds}, we can define rescaled maps $\tilde{v}_k$
on an increasing sequence of either disks or half-disks, depending on
whether and how fast $z_k$ approaches the boundary.  These then
have a subsequence convergent to a non-constant 
finite energy plane or half-plane
$\tilde{v}_\infty$.  The usual arguments show that if $\tilde{v}_\infty$
is a plane, it cannot be extended to a sphere, and the linking conditions
on its asymptotic orbit force it to intersect either $P_\infty$ or $K$, neither
of which is allowed.  For the half-plane case, $\tilde{v}$ extends to
a non-constant pseudoholomorphic disk, and the same argument as before 
shows that $\tilde{v}(\p \DD)$ is contractible on $L$, 
thus its Maslov index is~$0$, and it must therefore have vanishing
$\omega_\infty$--energy, another contradiction.

Given the uniform bound, there are constants $c_k \in \RR$ and a subsequence
of $(a_k + c_k,u_k)$ which converges in $\Cinftyloc$ to a
$\tilde{J}_\infty$--holomorphic finite energy map
$\tilde{u}_\infty : \dot{\Sigma} \to \RR\times M$, with the boundary
condition $\tilde{u}_\infty(\p\Sigma) \subset \{0\}\times L$.
Then it remains to
prove that $\tilde{u}_\infty$ has a positive puncture at $\infty$, asymptotic
to $P_\infty$ with covering number~$1$.
If the puncture is removable, we obtain a holomorphic disk
$$
\DD \to \RR\times S^3 : z \mapsto \tilde{u}_\infty(1/z)
$$
mapping $\p \DD$ to a meridian on $L$, thus $u_\infty$ must
intersect $K$.  But then $u_k$ for large $k$ would have to intersect the
interior of $N$, giving a contradiction.  Now suppose
$\tilde{u}_\infty$ is asymptotic to a periodic orbit $P$ at $\infty$.
If we extend $\tilde{u}_\infty$ to a smooth map over $\CC$, taking
$\dD$ into the solid torus $N$, then the same argument as in the plane case
shows that $P$ cannot be geometrically distinct from $P_\infty$.
Therefore $P$ is an $n$--fold cover of $P_\infty$
with $n \ne 0$.  (Here we adopt the convention of setting $n$ to
\emph{negative} the covering number if the puncture is negative; this
possibility is not excluded automatically when $\p\Sigma \ne \emptyset$.)
Now observe $u_k(\dot{\Sigma}) \cap K = \emptyset$ for all $k$, so if $k$ 
is sufficiently large, a small perturbation of $u_k$ realizes a
homology $\p[u_k] = n[P_\infty] + [\gamma]$ in $S^3\setminus K$,
where $\gamma$ is a negatively oriented meridian on $L$.
Consequently
$$
n \cdot \lk(P_\infty, K) = -\lk(\gamma, K) = 1,
$$
and since $\lk(P_\infty,K) > 0$ by assumption, $n$ can only be $1$.  
\end{proof}

\subsection{Convergence of conformal structures}
\label{subsec:conformal}

We now show
that in the case $\chi(\dot{\Sigma}) < 0$, the induced sequence of
conformal structures is compact.

\begin{prop}
\label{prop:jCompact}
Given the hypotheses of Theorem~\ref{thm:compactness},
there is a smooth complex structure $j_\infty$ on $\Sigma$ and a
sequence of diffeomorphisms $\varphi_k : \Sigma \to \Sigma$ fixing
$\infty$ and preserving each component of $\p\Sigma$, such that
a subsequence of $\varphi_k^* j_k$ converges to $j_\infty$ in 
the $C^\infty$--topology.  

There are also constants $c_k \in \RR$ such that a subsequence of
$(a_k + c_k,u_k) \circ \varphi_k$ converges in 
$\Cinftyloc(\dot{\Sigma},\RR\times M)$ to a map
$\tilde{u}_\infty \in \mM_{\hH_\infty,\Lambda}$, which is positively
asymptotic to $P_\infty$ at the puncture.
\end{prop}
\begin{proof}
A subsequence of $(\Sigma,j_k,\{\infty\})$ converges to a
stable nodal surface $\mathbf{S} = (S, j, \{p\}, \Delta, N)$.  
A choice of decoration $r$ defines
the compact connected surface $\overline{\mathbf{S}}_r$, with a singular
conformal structure $j_{\mathbf{S}}$ and singular
Poincar\'e metric $h_{\mathbf{S}}$, both of which degenerate on a finite set of
circles and arcs $\Theta_{\Delta,N} \subset \overline{\mathbf{S}}_r$.  Then
convergence means there is a sequence of diffeomorphisms
$$
\varphi_k : \overline{\mathbf{S}}_r \to \Sigma
$$
such that:
\begin{enumerate}
\item $\varphi_k(p) = \infty$.
\item $\varphi_k^*j_k \to j_{\mathbf{S}}$ in 
 $\Cinftyloc(\overline{\mathbf{S}}_r \setminus \Theta_{\Delta,N})$.
\item All circles in $\varphi_k(\Theta_{\Delta,N})$ are closed geodesics
 in $(\dot{\Sigma},h_k)$, and all arcs in $\varphi_k(\Theta_{\Delta,N})$
 are geodesic arcs in $(\dot{\Sigma},h_k)$ that intersect $\p\Sigma$
 transversely.
\end{enumerate}
We can assume without loss of generality that the diffeomorphisms
$\varphi_k$ map a given component of $\p (\overline{\mathbf{S}}_r)$
always to the \emph{same} component of $\p\Sigma$, i.e.~$\varphi_k \circ
\varphi_j^{-1}$ always 
maps each connected component $\gamma_j\subset \p\Sigma$ to itself.

If $S_j$ is a connected component of $S$, let $\dot{S}_j$ be the punctured
surface obtained by removing all points in the set $(\{p\} \cup \Delta
\cup N) \cap S_j$.  Note that the stability condition
implies $\chi(\dot{S}_j^D) < 0$.
There is a natural embedding $\dot{S}_j \hookrightarrow
\overline{\mathbf{S}}_r \setminus \Theta_{\Delta,N}$, which we use to
define the sequence of complex structures 
$\varphi_k^*j_k$ and metrics $\varphi_k^* h_k$ on $\dot{S}_j$.
Then passing to a subsequence, we have $\varphi_k^*j_k \to j$ and
$\varphi_k^* h_k \to h$ in $\Cinftyloc$ on $\dot{S}_j$, where
$h$ is the Poincar\'e metric for $(\dot{S}_j,j)$.
Since $d\tilde{u}_k$ is uniformly bounded on compact subsets,
we can then find constants $c_k^j \in \RR$ such that
$$
\tilde{v}_k^j = (b_k^j, v_k^j) = (a_k + c_k^j, u_k ) \circ 
 \varphi_k|_{\dot{S}_j} :
 (\dot{S}_j, \varphi_k^*j_k) \to  (\RR\times S^3, \tilde{J}_k)
$$
is a sequence of pseudoholomorphic maps satisfying the appropriate boundary
conditions and a uniform $C^1$--bound.
Thus $\tilde{v}_k^j$ has a $\Cinftyloc$--convergent subsequence
$$
\tilde{v}_k^j \to \tilde{v}^j = (b^j,v^j) : \dot{S}_j \to \RR\times S^3,
$$
where $\tilde{v}^j$ satisfies $T\tilde{v}^j \circ j = \tilde{J}_\infty \circ
T\tilde{v}^j$.  Due to the uniform energy bound for $\tilde{u}_k$, we see also
that $E_\infty(\tilde{v}^j) < \infty$.
Repeating this process for every component
$S_j \subset S$, we obtain a set of $\tilde{J}_\infty$--holomorphic maps
\begin{equation*}
\begin{split}
\tilde{v}^1 : \dot{S}_1 &\to \RR\times S^3, \\
            &\vdots \\
\tilde{v}^N : \dot{S}_N &\to \RR\times S^3.
\end{split}
\end{equation*}

Our main goal now is to show that $\mathbf{S}$ is actually a smooth
Riemann surface with boundary, i.e.~$\Delta$ and $N$ are empty sets and
$S$ has only one component.  Then the set of solutions above reduces to
a single solution $\tilde{u}_\infty : \dot{\Sigma} \to \RR\times M$,
which we must show is positively asymptotic to $P_\infty$ at the puncture.
As with the bubbling off arguments in the
previous section, these results 
will follow mainly from topological considerations.

Recall from Remark~\ref{remark:repeatedKnots} that our notation allows
some of the components $K_i$ and $K_j$ of $K$ to be identical; in particular,
topological considerations require a given component $K_i\subset K$ to 
repeat $n$ times in the list $K_1,\ldots,K_m$ if $\lk(K_i,P_\infty) = n$.
The lists of components 
$N = N_1 \cup \ldots \cup N_m$ and $L = L_1 \cup \ldots \cup L_m$
are then defined with similar repetitions.  If $\gamma_1,\ldots,\gamma_m$
are the connected components of $\p\Sigma$ (\emph{not} repeated),
then the oriented loop $u_k(\gamma_j)$ is a meridian on $L_j = \p N_j$
with $\lk(u_k(\gamma_j), K_j) = -1$.  Thus the linking number
$\lk(K_j,P_\infty)$ is precisely 
the number of distinct components of $\p\Sigma$
mapped into the same torus $L_j$, and we have also $\lk(u_k(\gamma_j), K) = -1$
since $u_k(\gamma_j)$ is unlinked with all components of $K$ that are
distinct from $K_j$.  Adding this up for all $\gamma_j\subset \p\Sigma$,
we see that the expression
$$
- \lk(u_k(\p\Sigma), K)
$$
counts the connected components of $\p\Sigma$.
Also, the map $u_k$ realizes a homology
$\p[u_k] = [P_\infty] + [u_k(\p\Sigma)]$ in $S^3 \setminus K$, which gives
the useful formula
\begin{equation}
\label{eqn:linkIdentity}
\lk(K_j,P_\infty) = - \lk(K_j, u_k(\p\Sigma)).
\end{equation}

In light of this topological setup, $u_k$ extends to a smooth map
$$
\bar{u}_k : \CC \to \RR\times S^3
$$
which satisfies $T\bar{u}_k \circ j_k = \tilde{J}_k \circ T\bar{u}_k$
in $\dot{\Sigma} = \CC\setminus (\dD_1 \cup
\ldots \cup \dD_m) \subset \CC$, and maps each of the disks $\dD_j$
into $N_j$.
We may assume that $\bar{u}_k|_{\dD_j}$ has a single transverse positive
intersection with $K_j$.

Let $S_1 \subset S$ be the connected component that contains the marked
point~$p$.

\textbf{Claim:} \textit{$\tilde{v}^1$ is positively asymptotic
to $P_\infty$ at $p$.}  
If the puncture is removable,
then we can find an oriented circle $C \subset \dot{S}_1$ winding
clockwise around $p$ such that $v^1(C)$ lies in an arbitrarily small
neighborhood of some point in $S^3\setminus K$.  Then this neighborhood
also contains $v^1_k(C) = u_k(\varphi_k(C))$ for sufficiently large $k$,
and $\varphi_k(C)$ is a large circle in $\CC$, bounding a simply connected
region $\Omega$.
One can then define a smooth map
$$
\hat{u}_k : \CC \to S^3\setminus K,
$$
which matches $u_k$ outside of $\Omega$, so that the
loops $\hat{u}_k(\p \DD_R)$ approach $P_\infty$ as $R\to\infty$.
This implies that for any component $K_j \subset K$, $\lk(P_\infty, K_j) = 0$,
a contradiction.

If $p$ is a nonremovable puncture and $\tilde{v}^1$ is asymptotic to an
orbit $P$ that is geometrically distinct from $P_\infty$, we similarly
find a large clockwise oriented circle $\varphi_k(C)\subset \CC$, bounding 
a region $\Omega$, such that $u_k(\varphi_k(C))$
is close to $P$.  Then the existence of the map 
$\bar{u}_k|_\Omega : \Omega \to S^3\setminus P_\infty$ implies
$\lk(P,P_\infty) = 0$, and this is impossible.  
The alternative is that $P$ could be an $n$--fold
cover of $P_\infty$ for some integer $n \ne 0$.  (Negative $n$ would again
mean the puncture is negative.)  But then restricting $u_k$ to the region
outside of $\Omega$ gives a homotopy of $u_k(\varphi_k(C))$ to $P_\infty$
in $S^3\setminus K$, implying that for any component $K_j \subset K$,
$$
n \cdot \lk(P_\infty, K_j) = \lk(u_k\circ\varphi_k(C),K_j) = \lk(P_\infty,K_j),
$$
so $n = 1$.  This proves the claim.

With the asymptotic behavior at $p$ understood, it remains to prove that
$\mathbf{S}$ has no double points or unpaired nodes.  Note that it
suffices to prove this for the component $S_1 \subset S$.  
We shall set up the
discussion in a slightly more general way than is immediately necessary,
since it will also be useful for the degeneration argument in the
next section.

First some notation.  The $m$ connected components of $\p\Sigma$ are denoted
$$
\p\Sigma = \gamma_1 \cup \ldots \cup \gamma_m,
$$
and let us write the components of $\p S_1$ as
$$
\p S_1 = \alpha_1 \cup \ldots \cup \alpha_s.
$$
Note that $m \ge 2$ by assumption, but $\p S_1$ could conceivably be
empty.  Assume $S_1$ has a (possibly empty) set of unpaired nodes
$$
N \cap S_1 = \{ w_1,\ldots,w_\ell \},
$$
interior double points
$$
\Delta \cap \interior{S_1} = \{ z_1,\ldots,z_q \},
$$
and boundary double points
$$
\Delta \cap \p S_1\supset \Delta \cap \alpha_j = \{ \zeta_j^1,\ldots,
 \zeta_j^{r_j} \} \quad \text{ for $j = 1,\ldots,s$},
$$
where we are regarding $\Delta$ for the moment as a set of \emph{points}
in $S$ rather than pairs of points.
We know from Theorem~\ref{thm:removeSing} that $\tilde{v}^1$ extends
smoothly over each boundary double point $\zeta_j^i \in \Delta\cap \p S_1$,
and at each $z_j \in \Delta\cap \interior{S_1}$ and $w_j \in N\cap S_1$, 
$\tilde{v}_1$ either has a
removable singularity or is asymptotic to some periodic orbit of
$X_\infty$.

Let $\overline{S}_1$ denote the compact surface with piecewise smooth
boundary
obtained from $S_1 \setminus ((\Delta \cup N) \cap S_1)$ by replacing the
interior punctures $z_j,w_j$ with circles at infinity $\delta_{z_j},
\delta_{w_j}$ and the boundary punctures $\zeta_j^i$ with arcs at infinity
$\delta_{\zeta_j^i}$.  Each component $\alpha_j
\subset \p S_1$ then gives rise to a piecewise smooth circular component
$\bar{\alpha}_j \subset \p\overline{S}_1$.  There is a natural map
$\overline{S}_1 \to \overline{\mathbf{S}}_r$, which is an inclusion
except possibly on $\p\overline{S}_1$, where two distinct circles
$\delta_{z_j}$ or arcs $\delta_{\zeta_j^i}$ may have the same image; this
corresponds to the identification of double points in a pair.
Since $\overline{\mathbf{S}}_r$ is diffeomorphic to
$$
\Sigma = \CC\cup \{\infty\} \setminus (\dD_1 \cup \ldots \cup \dD_m),
$$
we can visualize $\overline{\mathbf{S}}_r \setminus \{p\}$ as the plane
with a finite set of disks removed.  An example of this is shown in
Figure~\ref{fig:compactness}.  Here we settle on the convention
that the circles $\delta_{z_j}$ are always oriented
as components of $\p S_1$.  Thus they appear as embedded loops
winding \emph{clockwise} in the plane, and each encloses a bounded region which
may contain some of the disks $\dD_i$.  Let $m_j$ be the number of such disks
enclosed by $\delta_{z_j}$.  Similarly, for $j = 1,\ldots,s$, denote by
$\widehat{m}_j$ the number of disks in the compact region enclosed by 
$\bar{\alpha}_j$; this number is always at least~$1$.
Figure~\ref{fig:compactness} shows a compact subset of
$\overline{\mathbf{S}}_r$
which contains the entire boundary of $\overline{S}_1$.
Here the closure of the white area is $\overline{S}_1$, and the lightly
shaded regions constitute the rest of $\overline{\mathbf{S}}_r$, while 
the darkly shaded regions are the disks $\dD_j$.

\begin{figure}
\begin{center}
\includegraphics{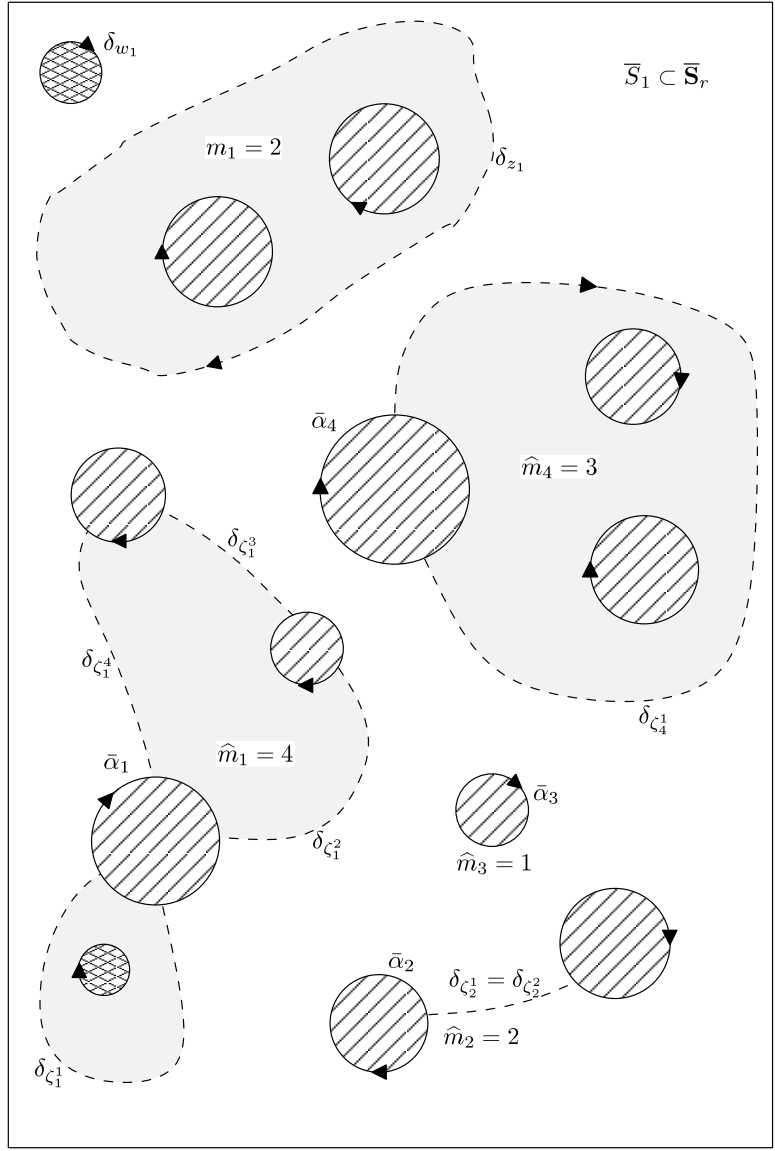}
\end{center}
\caption[Piecewise smooth circles in $\overline{\mathbf{S}}_r$]{
\label{fig:compactness} A compact subset of $\overline{\mathbf{S}}_r$
showing the piecewise smooth boundary of $\p \overline{S}_1$.  Here we
assume $\p S_1$ has four components $\alpha_1,\ldots,\alpha_4$, 
$S_1$ has one interior double point
$\Delta\cap \interior{S}_1 = \{z_1\}$, seven boundary double points
$\Delta\cap \alpha_1 = \{\zeta_1^1,\ldots,\zeta_1^4 \}$,
$\Delta\cap \alpha_2 = \{\zeta_2^1,\zeta_2^2 \}$,
$\Delta\cap \alpha_3 = \emptyset$,
$\Delta\cap \alpha_4 = \{\zeta_4^1\}$,
and one unpaired node $N \cap S_1 = \{w_1\}$.}
\end{figure}

The integers defined above are related by
\begin{equation}
\label{eqn:m}
m = \ell + \sum_{j=1}^q m_j + \sum_{j=1}^s \widehat{m}_j,
\end{equation}
and as remarked already,
\begin{equation}
\label{eqn:mhat}
\widehat{m}_j \ge 1 \quad \text{ for all $j = 1,\ldots,s$}.
\end{equation}
There are also constraints imposed by the stability condition for
each component of $S$: the double of $\dot{S}_1$ must have
negative Euler characteristic, thus
\begin{equation}
\label{eqn:S_1stability}
2 (s + q + \ell) + \sum_{j=1}^s r_j > 2,
\end{equation}
and applying similar reasoning to the portions of $\overline{\mathbf{S}}_r$ 
inside the loops $\delta_{z_j}$, we have
\begin{equation}
\label{eqn:m_j}
m_j \ge 2 \quad \text{ for all $j = 1,\ldots,q$}.
\end{equation}

We now transfer this picture onto $\dot{\Sigma}$ via the diffeomorphism
$$
\varphi_k : \overline{\mathbf{S}}_r \setminus \{p\} \to \dot{\Sigma}
$$
for large $k$ (see Figure~\ref{fig:compactness2}).  For $j=1,\ldots,q$,
denote by $\p_j\Sigma$ the $m_j$ components of $\p\Sigma$ that are enclosed 
within $\varphi_k(\delta_{z_j})$, and for $j=1,\ldots,s$ let $\hat{\p}_j\Sigma$
be the $\widehat{m}_j$ components in the closed region bounded by
$\varphi_k(\bar{\alpha}_j)$.
Now for each component $\alpha_j \subset \p S_1$, we define a perturbed
loop $\alpha_j' \subset \interior{S_1}$ which misses the double points.
The images $\varphi_k(\alpha_j') \subset \dot{\Sigma}$ are represented as
dotted loops in Figure~\ref{fig:compactness2}; each encloses a bounded 
region that contains $\hat{\p}_j\Sigma$.
Similarly, for each interior
double point $z_j$ we choose a perturbed loop 
$C_j \subset \interior{\overline{S}_1}$
near $\delta_{z_j}$, so $\varphi_k(C_j)$ encloses $\p_j\Sigma$.  Define
also the loops $\beta_j \subset \interior{\overline{S}_1}$ as
perturbations of $\delta_{w_j}$ for unpaired nodes $w_j \in N\cap S_1$: thus
each $\varphi_k(\beta_j)$ encloses a unique connected component 
$\gamma_{g(j)} \subset \p\Sigma$.
Observe that $\p\Sigma$ is now the disjoint union
$$
\p\Sigma = \left(\bigcup_{j=1}^q \p_j\Sigma\right) \cup
\left(\bigcup_{j=1}^s \hat{\p}_j\Sigma\right) \cup
\left(\bigcup_{j=1}^\ell \gamma_{g(j)} \right).
$$
The images under $\varphi_k$ of the various 
perturbed loops are shown with dotted lines in Figure~\ref{fig:compactness2}.

\begin{figure}
\begin{center}
\includegraphics{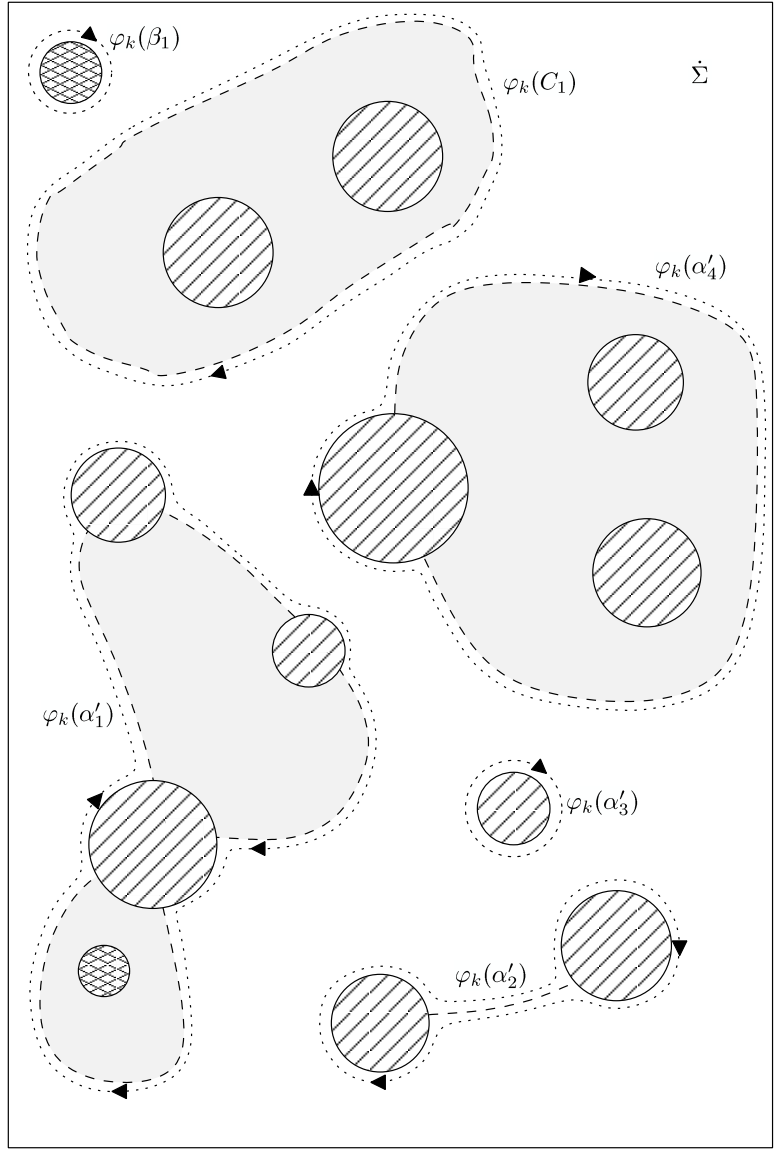}
\end{center}
\caption[Perturbed circles in 
 $\varphi_k(\overline{\mathbf{S}}_r)$]{\label{fig:compactness2} The 
 image of Figure~\ref{fig:compactness}
under $\varphi_k : \overline{\mathbf{S}}_r\setminus\{p\} \to \dot{\Sigma}$, 
showing the
perturbed loops $\alpha_1',\ldots,\alpha_4'$, $\beta_1$ and $C_1$ 
as dotted lines.}
\end{figure}

From this picture we can deduce some topological facts about the
behavior of $v^1 : \dot{S}_1 \to S^3$ at its boundary and punctures.
For a component $\alpha_j \subset \p S_1$, we have $v^1(\alpha_j)
\subset L_{f(j)}$ for some $f(j) \in \{1,\ldots,m\}$, and we can assume 
$u_k\circ\varphi_k(\alpha_j')$ is 
$C^0$--close to $v^1(\alpha_j)$.  Then restricting $u_k$ to the bounded
region inside $\varphi_k(\alpha_j')$ realizes a homology
$$
\p[u_k] = -[u_k\circ\varphi_k(\alpha_j')] + [u_k(\hat{\p}_j\Sigma)]
$$
in both $S^3\setminus P_\infty$ and $S^3\setminus K$.  This implies
$$
\lk(u_k\circ\varphi_k(\alpha_j'), P_\infty)
= \lk(u_k(\hat{\p}_j\Sigma), P_\infty) = 0,
$$
and thus
\begin{equation}
\label{eqn:vBoundaryMeridian}
\lk(v^1(\alpha_j), P_\infty) = 0.
\end{equation}
This means $v^1(\alpha_j)$ covers a meridian on $L_{f(j)}$, and its homotopy 
class can be deduced exactly via the linking number with $K$:
$$
\lk(v^1(\alpha_j), K) = \lk(u_k\circ\varphi_k(\alpha_j'), K) 
= \lk(u_k(\hat{\p}_j\Sigma), K) = -\widehat{m}_j.
$$
Since $v^1(\alpha_j)$ is only linked with one component of $K$,
\begin{equation}
\label{eqn:vBoundaryHomotopy}
\lk(v^1(\alpha_j), K_{f(j)}) = -\widehat{m}_j.
\end{equation}

Turning our attention next to the unpaired nodes, let us assume there is
a simply covered orbit $P_j \subset S^3$ of $X_{\infty}$ such
that $\tilde{v}^1$ is asymptotic to an $|n_j|$--fold cover of $P_j$ at
$w_j \in N\cap S_1$, for some $n_j\in\ZZ$.  Again, we're using the convention 
that the sign of $n_j$ matches the sign of the puncture at $w_j$, and
we set $n_j = 0$ if the puncture is removable (in which case it doesn't
matter what $P_j$ is).  Now, restricting $u_k$ to the region between
$\gamma_{g(j)}$ and $\varphi_k(\beta_j)$, we have a homology
$$
\p[u_k] = [u_k(\gamma_{g(j)})] - [u_k\circ\varphi_k(\beta_j)],
$$
in both $S^3\setminus P_\infty$ and $S^3\setminus K$,
and we can assume $[u_k\circ\varphi_k(\beta_j)]$ is homologous to
$n_j[P_j]$.  Thus for every component $K_i \subset K$,
\begin{equation}
\label{eqn:vUnpairedK}
n_j \lk(P_j, K_i) = \lk(u_k(\gamma_{g(j)}), K_i).
\end{equation}
Adding these up for all components of $K$, we find
$$
n_j \lk(P_j, K) = -1,
$$
implying that the puncture is nonremovable and the orbit is simply covered.
If $P_j = P_\infty$ this gives $n_j m = -1$, which cannot be true since
$m \ge 2$ by assumption.  Thus $P_j$ is geometrically distinct from $P_\infty$,
and using the homology in $S^3\setminus P_\infty$, we have
$n_j \lk(P_j, P_\infty) = \lk(u_k(\gamma_{g(j)}), P_\infty) = 0$, implying
\begin{equation}
\label{eqn:vUnpairedPinfty}
\lk(P_j, P_\infty) = 0.
\end{equation}

We can reach similar conclusions about the behavior of $\tilde{v}^1$ at
an interior double point $z_j \in \Delta\cap \interior{S_1}$.  Using the
same convention as above, assume $v^1$ approaches an $|n_j'|$--fold cover
of some simply covered orbit $P_j'$ at $z_j$.  Then we may assume
$[u_k\circ\varphi_k(C_j)]$ is homologous to $n_j' [P_j']$, and by restricting
$u_k$ over the bounded region inside $\varphi_k(C_j)$,
$$
\p[u_k] = [u_k(\p_j\Sigma)] - [u_k\circ\varphi_k(C_j)]
$$
in both $S^3\setminus K$ and $S^3\setminus P_\infty$.
This implies for all components $K_i \subset K$,
\begin{equation}
\label{eqn:vDoubleK}
n_j'\lk(P_j', K_i) = \lk(u_k(\p_j\Sigma), K_i),
\end{equation}
and summing this over the components of $K$, we have
$$
n_j' \lk(P_j', K) = - m_j \le -2,
$$
so $n_j'$ cannot be zero, i.e.~the puncture is not removable.
If $P_j' = P_\infty$, we have $n_j' m = - m_j$, then $m_j \le m$ implies
$n_j' = -1$ and $m = m_j$.  But this contradicts the stability assumption;
indeed, combining \eqref{eqn:mhat}, \eqref{eqn:m_j} and
\eqref{eqn:m}, we find $q = 1$ and $s = \ell = 0$, violating
\eqref{eqn:S_1stability}.  Therefore $P_j'$ is geometrically distinct
from $P_\infty$, and the homology in $S^3\setminus P_\infty$ gives
$n_j' \lk(P_j', P_\infty) = \lk(u_k(\p_j\Sigma), P_\infty) = 0$, thus
\begin{equation}
\label{eqn:vDoublePinfty}
\lk(P_j', P_\infty) = 0.
\end{equation}
At this point all the vital ingredients are in place.

\textbf{Claim:} \textit{$N \cap S_1$ and $\Delta \cap S_1$ are empty.}
If $w_j \in N\cap S_1$, then $v^1$ is asymptotic to a periodic orbit
$P_j$ which is geometrically distinct from $P_\infty$, and is also unlinked
with it, by \eqref{eqn:vUnpairedPinfty}.  But we have assumed there is
no such orbit, therefore $N\cap S_1 = \emptyset$.  The same argument
proves $\Delta\cap \interior{S_1} = \emptyset$, using
\eqref{eqn:vDoublePinfty}.  

It remains only to exclude double points on
the boundary.  We now can assume that $\p S_1 \ne \emptyset$ and the only
puncture of $\tilde{v}^1$ is at $p$, 
where it is positively asymptotic to $P_\infty$.  By assumption, there is
a trivialization $\Phi_\infty$ of $(v^1)^*\xi_\infty|_M$ 
for which $\muCZ^{\Phi_\infty}(P_\infty) = 3$ and,
using \eqref{eqn:vBoundaryHomotopy} and the fact that $v^1(\alpha_j)$ covers
a meridian for each component $\alpha_j \subset \p S_1$, 
the Maslov index along $\alpha_j$ is 
$2 \lk(v^1(\alpha_j),K_{f(j)}) = - 2 \widehat{m}_j$.  Thus we compute
$\mu(\tilde{v}^1) = 3 - 2 \sum_{j=1}^s \widehat{m}_j$, and
\begin{multline*}
\ind(\tilde{v}^1) = \mu(\tilde{v}^1) - \chi(\dot{S}_1) + s \\
= 3 - 2\sum_{j=1}^s \widehat{m}_j - (1 - s) + s 
= 2 + 2\sum_{j=1}^s \left(1 - \widehat{m}_j\right).
\end{multline*}
The $\omega_\infty$--energy 
of $\tilde{v}^1$ is clearly nonzero since $v^1(\p S_1)$
and the image of $v^1$ near $p$ cannot belong to the same orbit.
Thus Theorem~\ref{thm:windpi} gives
$$
0 \le 2\windpi(\tilde{v}^1) \le 2\sum_{j=1}^s \left(1 - \widehat{m}_j \right).
$$
Since $\widehat{m}_j \ge 1$ for all $j$, the right hand side of this expression
is never positive, and is zero if and only if $\widehat{m}_j = 1$ for all $j$.
This excludes situations such as $\bar{\alpha}_1$ and $\bar{\alpha}_2$ in 
Figure~\ref{fig:compactness}, where double points give rise to arcs that
connect two distinct disks.  All the arcs in $\delta_{\zeta^i_j} \subset
\p\overline{S}_1$ must therefore begin and end on the
same circle, enclosing a region of the plane as with $\bar{\alpha}_4$ in
the figure.  But now the stability condition requires this enclosed region
to have negative Euler characteristic after doubling, which can only be true
if it contains at least one disk, 
contradicting the fact that $\widehat{m}_j = 1$.
We conclude that there are no such arcs $\delta_{\zeta_j^i}$, and hence
no double points $\zeta_j^i \in \Delta \cap \p S_1$.

It follows now that $\mathbf{S}$ has no double points or unpaired nodes
at all, thus the convergence $(\Sigma,j_k,\{\infty\}) \to (S,j,\{p\},\Delta,N)$
simply means there are diffeomorphisms $\varphi_k : S \to \Sigma$ such that
$\varphi_k(p) = \infty$ and
$\varphi_k^*j_k \to j$ in $C^\infty(S)$.  Then after $\RR$--translation,
$\tilde{u}_k\circ \varphi_k \to \tilde{v}^1$ in $\Cinftyloc(S \setminus \{p\},
\RR\times S^3)$, and $\tilde{v}^1$ has the same asymptotic limit as 
$\tilde{u}_k$.  This completes the proof of Prop.~\ref{prop:jCompact}.
\end{proof}

\subsection{Degeneration at the boundary}
\label{subsec:noncompactness}

The proof of $\Cinftyloc$--convergence in 
Theorem~\ref{thm:degeneration} uses many of the same
arguments as Theorem~\ref{thm:compactness}, so we will not repeat these
in any detail, but rather emphasize the aspects that change when
the new orbits appear at $\p M$ in the limit.  As before, it's convenient
to treat the stable and non-stable cases separately.

\subsubsection*{The non-stable case}

The assumptions of Theorem~\ref{thm:degeneration} require that
$\p\Sigma$ be nonempty, so the only non-stable case to consider is
$m=1$: then $K$ is a knot with $\lk(P_\infty,K) = 1$,
and we may assume
$(\dot{\Sigma},j_k) = (\CC\setminus \dD, i)$ where $\dD =
\interior{\DD}$.  It will be convenient to use the biholomorphic map
$$
\psi : \RR \times S^1 \to \CC\setminus \{0\} :
(s,t) \mapsto e^{2\pi(s+it)}
$$
and consider the sequence of $\tilde{J}_k$--holomorphic half-cylinders
$$
\tilde{v}_k = (b_k, v_k) = \tilde{u}_k\circ \psi : [0,\infty)\times S^1
\to \RR\times S^3,
$$
with $v_k(\{s\}\times S^1) \to P_\infty$ as $s \to \infty$.
We claim $|d\tilde{v}_k|$ is uniformly bounded.  The proof is almost
identical to what was done in Prop.~\ref{prop:stableCompactness}: a
sequence $z_k$ with $|d\tilde{v}_k(z_k)| \to \infty$ gives rise to
a non-constant finite energy plane or disk.  A disk is impossible for the
same reasons as before: its boundary would have to be contractible on $L$,
leading to the conclusion that the map is constant.  A plane cannot be
asymptotic to any cover of $P_\infty$ or any orbit that is linked with it.
The only new feature is that \textit{a priori} the plane could be asymptotic 
to one of the orbits on $L$, but this would 
imply that $v_k$ intersects $K$ for 
large $k$, and is thus also excluded.

Pick an open neighborhood $\uU$ of $P_\infty$ in $M$, small enough so
that its closure does not intersect $L$.  Then define
$$
s_k = \min \{s \in [0,\infty)\ |\ v_k((s,\infty)\times S^1) 
 \subset \uU \}.
$$
\textbf{Claim:} \textit{$s_k \to \infty$.}
If not, there is a subsequence for which
$s_k \to s_\infty \in [0,\infty)$ and (in light of the gradient bound),
there are real numbers $c_k$ such that $(b_k + c_k, v_k)$ is
$\Cinftyloc$--convergent to a $\tilde{J}_\infty$--holomorphic half-cylinder
$$
\tilde{v} = (b,v) : [0,\infty)\times S^1 \to \RR\times S^3
$$
with finite energy.  Observe that $v(\{0\}\times S^1)$ is a
meridian on $L$.  Then if the puncture of $\tilde{v}$ is removable,
$v$ extends to a disk which must intersect $K$, implying that
some part of the image of $v_k$ lies inside the solid torus $N$ for
sufficiently large $k$, a contradiction.  

Since $\infty$ is not a removable puncture,
denote by $P$ its asymptotic orbit.  Now the usual linking
arguments imply
$$
\lk(P,K) = \lk(P_\infty,K) = 1,
$$
and if $P$ and $P_\infty$ are geometrically distinct,
$$
\lk(P,P_\infty) = 0.
$$
This leaves only two possibilities: $P$ is $P_\infty$ (positive puncture,
simply covered) or it is a simply covered orbit on $L$ (negative puncture).
In the latter case, the fact that $\omega_\infty$ is exact and vanishes
on $L$ implies $E_{\omega_\infty}(\tilde{v}) = 0$, thus
the image of $v$ is contained
in a periodic orbit, and therefore in $L$.  But this gives a contradiction,
because $v(\{s_\infty\}\times S^1)$ is in the closure of $\uU$, which is
disjoint from $L$.

There remains the possibility that $\tilde{v}$ is positively asymptotic
to $P_\infty$.  But then $\tilde{v}$ has precisely the same asymptotic
and boundary conditions as $\tilde{v}_k$, hence $\ind(\tilde{v}) = 2$
and Theorem~\ref{thm:windpi} gives $\windpi(\tilde{v}) = 0$.  This
implies $v$ is transverse to $X_\infty$, which is impossible
at $\{0\}\times S^1$ because both the image of $v$ and the orbits of
$X_\infty$ on $L$ are meridians.  This contradiction proves the claim that
$s_k \to \infty$.

With this established, define a sequence
$$
\tilde{w}_k = (\beta_k,w_k) : [-s_k,\infty) \times S^1 \to \RR\times S^3
$$
by $\tilde{w}_k(s,t) = \tilde{v}_k(s + s_k,t)$.  Then a subsequence of
$(\beta_k + c_k, w_k)$ converges in $\Cinftyloc(\RR\times S^1,\RR\times S^3)$ to a 
$\tilde{J}_\infty$--holomorphic finite energy cylinder
$$
\tilde{w}_\infty = (\beta_\infty,w_\infty) : \RR\times S^1 \to \RR\times S^3.
$$
The loop $\gamma := w_\infty(\{0\}\times S^1)$ is now the uniform limit of
$v_k(\{s_k\}\times S^1)$, and the usual arguments show that
$$
\lk(\gamma,P_\infty) = 0
\qquad\text{ and }\qquad
\lk(\gamma,K) = 1,
$$
thus $\tilde{w}_\infty$ cannot be a constant map.  If 
$E_{\omega_\infty}(\tilde{w}_\infty) = 0$, then these linking conditions
and the fact that $\gamma$ is in the closure of $\uU$ imply that
$\tilde{w}_\infty$ parametrizes $\RR\times P_\infty$.  However, there
exists a sequence $(s_k',t_k')$ with $s_k' < s_k$ and $s_k - s_k' \to 0$
such that $v_k(s_k',t_k') \not\in \uU$, implying that $\gamma$ also meets
the boundary of $\uU$, a contradiction.  Therefore 
$E_{\omega_\infty}(\tilde{w}_\infty) > 0$.

We shall now show that both punctures of $\tilde{w}_\infty$ are
positive and asymptotic to the appropriate orbits.  If both are
removable, we obtain a nonconstant holomorphic sphere of index~$-2$, 
contradicting Theorem~\ref{thm:windpi} as before.  If only one is removable,
then we can define a smooth map of a disk into $S^3\setminus K$ sending
the boundary to $\gamma$, implying the contradiction $\lk(\gamma,K) = 0$.
Now denote the two asymptotic orbits by 
$w_\infty(\{\pm\infty\} \times S^1) = P_\pm$.  We find,
$$
\lk(P_\pm, K) = \lk(P_\infty, K) = 1,
$$
and if $P_\pm$ is geometrically distinct from $P_\infty$,
$$
\lk(P_\pm,P_\infty) = 0.
$$
Therefore each orbit $P_\pm$ is either $P_\infty$ or is contained
in $L$, simply covered in either case.  We can determine the sign of each
puncture by comparing the orientations of $w_\infty(\{\pm\infty\}\times S^1)$
with the orientations of the orbits.  This allows four possibilities:
\begin{enumerate}
\setlength{\itemsep}{0in}
\item[(i)] $P_+ = P_\infty$ (positive puncture) and $P_- = P_\infty$ 
 (negative puncture)
\item[(ii)] $P_+ \subset L$ (negative puncture) and $P_- \subset L$ (positive 
 puncture)
\item[(iii)] $P_+ = P_\infty$ (positive puncture) and $P_- \subset L$ 
 (positive puncture)
\item[(iv)] $P_+ \subset L$ (negative puncture) and $P_- = P_\infty$
 (negative puncture)
\end{enumerate}
Case~(iv) is immediately excluded because both punctures can't be negative.
Cases~(i) and~(ii) would both imply $E_{\omega_\infty}(\tilde{w}_\infty)
= 0$, using again the fact that $\omega_\infty$ is exact and vanishes on $L$.
We conclude that both punctures are
positive, with $P_+ = P_\infty$ and $P_- \subset L$.

To apply this result to the sequence $\tilde{u}_k$, define a sequence of
diffeomorphisms
$$
\varphi_k : \CC \setminus\{0\} \to \CC\setminus \DD
$$
such that $\varphi_k(z) = e^{2\pi s_k}z$ for all $z$
with $|z| \ge 2 e^{-2\pi s_k}$.
Then observe that $\tilde{w}_k\circ\psi^{-1}(z) = \tilde{u}_k\circ\varphi_k(z)$
whenever $|z| \ge 2 e^{-2\pi s_k}$,
thus after $\RR$--translation, a subsequence of $\tilde{u}\circ\varphi_k$
converges in $\Cinftyloc(\CC\setminus\{0\},\RR\times S^3)$ to 
$$
\tilde{u}_\infty = \tilde{w}_\infty\circ\psi^{-1} : \CC\setminus\{0\} \to
 \RR\times S^3,
$$
which is asymptotic to $P_\infty$ at $\infty$ and an orbit on $L$ at $0$.
Clearly also $\varphi_k^* i \to i$ in $\Cinftyloc(S^2\setminus\{0\})$.
We have thus proved $\Cinftyloc$--convergence for 
Theorem~\ref{thm:degeneration} in the case $\chi(\dot{\Sigma}) \ge 0$.

\subsubsection*{The stable case}

Now assume $\chi(\dot{\Sigma}) < 0$.  The proof of 
Theorem~\ref{thm:degeneration} in this case will follow roughly the same 
sequence of steps as in
Theorem~\ref{thm:compactness}, with a few important differences.

\subsubsection*{Step~1: Gradient bounds}

We begin by establishing a bound
$$
| d\tilde{u}_k(z) |_k \le \frac{C}{\inj_{k}(z)}.
$$
The proof is mostly the same as in Prop.~\ref{prop:gradBounds}.
If a finite energy plane bubbles off, then it is asymptotic to an orbit
$P$ which (for topological reasons) cannot be a cover of $P_\infty$, and 
$\lk(P,P_\infty) = 0$.  The only remaining alternative
(which is new in this situation) is that $P$ is a meridian on one of the tori
$L_j$, but this would imply $\lk(P,K_j) \ne 0$, so $u_k(\dot{\Sigma})$ would
have to intersect $K_j$ for some large $k$.
The argument excluding disk bubbles is the same as before.

As in the proof of Prop.~\ref{prop:jCompact}, a subsequence of 
$(\Sigma,j_k,\{\infty\})$ converges
to a stable nodal surface $\mathbf{S} = (S,j,\{p\},\Delta,N)$.  
We again denote the connected components by
$S = S_1 \cup \ldots \cup S_N$ with corresponding punctured surfaces 
$\dot{S}_j$, choosing the labels so that
$p \in S_1$.  The gradient bound above implies that we
can find constants $c_k^j \in \RR$ such that
$$
(a_k + c_k^j, u_k) \circ \varphi_k|_{\dot{S}_j} \to
\tilde{v}^j : \dot{S}_j \to \RR\times S^3
$$
in $\Cinftyloc(\dot{S}_j, \RR\times S^3)$, where
$T\tilde{v}^j \circ j = \tilde{J}_\infty \circ T\tilde{v}^j$.
Our main goal will be to show that $\mathbf{S}$ has no double points and
no boundary, but does have $m$ unpaired nodes, one corresponding to 
each component of $\p\Sigma$.

\subsubsection*{Step~2: Asymptotic behavior at $p$}

The same arguments as in Prop.~\ref{prop:jCompact} show that
$p$ is a nonremovable puncture for $\tilde{v}^1 : \dot{S}_1 \to \RR\times S^3$,
and if $P$ is an asymptotic limit then either $P = P_\infty$ (simply covered)
or $P$ is geometrically distinct from $P_\infty$, with $\lk(P,P_\infty) = 0$.
In the present context this last possibility implies that $P$ is an $n$--fold
cover of some orbit $P_1$ on one of the tori $L_i$, with
$\lk(P, K_i) = n\cdot \lk(P_1, K_i) = -n$.  (As always, $n \ne 0$ and is
negative if the puncture is negative.)
Then we can choose a small circle $C$
about $p$ such that $u_k(\varphi_k(C))$ is close to $P$ for some large $k$,
and thus construct a homotopy from $P$ to $P_\infty$ through $S^3\setminus K$,
implying
$$
\lk(P, K_j) = \lk(P_\infty, K_j) > 0
$$
for each component $K_j \subset K$.  The left hand side is $0$ if 
$K_j \ne K_i$, so this alternative can only happen if $K$ is connected:
in that case $-n = \lk(P_\infty, K) = m$, so $p$ is a negative puncture
and $P$ is an $m$--fold cover of $P_1$.  We shall use arguments similar
to the proof of the non-stable case to show that this is also impossible.

Identify a punctured neighborhood of $p$ in $\dot{S}_1$ with the
positive half-cylinder via a holomorphic embedding
$$
\psi : [0,\infty)\times S^1 \hookrightarrow S_1\setminus\{p\},
$$
and define $\tilde{w}_k = (\beta_k,w_k) = \tilde{u}_k\circ\varphi_k \circ \psi 
: [0,\infty)\times S^1 \to \RR\times S^3$.  These half-cylinders are
$\tilde{J}_k$--holomorphic with the varying complex structures
$(\varphi_k \circ \psi)^*j_k = \psi^*\varphi_k^*j_k$ on the domain, and they
converge in $\Cinftyloc$ (possibly after translation in $\RR\times M$)
to $\tilde{v}^1 \circ \psi$.  Observe that $\varphi_k^*j_k \to j$ in 
$C^\infty$ on any compact
neighborhood of $p$, thus $\psi^*\varphi_k^*j_k \to \psi^*j = i$ 
in $C^\infty([0,\infty)\times S^1)$, not just on compact subsets;
this follows from Lemma~\ref{lemma:notJustCinftyloc} below.
The familiar argument then establishes a uniform bound on
$| d\tilde{w}_k |$ over $[0,\infty)\times S^1$: the alternative is that
$| d\tilde{w}_k(s_k,t_k)|$ blows up on some sequence with $s_k \to \infty$,
in which case a $\tilde{J}_\infty$--holomorphic finite energy
plane bubbles off, leading to the usual contradictions.

Due to the asymptotic behavior of $\tilde{v}^1$, there exists a sequence
$s_k \to \infty$ such that $w_k(s_k,\cdot)$ converges in $C^\infty(S^1,S^3)$
to a negatively oriented $m$--fold cover of the orbit $P_1$.  
But the 
half-cylinders $\tilde{w}_k$ are each asymptotic to $P_\infty$, thus we can
(as in the stable case) pick a small open neighborhood $P_\infty \subset \uU 
\subset M$ and define
$$
s_k' = \min \{s \in [0,\infty)\ |\ w_k((s,\infty)\times S^1) 
 \subset \uU \}.
$$
Clearly $s_k' > s_k$, thus $s_k' \to \infty$.  Now define
$\tilde{v}_k : [-s_k',\infty)\times S^1 \to \RR\times M$ by
$$
\tilde{v}_k(s,t) = (b_k(s,t),v_k(s,t)) = (\beta_k(s + s_k',t) - 
\beta_k(s_k',0), w_k(s + s_k',t)).
$$
These satisfy a uniform $C^1$--bound and are $\tilde{J}_k$--holomorphic
with respect to a sequence of complex structures which converge to
$i$ in $\Cinftyloc(\RR\times S^1)$, hence a subsequence converges to
a $\tilde{J}_\infty$--holomorphic finite energy cylinder
$$
\tilde{v}_\infty = (b_\infty,v_\infty) : \RR\times S^1 \to \RR\times M.
$$
As in the stable case, the loop $\gamma := v_\infty(\{0\}\times S^1)$
is necessarily nontrivial and not contained in a periodic orbit,
thus $\tilde{v}_\infty$ is nonconstant and 
$E_{\omega_\infty}(\tilde{v}_\infty) > \infty$.  The usual topological
constraints now imply that both punctures are nonremovable: in particular 
$\tilde{v}_\infty$ is asymptotic to $P_\infty$
at $+\infty$, and an $m$--fold covered orbit
on $L$ at $-\infty$, with both punctures positive.  
Denote the $m$--fold covered orbit on $L$ by $P_-$.

This leads to the following contradiction.
Let $\Psi$ denote the natural trivialization of $\xi_\infty$ along
$P_-$ defined by the intersection $TL \cap \xi_\infty$.  Then if
$e_-$ is the asymptotic eigenfunction at the puncture,
we claim $\wind^{\Psi}(e_-) = 0$.  Otherwise, we could find some
$s_0$ near $-\infty$ such that for large $k$, the loop 
$u_k\circ\varphi_k\circ\psi(s_0,\cdot)$ winds nontrivially around $P_-$, 
and must therefore
intersect $L$, which is a contradiction.  Then Lemma~\ref{lemma:MBparity}
gives $\muCZ^{\Psi-}(P_-) = 1$.  In terms of the given trivialization
$\Phi_\infty$ of $\xi_\infty|_M$, we have
$\wind^{\Phi_\infty}_{P_-}(\Psi) = -m$ and thus
$\muCZ^{\Phi_\infty-}(P_-) = 1 - 2m$.  Now
$\ind(\tilde{w}) = \muCZ(\tilde{w}) = 3 + 1 - 2m = 4 - 2m$, and
Theorem~\ref{thm:windpi} gives
$$
0 \le 2\windpi(\tilde{w}) \le \ind(\tilde{w}) - 2 + \#\Gamma_0
 = 2 - 2m.
$$
This is impossible, since we've assumed $m \ge 2$.
We're left with the alternative that $\tilde{v}^1$ is positively 
asymptotic to $P_\infty$ at the marked point $p$.

Before moving on, we should note the following lemma, which was used in
the argument above to prove $C^\infty$--convergence on the noncompact set
$[0,\infty)\times S^1$.
\begin{lemma}
\label{lemma:notJustCinftyloc}
Let $A_k : \DD \to \End(T\DD)$ be a sequence of smooth sections of the
tensor bundle $\End(T\DD) \to \DD$ such that $A_k \to 0$ in $C^\infty(\DD)$.
Then if $\psi : [0,\infty)\times S^1 \to \DD \setminus\{0\}$ is the 
biholomorphic map $\psi(s,t) = e^{-2\pi(s+it)}$, the tensors 
$\psi^*A_k$ on $[0,\infty)\times S^1$ converge uniformly to $0$ with all
derivatives.
\end{lemma}
\begin{proof}
Define the Euclidean metric on both $\DD$ and $[0,\infty)\times S^1$, and
use the natural coordinates on each to write
sections of $\End(T\DD)$ or $\End(T([0,\infty)\times S^1))$ as smooth real 
$2$-by-$2$ matrix valued functions.  If $\psi(s,t) = z$, then the first
derivative of $\psi$ at $(s,t)$ and its inverse can be written as
\begin{equation}
\begin{split}
\label{eqn:psiDerivs}
D\psi(s,t) &= -2\pi e^{-2\pi(s+it)} = -2\pi z, \\
D\psi^{-1}(z) &= -\frac{1}{2\pi z} = - \frac{1}{2\pi} e^{2\pi(s+it)},
\end{split}
\end{equation}
using the natural inclusion of $\CC$ in the space of real $2$-by-$2$ matrices.
Then
$$
(\psi^*A_k)(s,t) = D\psi^{-1}(z) \circ A_k(z) \circ D\psi(s,t) =
 e^{2\pi it} A_k(z) e^{-2\pi it},
$$
so $\| \psi^*A_k \|_{C^0} = \| A_k \|_{C^0} \to 0$ since the matrices on
either side of $A_k(z)$ are orthogonal.  We obtain convergence for all
derivatives by observing that for any multiindex $\alpha$,
$\p^\alpha (\psi^* A_k)(s,t)$ is a finite sum of expressions of the form
$$
c \cdot U \cdot e^{2\pi it} \cdot D^j A_k(z)(z,\ldots,z) \cdot e^{-2\pi it}
 \cdot V
$$
where $c$ is a real constant, $U$ and $V$ are constant unitary matrices
(i.e.~complex numbers of modulus $1$), and $j \le |\alpha|$.
This is clearly true for $|\alpha| = 0$ and follows easily for all $\alpha$
by induction, using \eqref{eqn:psiDerivs}.  The norm of this expression
clearly goes to $0$ uniformly in $(s,t)$ as $k \to \infty$.
\end{proof}

\subsubsection*{Step~3: Degeneration of $j_k$}

Most of the hard work for this step
was done in the proof of Prop.~\ref{prop:jCompact}; in particular,
the discussion surrounding Figures~\ref{fig:compactness} 
and~\ref{fig:compactness2} applies in the present situation as 
well.  The main difference here is that, since there are now orbits that are
unlinked with $P_\infty$, it is not so trivial to exclude interior
double points.  Unpaired nodes, of course, will not be excluded at all;
they will \emph{replace} the boundary.

\textbf{Claim:} \textit{$\Delta \cap S_1$ is empty.}
This will follow from similar algebraic relations to the ones 
that were previously used only to exclude boundary double points.
At any component $\alpha_j \subset \p S_1$, the homotopy class of
$v^1(\alpha_j)$ in $L_{f(j)}$ is fully determined by 
\eqref{eqn:vBoundaryHomotopy}, giving the Maslov index $-2 \widehat{m}_j$
with respect to the given trivialization $\Phi_\infty$ of $\xi_\infty|_M$.

The behavior at an unpaired node $w_j \in N\cap S_1$ is similarly 
constrained: by \eqref{eqn:vUnpairedPinfty}, the asymptotic limit
$P_j$ can only be one of the Morse-Bott orbits on some torus $L_i$.
Then \eqref{eqn:vUnpairedK} tells us the torus in question must be
$L_{g(j)}$, and since $\lk(P_j, K_{g(j)}) = -1$, the covering number
$n_j = 1$.  So $w_j$ is a positive puncture, and repeating the argument
from Step~2, the asymptotic eigenfunction has
zero winding relative to the natural framing determined by
$TL_{g(j)} \cap \xi_\infty$.
Lemma~\ref{lemma:MBparity} then gives
Conley-Zehnder index $1$ with respect to this framing.
The framing itself has winding number $-1$ along
$P_j$ with respect to the trivialization $\Phi_\infty$, which changes the
Conley-Zehnder index to $-2 + 1 = -1$.

Likewise at an interior double point $z_j \in \Delta\cap \interior{S_1}$,
the asymptotic limit $P_j'$ must belong to a Morse-Bott torus,
and summing \eqref{eqn:vDoubleK} over all components $K_i \subset K$ we have
$$
-n_j' = n_j' \lk(P_j', K) = \lk(u_k(\p_j\Sigma),K) = -m_j,
$$
so $z_j$ is a positive puncture with covering number $m_j$.  The 
Conley-Zehnder index with respect to the natural framing on the torus
is again $1$, but now the framing winds $-m_j$ times with respect to
$\Phi_\infty$, giving index $-2m_j + 1$.

We now compute the Maslov index
\begin{multline*}
\mu(\tilde{v}^1) = 3 + \ell(-1) + \sum_{j=1}^q (1 - 2m_j)
 - 2 \sum_{j=1}^s \widehat{m}_j \\
= 3 - \ell + q - 2\left( \sum_{j=1}^q m_j + \sum_{j=1}^s \widehat{m}_j\right),
\end{multline*}
and \eqref{eqn:FredholmIndex} gives
\begin{equation*}
\begin{split}
\ind(\tilde{v}^1) &= \mu(\tilde{v}^1) - \chi(\dot{S}_1) + s \\
&= 3 - \ell + q - 2\left( \sum_{j=1}^q m_j + \sum_{j=1}^s \widehat{m}_j\right)
- (1 - s - \ell - q) + s \\
&= 2 + 2\sum_{j=1}^q (1 - m_j) + 2\sum_{j=1}^s (1 - \widehat{m}_j).
\end{split}
\end{equation*}
We can assume that at least one of the sets $\p S_1$, $N\cap S_1$ and
$\Delta\cap \interior{S_1}$ is nonempty, in which case $v^1$ approaches
one of the tori $L_i$ \emph{somewhere}, while approaching $P_\infty$ at
the marked point $p$.  It follows that the image of $v^1$ is not contained
in any single periodic orbit, so $E_{\omega_\infty}(\tilde{v}^1) > 0$.
Thus Theorem~\ref{thm:windpi} gives
\begin{equation}
\label{eqn:vWindpi}
0 \le 2\windpi(\tilde{v}^1) \le
2\left( \sum_{j=1}^q (1 - m_j) + \sum_{j=1}^s (1 - \widehat{m}_j) \right).
\end{equation}
Recalling that always $m_j \ge 2$ and $\widehat{m}_j \ge 1$, we conclude
$q = 0$ and $\widehat{m}_j = 1$ for each $j$, so $\Delta \cap \interior{S}_1$
is empty, and by the same argument as in the proof of
Theorem~\ref{thm:compactness}, so is $\Delta \cap \p S_1$.

\textbf{Claim:} \textit{$\p S_1 = \emptyset$ and $\#N = m$.}
We've now established that $S$ can have only one connected component
(there are no double points to connect $S_1$ with anything else), thus
$\overline{\mathbf{S}}_r = \overline{S}_1 \cong \Sigma$, and
$m = s + \ell$.  We need to prove $s = 0$.  Having just shown that everything
on the right hand side of \eqref{eqn:vWindpi} vanishes, we have
$\windpi(\tilde{v}^1) = 0$, so $v : \dot{S}_1 \to S^3$ is immersed and
transverse to $X_{\infty}$.  But if $\p S_1 \ne \emptyset$ this
cannot be true, because $v^1(\p S_1)$ and all orbits of 
$X_{\infty}$ on $L_j$ are meridians.

By the above results, $S$ is a sphere with one marked point $p$ and
unpaired nodes $N = \{w_1,\ldots,w_m\} \subset S \setminus\{p\}$, so we
can identify it holomorphically with the Riemann sphere $(S^2,i)$, setting
$\infty := p$ and $\Gamma' := N$.  The diffeomorphisms $\varphi_k :
\overline{\mathbf{S}}_r \to \Sigma$ preserve $\infty$, and restricting them
to the interior they define diffeomorphisms
$$
\varphi_k : S \setminus \Gamma' \to \interior{\Sigma},
$$
with $\varphi_k^* j_k \to i$ in $\Cinftyloc(S\setminus\Gamma')$.
Moreover, after $\RR$--translation, 
$\tilde{u}_k \circ \varphi_k \to \tilde{v}^1$ in
$\Cinftyloc(S\setminus (\{\infty\}\cup\Gamma'), \RR\times S^3)$, and
$\tilde{v}^1$ has precisely the required asymptotic behavior at the punctures
$\infty$ and $w_j\in \Gamma'$.
This concludes the proof of $\Cinftyloc$--convergence for 
Theorem~\ref{thm:degeneration}.

\subsection{Convergence at the punctures}
\label{subsec:punctures}

To finish proving Theorems~\ref{thm:compactness} 
and~\ref{thm:degeneration}, it remains only to establish that the sequences
of maps $(a_k + c_k, u_k) \circ \varphi_k$ behave well on small neighborhoods
of the punctures and boundary.  This follows from the next three results.

\begin{lemma}
\label{lemma:cylCoords}
Let $j_k$ be a sequence of complex structures on $\dot{\DD} := 
\DD\setminus\{0\}$ such that
$j_k \to i$ in $\Cinftyloc(\dot{\DD})$, and take a sequence of
biholomorphic maps
$$
\psi_k : ([0,R_k) \times S^1,i) \to (\dot{\DD},j_k)
$$
for $R_k \in (0,\infty]$.
Then after passing to a subsequence, $R_k \to \infty$ and $\psi_k$ converges in 
$\Cinftyloc([0,\infty)\times S^1,\dot{\DD})$ to a biholomorphic map
$\psi : ([0,\infty) \times S^1,i) \to (\dot{\DD},i)$.
\end{lemma}
This can be proved by a routine bubbling off analysis for the embedded
holomorphic maps $\psi_k^{-1} : (\dot{\DD},j_k) \hookrightarrow 
(\RR\times S^1,i)$; we refer to \cite{Wendl:thesis} for the details.
With this preparation, we can reduce the problem of convergence at the
boundary and ends to the following two statements; we'll prove only the
second, since both use almost identical arguments.  Let
$\hH_k = (\xi_k,X_k,\omega_k,J_k)$ be a sequence of stable Hamiltonian
structures on a compact $3$--manifold $M$ with boundary, 
converging in $C^\infty$ to $\hH_\infty = 
(\xi_\infty,X_\infty,\omega_\infty,J_\infty)$, with associated almost complex
structures $\tilde{J}_k \to \tilde{J}_\infty$.  Assume also that the 
taming forms $\omega_k$ are exact.

\begin{prop}
\label{prop:endConvergence}
Assume $P \subset M$ is a nondegenerate periodic orbit of $X_k$ for all
$k \le \infty$.  Suppose $\tilde{v}_k = (b_k,v_k) : [0,\infty)\times S^1 \to
\RR\times M$ is a sequence of finite energy $\tilde{J}_k$--holomorphic maps
asymptotic to $P$, with uniformly bounded energy $E_k(\tilde{v}_k) < C$ 
and converging in $\Cinftyloc([0,\infty)\times S^1,
\RR\times M)$ to a $\tilde{J}_\infty$--holomorphic map
$\tilde{v}_\infty = (b_\infty,v_\infty) : [0,\infty)\times S^1 \to \RR\times M$,
also asymptotic to $P$.  Then for every sequence $s_k \to \infty$,
the loops $v_k(s_k,\cdot)$ converge in $C^\infty(S^1,M)$ 
to a parametrization of $P$.
\end{prop}

\begin{prop}
\label{prop:bndryConvergence}
Assume $L \subset M$ is a
$2$--torus which is tangent to all $X_k$ and is a Morse-Bott torus for 
$X_\infty$.  Let $R_k \to \infty$, and suppose
$\tilde{v}_k = (b_k,v_k) : [0,R_k] \times S^1 \to \RR\times M$ is a
sequence of $\tilde{J}_k$--holomorphic maps converging 
in $\Cinftyloc([0,\infty)\times S^1,\RR\times M)$ to a 
$\tilde{J}_\infty$--holomorphic half-cylinder
$\tilde{v}_\infty = (b_\infty,v_\infty) : [0,\infty)\times S^1 \to 
\RR\times M$, and satisfying a uniform energy bound
$E_k(\tilde{v}_k) < C$.  Assume also 
that $\tilde{v}_k(\{R_k\}\times S^1) \subset
\{c_k\}\times L$ for some sequence $c_k \in \RR$, and $\tilde{v}_\infty$
is asymptotic to a periodic orbit $P \subset L$ which is homotopic along
$L$ to each of the loops $v_k(\{R_k\}\times S^1)$.
Then for every sequence $s_k \in [0,R_k]$ with $s_k \to \infty$, the loops
$v_k(s_k,\cdot)$ have a subsequence convergent in $C^\infty(S^1,M)$ to a 
closed orbit homotopic to $P$ in $L$.
\end{prop}
\begin{proof}
We claim first that for any sequence $s_k \le R_k$ with $s_k \to \infty$,
$$
\int_{[s_k,R_k]\times S^1} v_k^*\omega_k \to 0.
$$
Indeed, the loop $v_k(s_0,\cdot)$ can be made arbitrarily close in 
$C^\infty(S^1,M)$ to a parametrization of $P$ by choosing $s_0$ and $k$ 
large enough.
Then for any $\epsilon > 0$, the exactness of $\omega_k$ 
implies that we can find $k_0 \in \NN$ and $s_0 \in [0,R_{k_0}]$ such that
$$
\int_{[s_0,R_k] \times S^1} v_k^*\omega_k < \epsilon
$$
for all $k \ge k_0$.  Since $v_k^*\omega_k$ is positive,
$\int_{[s_k,R_k]\times S^1} v_k^*\omega_k$ is bounded by this as soon as
$s_k \ge s_0$, proving the claim.

From this and the uniform energy bound, we use the same argument as in
Theorem~\ref{thm:removeSing}
to derive a uniform bound on $|d\tilde{v}_k|$: else a nonconstant
finite energy plane or half-plane with zero $\omega_\infty$--energy bubbles
off, contradicting Prop.~\ref{prop:constants}.  

Consider now a sequence
$s_k \le R_k$ with $s_k \to \infty$, and suppose $R_k - s_k$ is unbounded.
Then a subsequence of
$$
\tilde{w}_k : [-s_k,R_k - s_k] \times S^1 \to \RR\times M :
(s,t) \mapsto \tilde{v}_k(s + s_k,t)
$$
converges (after $\RR$--translation) in $\Cinftyloc(\RR\times S^1,\RR\times M)$
to a $\tilde{J}_\infty$--holomorphic finite energy cylinder
$\tilde{w}_\infty = (\beta_\infty,w_\infty) : \RR\times S^1 \to \RR\times M$
with $E_{\omega_\infty}(\tilde{w}_\infty) = 0$.  By Prop.~\ref{prop:constants},
such an object is either constant or an orbit cylinder.  
To rule out the former, we claim
$$
\int_{\{0\}\times S^1} w_\infty^*\lambda_\infty =
\lim \int_{\{s_k\}\times S^1} v_k^*\lambda_k = \int_P \lambda_\infty =: Q_0
\ne 0.
$$
Indeed, since $v_k(R_k,\cdot)$ is homotopic to $P$ and $d\lambda_k$ vanishes
along $L$, we can assume for $k$ sufficiently large that
$$
\left| Q_0 - \int_{\{R_k\}\times S^1} v_k^*\lambda_k \right| < \epsilon.
$$
Then assuming also $\int_{[s_k,R_k]\times S^1} v_k^*\omega_k < \epsilon$ and
using Lemma~\ref{lemma:omegaMetric},
\begin{equation*}
\begin{split}
\left| Q_0 - \int_{\{s_k\}\times S^1} v_k^*\lambda_k \right| &=
\left| Q_0 - \int_{\{R_k\}\times S^1} v_k^*\lambda_k
 + \int_{[s_k,R_k]\times S^1} v_k^*d\lambda_k \right| \\
&\le (1 + C)\epsilon.
\end{split}
\end{equation*}
Consequently, $\tilde{w}_\infty$ is a trivial cylinder over a closed orbit with
the same period as $P$.

If instead $R_k - s_k$ remains bounded as $k \to \infty$, we consider
the maps 
$$
\tilde{w}_k : [-R_k,0]\times S^1 \to \RR\times M :
(s,t) \mapsto \tilde{v}_k(s + R_k,t),
$$
which satisfy the boundary condition $\tilde{w}_k(\{0\}\times S^1) \subset
\{c_k\} \times L$.  Then a subsequence converges in 
$\Cinftyloc((-\infty,0]\times S^1,\RR\times M)$ after 
$\RR$--translation
to a $\tilde{J}_\infty$--holomorphic finite energy half-cylinder
$\tilde{w}_\infty : (-\infty,0]\times S^1 \to \RR\times M$, with
$E_{\omega_\infty}(\tilde{w}_\infty) = 0$.  Repeating the argument above,
$\tilde{w}_\infty$ parametrizes half of an orbit cylinder.

We've shown now that for every sequence $s_k \to \infty$ with $s_k \le R_k$,
the sequence of loops $v_k(s_k,\cdot)$ has a subsequence converging in
$C^\infty(S^1,M)$ to a closed orbit of $X_\infty$.  We claim finally that
this orbit lies in $L$.  If not, then we can find a sequence
$s_k' \in (s_k,R_k)$ such that the loops $v_k(s_k',\cdot)$ touch the
boundary of a small neighborhood of $L$.  But by the Morse-Bott condition,
we may assume this neighborhood contains no other closed orbits of the same
period, thus $s_k'$ can have no subsequence for which 
$v_k(s_k',\cdot)$ converges to an appropriate orbit, giving a
contradiction.
\end{proof}

\section{The main construction}
\label{sec:mainConstruction}

\subsection{Surgery and Lutz twists on transverse links}
\label{subsec:DehnLutz}

We now define precisely the type of 
surgery on contact manifolds that we wish to perform.
In the following, $S^1$ is always defined to be the quotient
$\RR / \ZZ$.

\begin{lemma}
\label{lemma:lambdas}
Let $(\theta,\rho,\phi)$ be the standard cylindrical polar coordinates
on $S^1 \times \RR^2$, oriented by the basis $(\p_\theta,\p_\rho,\p_\phi)$.
For $\rho \ge 0$, choose real-valued 
functions $f(\rho)$ and $g(\rho)$ such that
$(\rho,\phi) \mapsto f(\rho)$ and $(\rho,\phi) \mapsto g(\rho) / \rho^2$
define smooth functions on $\RR^2$.  
Then 
$$
\lambda := f(\rho)\ d\theta + g(\rho)\ d\phi
$$
defines a smooth $1$--form on $S^1\times \RR^2$, which is a
positive contact form if
and only if the following two conditions are met:
\begin{itemize}
\item[(i)] The Wronskian $D(\rho) := f(\rho)g'(\rho) - f'(\rho)g(\rho) > 0$
for all $\rho > 0$.
\item[(ii)] $f(0) g''(0) > 0$.
\end{itemize}
In that case, the corresponding Reeb vector field is given by
\begin{equation}
\label{eqn:Reeb}
X(\theta,\rho,\phi) = \frac{1}{D(\rho)} 
(g'(\rho) \p_\theta - f'(\rho) \p_\phi).
\end{equation}
\end{lemma}
\begin{proof}
A simple calculation shows that
$$
\lambda \wedge d\lambda = D(\rho)\ d\theta \wedge d\rho \wedge d\phi
 = \frac{D(\rho)}{\rho}\ d\theta \wedge dx \wedge dy,
$$
where $(x,y)$ are Cartesian coordinates on $\RR^2$.  Then
$\lim_{\rho\to 0} D(\rho)/\rho = D'(0) = f(0)g''(0)$, and it is
straightforward to verify that the expression for $X$ above satisfies
$d\lambda(X,\cdot) \equiv 0$ and $\lambda(X) \equiv 1$.
\end{proof}
Intuitively, these conditions on $f$ and $g$ mean that the curve
$\rho\mapsto(f(\rho),g(\rho))$ always winds counterclockwise around
the origin in the $xy$--plane, beginning on the $x$--axis with zero
velocity and nonzero angular acceleration.  

\begin{figure}
\includegraphics{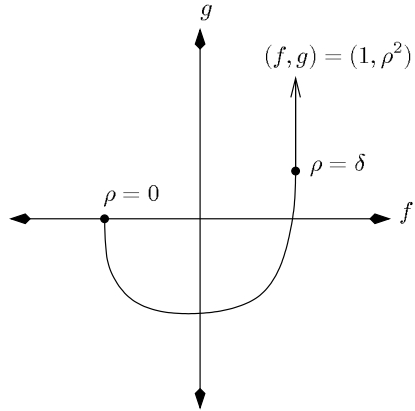}
\caption[Half Lutz twist]{\label{fig:halfLutz} Half Lutz twist of $\lambda_0$}
\end{figure}

\begin{figure}
\includegraphics{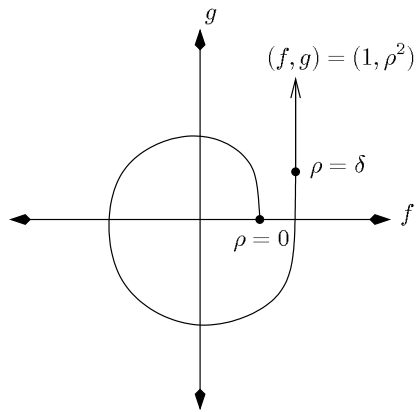}
\caption[Full Lutz twist]{\label{fig:fullLutz} Full Lutz twist of $\lambda_0$}
\end{figure}

Let $(M,\xi)$ be an oriented $3$--manifold with a positive and 
cooriented contact structure, and
suppose $K \subset M$
is an oriented knot which is \emph{positively transverse} to $\xi$,
i.e.~its orientation matches the coorientation of $\xi$.
A \emph{Lutz twist} along $K$ is defined as follows.
By the contact neighborhood theorem, $K$ has a solid torus neighborhood
$N_K \subset M$ which can be identified with 
$S^1\times \overline{B^2_\epsilon(0)}$,
where $\overline{B^2_\epsilon(0)}$ is the closed ball of radius
$\epsilon$ around the origin in $\RR^2$, such that
$K = S^1 \times \{0\}$ and $\xi|_{N_K}$ is the kernel of
\begin{equation}
\label{eqn:lambda0}
\lambda_0 := d\theta + \rho^2 \ d\phi.
\end{equation}
using the cylindrical polar coordinates of Lemma~\ref{lemma:lambdas}.  
We then change $\xi$ by replacing $\lambda_0$ on $N_K$ with
$\lambda_K := f(\rho)\ d\theta + g(\rho)\ d\phi$,
where $f$ and $g$ are functions chosen as in Lemma~\ref{lemma:lambdas}
so that $\lambda_K$ is a contact form, and furthermore:
\begin{enumerate}
\setlength{\itemsep}{0in}
\item There exists $\delta \in (0,\epsilon)$ such that 
$(f(\rho),g(\rho)) = (1,\rho^2)$ for $\rho \ge \delta$.
\item The trajectory $\rho \mapsto (f(\rho),g(\rho))$ rotates at least
halfway around the origin for $\rho \in [0,\delta]$
(see Figures~\ref{fig:halfLutz} and~\ref{fig:fullLutz}).
\end{enumerate}
This operation produces a new contact structure $\xi_K = \ker\lambda_K$,
such that there exists at least one radius $\rho_0 \in (0,\delta)$ at which
$g(\rho_0) = 0$.  This means
the meridian $\{ (0,\rho_0,\phi)\ |\ \phi \in \RR / 2\pi\ZZ \}$
is Legendrian and forms the boundary of an overtwisted disk, so that
$\xi_K$ is necessarily overtwisted.  Relatedly, there is at least one
radius $\rho_1 \in (0,\rho_0)$ at which $g'(\rho_1) = 0$ and $f'(\rho_1) > 0$,
so that the Reeb vector field on the torus $\{ \rho = \rho_1 \}$ generates
periodic orbits which are negatively oriented meridians.  This 
detail will be important for
constructing finite energy foliations in such neighborhoods.

The twists shown in Figures~\ref{fig:halfLutz} and~\ref{fig:fullLutz} may
be called the ``half Lutz twist'' and ``full Lutz twist'' respectively:
the former changes the homotopy class of $\xi$ as a $2$--plane distribution,
while the latter does not (see \cite{Bennequin} for an explicit homotopy). 
The half Lutz twist is particularly important for the following reason:
by an obstruction theory argument due originally to
Lutz \cites{Lutz:thesis,Lutz:77}, any homotopy class of cooriented 
$2$--plane distributions on $M$ admits a positive contact structure,
which can be obtained from any other $\xi$ 
by half Lutz twists along some positively transverse link.
See \cite{Geiges:contact} for a fuller discussion of this result.

We next generalize this to a twisting version of nontrivial Dehn surgery on 
contact manifolds.
Assume $M = S^3$ with positive contact structure $\xi$ and
positively transverse knot $K \subset S^3$.  Identify a neighborhood
$N_K$ of $K$ once more with $S^1 \times B^2_\epsilon(0)$, requiring
in particular that $\xi|_{N_K}$ be the kernel of $\lambda_0
= d\theta + \rho^2\ d\phi$ and that the longitude
$\{ (\theta,\epsilon,0) \ |\ \theta \in S^1 \} \subset S^3$ be
homologous to zero in $S^3 \setminus K$.  
Let $\lambda_1 = f_1(\rho)\ d\theta + g_1(\rho)\ d\phi$ be a contact
form on $N_K$ obtained from $\lambda_0$ by a Lutz twist as described above,
choosing $f_1$ and $g_1$ so that there is a radius $\rho_1 \in (0,\epsilon)$ 
at which $g_1'(\rho_1) = 0$ and $f_1'(\rho_1) > 0$, while $g_1'(\rho) > 0$ for
all $\rho \in (\rho_1,\epsilon]$.

A \emph{framing} of $K$ is a number $p/q \in \QQ \cup \{\infty\}$,
where we assume $p$ and $q$ are relatively prime integers.
Define now another solid torus $N' := S^1 \times \overline{B^2_\epsilon(0)}$,
with canonical cylindrical polar coordinates 
$(\theta',\rho',\phi')$.  Define also $\eta := \phi / 2\pi$,
$\eta' := \phi' / 2\pi \in S^1$, so that the pairs $(\theta,\eta)$,
$(\theta',\eta') \in S^1 \times S^1$ give coordinates on the tori 
$\p N_K$ and $\p N'$ respectively.
Now choose $\delta \in (0,\rho_1)$ and define an embedding
$$
\psi : N' \setminus (S^1 \times B^2_\delta(0)) \hookrightarrow N_K :
(\theta',\rho',\phi') \mapsto (\theta(\theta',\phi'),
\rho',\phi(\theta',\phi')),
$$
where the map $(\theta',\phi') \mapsto (\theta,\phi)$ is
determined by an orientation preserving 
diffeomorphism $\p N' \mapsto \p N_K$ of the form
\begin{equation}
\label{eqn:surgeryMatrix}
\begin{pmatrix} \theta \\ \eta \end{pmatrix} =
\begin{pmatrix} n & q \\ m & p \end{pmatrix}
\begin{pmatrix} \theta' \\ \eta' \end{pmatrix},
\end{equation}
for some matrix in $\SL(2,\ZZ)$.  A new manifold $M_K$ is defined by
removing $S^1 \times B^2_\delta(0)\subset N_K$ from $S^3$ and gluing 
in $N'$ via this embedding.  The topological type of $M_K$ depends only on
$p/q \in \QQ \cup \{\infty\}$ (see \cite{Saveliev}).

The contact form $\lambda_1$ on $N_K$ pulls back via $\psi$ to a
contact form $\lambda_K$ on $N' \setminus (S^1 \times B^2_\delta(0))$,
which in the coordinates $(\theta',\rho',\phi')$
has the form
\begin{equation*}
\begin{split}
\lambda_K &= f_K(\rho')\ d\theta' + g_K(\rho') d\phi' \\
&= \left[ n f_1(\rho') + 2\pi m g_1(\rho') \right]\ d\theta'
+ \left[ \frac{q}{2\pi} f_1(\rho') + p g_1(\rho') \right]\ d\phi'.
\end{split}
\end{equation*}
Clearly $\xi_K := \ker\lambda_K$ has a natural extension to $M_K \setminus N'$,
and we can extend $f_K$ and $g_K$ to $\rho \in [0,\delta]$ so that
$\lambda_K$ becomes a contact form on $N'$.  

We will refer to the operation described above as a \emph{rational}
\emph{twist surgery} along $K$, or \emph{integral} in the case where
$p/q \in \ZZ \cup \{\infty\}$, i.e.~$q = \pm 1$ or~$0$.
Observe that when $q = 0$, we can choose the surgery 
matrix to be the identity, which gives simply a Lutz twist along $K$.
By the theorem of Lickorish \cite{Lickorish} and Wallace \cite{Wallace},
every closed oriented $3$--manifold $M$ can be obtained by integral surgery
along some link $K \subset S^3$; then making $K$ positively transverse by
a $C^0$--perturbation, the procedure we've described proves the result of
Martinet \cite{Martinet} that every such $3$--manifold admits a 
cooriented and positive contact structure $\xi$.

Let $N \subset S^3$ be the union of all the solid tori removed from
$S^3$ in the above gluing, and let $N' \subset M$ be the corresponding
solid tori that are glued in, so there's a natural diffeomorphism
$S^3 \setminus N = M \setminus N'$.  The Lutz argument
now provides a link $K_1 \subset M$, positively transverse to $\xi_K$,
such that the contact structure $\xi_1$ obtained by half Lutz twists
along the components of $K_1$ may represent any desired homotopy class.
In particular, we can specify the homotopy class of $\xi_1$ over the
complement of a small $3$--ball $B \subset M$ by choosing any link $K'$ 
representing the appropriate class in $H_1(M)$, then perturbing it to be 
positively transverse---note that we're thus free to assume 
$K' \subset M\setminus N'$.
We can also assume $B \subset M\setminus N'$, and then the homotopy class 
of $\xi'$ over $B$ can be
changed as needed by twisting along a transverse link $K'' \subset B$
with the appropriate self-linking number.  In summary, we can assume
the transverse link $K_1$ needed to change the homotopy class of $\xi_K$
lies in $M \setminus N' = S^3 \setminus N$, and this leads to 
be following statement of the famous Lutz-Martinet theorem.

\begin{thm}[Lutz, Martinet]
\label{thm:LutzMartinet}
Given a closed oriented $3$--manifold $M$ with a cooriented $2$--plane
distribution $\alpha$, there exists a positive contact structure $\xi$
homotopic to $\alpha$, and $(M,\xi)$ can be obtained from the tight
three-sphere $(S^3,\xi_0)$ by a rational twist surgery along some
transverse link.
\end{thm}

\begin{remark}
By Eliashberg's classification theorem for overtwisted contact structures
\cite{Eliashberg:overtwisted}, the procedure above produces every
overtwisted contact structure on every closed $3$--manifold.
\end{remark}

The main goal of this paper is to construct finite energy foliations on
contact manifolds obtained from $(S^3,\xi_0)$ by twist surgeries.
The following result will be helpful for establishing that these foliations
can be made to have certain nice properties, e.g.~that all punctures are 
positive and all orbits are simply covered.

\begin{prop}
\label{prop:important}
For the surgery in Theorem~\ref{thm:LutzMartinet}, we can assume without
loss of generality any or all of the following:
\begin{enumerate}
\setlength{\itemsep}{0in}
\item The surgery is integral.
\item $\xi$ is overtwisted.
\item For each component $K_j \subset K$, let $\rho_1 \in (0,\epsilon)$ be
the largest radius where $g_1'(\rho_1) = 0$.  Then for all 
$\rho \in (0,\rho_1]$,
$$
f_K'(\rho) g_K''(\rho) - f_K''(\rho) g_K'(\rho) > 0.
$$
\item If the surgery at $K_j$ is topologically nontrivial (i.e.~$q \ne 0$), 
$\rho_1$ is the radius above and $r \in (0,\epsilon)$ is the smallest 
radius at which $f_K'(r) / g_K'(r) = f_K'(\rho_1) / g_K'(\rho_1)$, then
$$
\frac{f_K'(r)}{g_K'(r)} - \frac{f_K''(0)}{g_K''(0)}
$$
is positive and close to~$0$.
\end{enumerate}
\end{prop}
\begin{proof}
The Lickorish-Wallace theorem guarantees that integral surgeries are
sufficient.  If the surgery is topologically nontrivial (i.e.~not merely
a Lutz twist) at each component of $K$, then the contact structure $\xi$
defined as above need not generally be overtwisted.  We can, however,
\emph{make} it overtwisted by performing an extra \emph{full} Lutz twist
along a transverse knot disjoint from $K$; this doesn't change the homotopy
class of $\xi$.  Condition~(3) simply means that the trajectory
$\rho \mapsto (f(\rho),g(\rho))$ in $\RR^2$ has nonzero \emph{inward}
acceleration for all $\rho \in (0,\rho_1]$.
Condition~(4) is a relation between the slopes of the trajectory
at $\rho=0$ and $\rho=r$, and can always be achieved by changing
$f_K$ and $g_K$ near $\rho=0$.  This change may involve an extra half
Lutz twist, which changes the homotopy class of $\xi$, but it can be
changed back by adding Lutz twists along extra transverse knots.
\end{proof}

\begin{figure}
\includegraphics{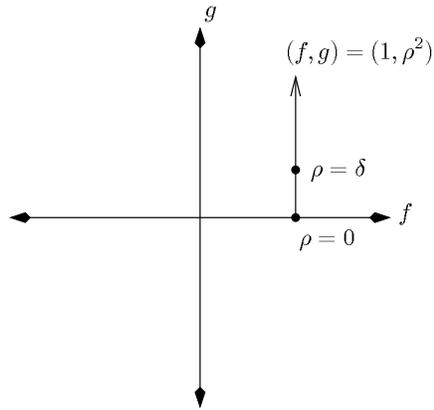}
\caption[The standard contact structure as a trajectory in 
 $\RR^2$]{\label{fig:stdlambda} The trajectory 
$\rho\mapsto (f(\rho),g(\rho))$ for the contact structure
$\lambda_0 = f(\rho)d\theta + g(\rho)d\phi = d\theta + \rho^2d\phi$}
\end{figure}

It will be helpful to know that transverse links in the tight $3$--sphere
$(S^3,\xi_0)$ can be assumed after transverse isotopy to approximate covers
of Hopf circles, as the latter admit coordinate neighborhoods in which
the standard contact form takes an especially simple form.
In the following, we view $S^3$ as the unit sphere in $\CC^2$, with
\begin{equation}
\label{eqn:standardS3}
\lambda_0(z)v := \frac{1}{2} \langle iz,v \rangle,
\end{equation}
for $z \in S^3 \subset \CC^2$ and $v \in T_zS^3 \subset \CC^2$,
where $\langle \ ,\ \rangle$ is the standard Euclidean inner product.
It was shown by Bennequin \cite{Bennequin} that all transverse links in
the standard contact $\RR^3$ are transversely isotopic to closed braids
about the $z$--axis.  Using a contact embedding of $\RR^3$ into the
tight $3$--sphere, one can then prove the following (see \cite{Wendl:thesis}
for details):

\begin{lemma}
\label{lemma:almostHopf}
Let $K \subset S^3$ be a link positively transverse to the standard
contact structure $\xi_0$.
Then for each component $K_j\subset K$ there is a smooth immersion
$F_j: [0,1] \times S^1 \to S^3$ such that 
$F_j(1,\cdot) : S^1\to S^3$ parametrizes $K_j$, 
$F_j(0,t) = (e^{2\pi i k_j t},0)$ for some $k_j\in\NN$, and for all fixed
$\tau\in (0,1]$, the collection of maps $F_j(\tau,\cdot) : S^1 \to S^3$ 
parametrizes a transverse link.
\end{lemma}

\subsection{Some simple foliations in $S^1 \times \RR^2$}
\label{subsec:local}

In this section we construct stable finite energy foliations on
the local neighborhoods that arise from twist surgery on transverse links.
Let $M = S^1 \times \RR^2$, with cylindrical polar coordinates
$(\theta,\rho,\phi)$ as in the previous section.  Then
using Lemma~\ref{lemma:lambdas}, define a positive contact form
$$
\lambda = f(\rho)\ d\theta + g(\rho)\ d\phi,
$$
where
$$
D(\rho) := f(\rho)g'(\rho) - f'(\rho)g(\rho) > 0
\text{ for all $\rho > 0$, and }
f(0)g''(0) > 0.
$$
Notice that smoothness at $\rho=0$ requires
$g(0) = f'(0) = g'(0) = 0$.

The contact structure $\xi = \ker\lambda$ is spanned for all 
$\rho > 0$ by the two vector fields
\begin{equation}
\label{eqn:v1v2}
v_1(\theta,\rho,\phi) = \p_\rho, \quad\quad\quad
v_2(\theta,\rho,\phi) = \frac{1}{D(\rho)}
( -g(\rho) \p_\theta + f(\rho) \p_\phi ),
\end{equation}
and the Reeb vector field $X$ is given by \eqref{eqn:Reeb}.
The flow of $X$ and its linearization are quite easy to compute,
leading to the following characterization of periodic orbits:

\begin{prop}
\label{prop:symOrbits}
Suppose $r > 0$ and $f'(r) / 2\pi g'(r) = p/q \in \QQ \cup \{\infty\}$
for relatively prime integers $p$ and $q$, whose signs match the signs
of $f'(r)$ and $g'(r)$ respectively.  Then the torus
$$
L_r := \{\rho = r\} \subset M
$$ 
is foliated by closed orbits of the form
$$
x(t) = \left(\theta_0 + \frac{g'(r)}{D(r)}t,
 \ r,
 \ \phi_0 - \frac{f'(r)}{D(r)}t\right)
= \left(\theta_0 + \frac{q}{T}t, \ r, \ \phi_0 - \frac{2\pi p}{T} t \right),
$$
all having minimal period
\begin{equation}
\label{eqn:period}
T = q\frac{D(r)}{g'(r)} = 2\pi p\frac{D(r)}{f'(r)}
\end{equation}
(in the cases where $f'(r) = p = 0$ or $g'(r) = q = 0$, pick
whichever one of these expressions makes sense).
The torus is \emph{Morse-Bott} if and only if
the function $\rho\mapsto f'(\rho) / g'(\rho)$ (or its reciprocal)
has nonvanishing derivative at $r$.

Likewise, the circle
$$
P := \{ \rho = 0 \} \subset M
$$
is a closed orbit with minimal period $T = |f(0)|$.  For $k \in \NN$,
its $k$--fold cover $P^k$ is degenerate if and only if
$\frac{k f''(0)}{2\pi g''(0)} \in \ZZ$, and otherwise has
$$
\muCZ^{\Phi_0}(P^k) = 
2\left\lfloor - \frac{k f''(0)}{2\pi g''(0)}\right\rfloor + 1,
$$
where $\Phi_0$ is the natural symplectic trivialization of $\xi$ along 
$P$ induced by the coordinates.  
Here $\lfloor x \rfloor$ means the greatest integer $\le x$.
\end{prop}

In terms of the curve $\rho \mapsto (f(\rho),g(\rho))$ in $\RR^2$, 
this says that
the torus $L_r$ is Morse-Bott when the slope of the curve has nonvanishing
derivative at~$r$.  The nondegeneracy and index of any cover of $P$
depend similarly on the slope of this curve as it pushes off from the
$x$--axis at $\rho=0$.  

We've chosen the vector fields $v_1$ and $v_2$ above so that
$d\lambda(v_1,v_2) \equiv 1$, i.e.~they give a symplectic trivialization
of $\xi$ over $M\setminus P$.  Use these now to define an admissible
complex multiplication $J$ by
\begin{equation}
\label{eqn:defofJ}
Jv_1 = \beta(\rho) v_2, \quad\quad\quad J v_2 = -\frac{1}{\beta(\rho)}v_1
\end{equation}
for some smooth function $\beta(\rho)$.  The behavior of $\beta$ near $0$ 
can be chosen to ensure that $J$ is smooth at $\rho = 0$.
Then an $\RR$--invariant almost complex structure $\tilde{J}$ on $\RR\times M$
is defined in the standard way, and we
seek maps $\tilde{u} : (S,j) \to (\RR\times M, \tilde{J})$ 
defined on a Riemann surface $(S,j)$ and satisfying 
$T\tilde{u} \circ j = \tilde{J} \circ T\tilde{u}$.
Choosing conformal coordinates $(s,t)$ on $S$, the map $u$ can be written 
in coordinates as
$u(s,t) = (\theta(s,t),\rho(s,t),\phi(s,t))$, and then 
the Cauchy-Riemann equation becomes
\begin{equation}
\label{eqn:CR}
\begin{aligned}
a_s &= f\theta_t + g\phi_t \quad\quad\quad\quad &
\rho_s &= \frac{1}{\beta} (f'\theta_t + g'\phi_t) \\
a_t &= -f\theta_s - g\phi_s &
\rho_t &= -\frac{1}{\beta} (f'\theta_s + g'\phi_s) \\
\end{aligned}
\end{equation}

\begin{figure}
\begin{center}
\includegraphics{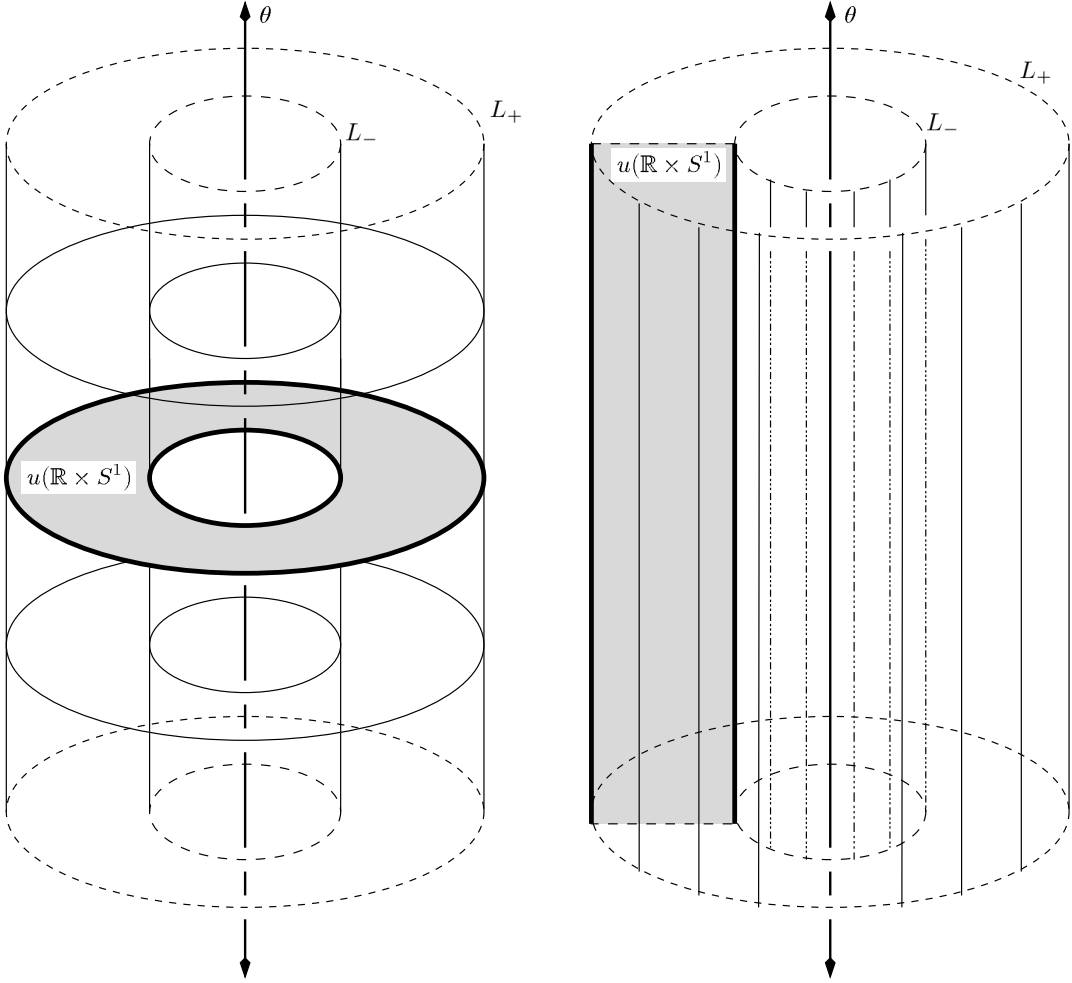}
\end{center}
\caption[Finite energy cylinders spanning concentric Morse-Bott 
tori]{\label{fig:concentric} Concentric tori with homologous
periodic orbits connected by a finite energy cylinder.  On the
left is the case where $g'(\rho_{\pm}) = 0$, so the orbits are
parallel to $\p_\phi$.  On the right, $f'(\rho_{\pm}) = 0$
gives orbits parallel to $\p_\theta$.}
\end{figure}

\begin{figure}
\begin{center}
\includegraphics{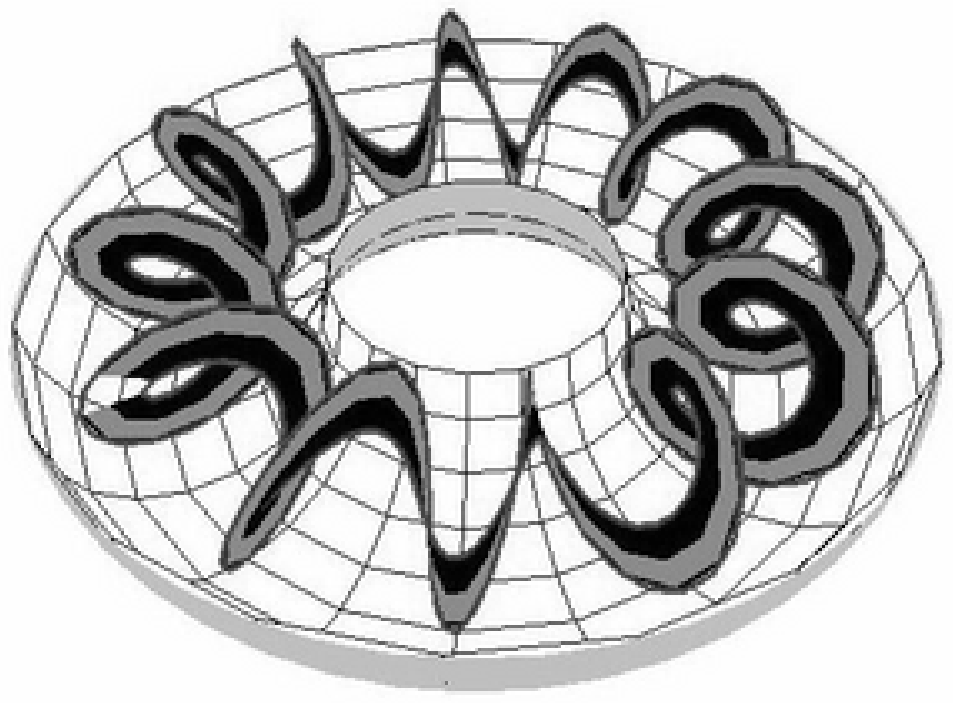}
\end{center}
\caption[Holomorphic cylinder of type $(p,q)$]{\label{fig:cylinderPQ} A 
cylinder of type $(p,q)$ in $S^1 \times \RR^2$ with 
$\rho_+ > \rho_- > 0$.}
\end{figure}

\begin{figure}
\begin{center}
\includegraphics{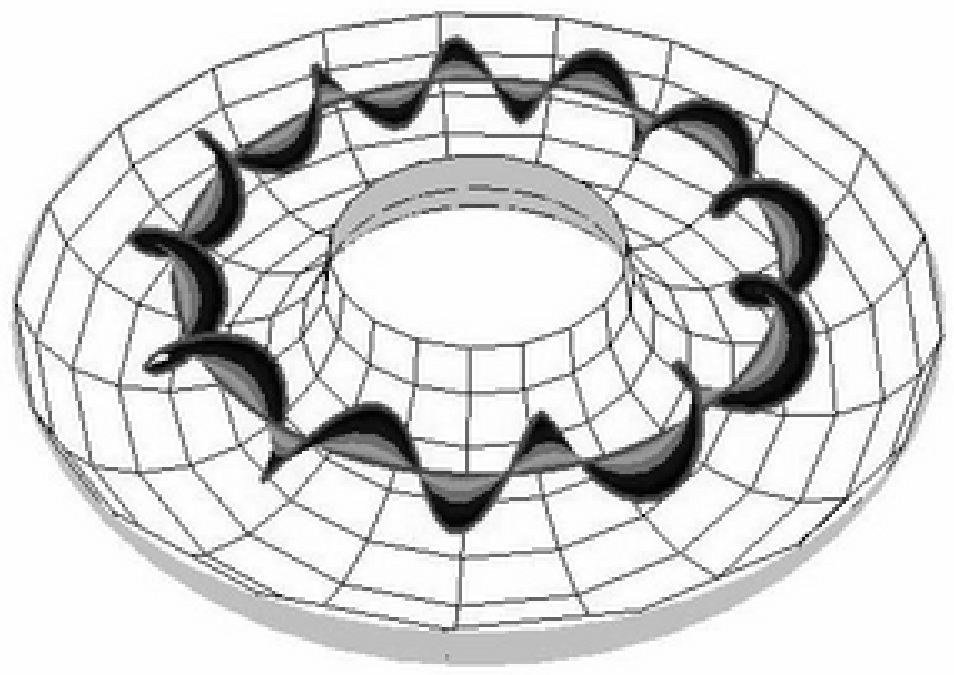}
\end{center}
\caption[Cylinder of type $(p,q)$ with $\rho_- = 0$]{\label{fig:cylinderPQa} 
A cylinder of type 
 $(p,q)$ in $S^1 \times \RR^2$ with $\rho_+ > \rho_- = 0$.}
\end{figure}

Given two concentric tori 
$L_\pm = \{\rho = \rho_{\pm}\}$, each foliated by periodic orbits that 
are homologous (up to a sign) in $H_1(M\setminus P)$, 
we shall now construct a stable finite energy
foliation of the region between them, each leaf being a
cylinder with ends asymptotic to orbits at $L_-$ and $L_+$
respectively (Figure~\ref{fig:concentric}).  
In particular, suppose there are 
two radii $\rho_{\pm}$ with $0 < \rho_- < \rho_+$,
such that
$$
\frac{f'(\rho_\pm)}{2\pi g'(\rho_\pm)} = \frac{p}{q} \in 
        \QQ \cup \{ \infty \}
\quad\text{ and }\quad
\frac{f'(\rho)}{2\pi g'(\rho)} \ne \frac{p}{q} \text{ for }
        \rho \in (\rho_-,\rho_+).
$$
A choice of sign must be made for $p$ and $q$: for reasons that will become
clear shortly, let us choose both so that the quantity
$qf'(\rho) - 2\pi pg'(\rho)$ is positive for $\rho \in (\rho_-,\rho_+)$.
The two tori $L_\pm$ are each foliated by families of periodic orbits, of
the form
$$
x_\pm(t) = \left( \theta_0 + \frac{q_\pm}{T_\pm}t, \rho_\pm,
        \phi_0 - \frac{2\pi p_\pm}{T_\pm}t \right).
$$
Here $p_\pm$ and $q_\pm$ are the same as $p$ and $q$ up to a sign, which
must be chosen so that the periods $T_\pm = 
\frac{q_\pm D(\rho_\pm)}{g'(\rho_\pm)} = 
\frac{2\pi p_\pm D(\rho_\pm)}{f'(\rho_\pm)}$ are positive.  Fixing values
of $\theta_0$ and $\phi_0$, suppose 
$\tilde{u} = (a,u): \RR\times S^1 \to \RR\times M$ is a map of the form
\begin{equation}
\label{eqn:ansatz}
(a(s,t),\theta(s,t),\rho(s,t),\phi(s,t)) = 
(a(s), \theta_0 + qt, \rho(s), \phi_0 - 2\pi pt ).
\end{equation}
Then using \eqref{eqn:CR}, the Cauchy-Riemann equation for $\tilde{u}$
reduces to the pair of ordinary differential equations
\begin{subequations}
\label{eqn:CR_ODE}
\begin{align}
\label{eqn:CR_ODE_rho}
\frac{d\rho}{ds} &= \frac{1}{\beta(\rho)}
        ( qf'(\rho) - 2\pi p g'(\rho) ), \\
\label{eqn:CR_ODE_a}
\frac{da}{ds} &= qf(\rho) - 2\pi p g(\rho).
\end{align}
\end{subequations}
These have unique solutions for any choice of $\rho(0) \in (\rho_-,\rho_+)$
and $a(0) \in \RR$.  Notice that due to our sign convention for $p$ and $q$,
the right hand side of \eqref{eqn:CR_ODE_rho} is always positive,
thus $\lim_{s\rightarrow\pm\infty} \rho(s) = \rho_\pm$, and
we see that $u(s,\cdot)$ converges in $C^\infty$ to parametrizations
of the orbits $x_\pm$ as $s \to \pm\infty$.  It follows then from
\eqref{eqn:energy} and Stokes' theorem that $\tilde{u}$ has finite energy
$E(\tilde{u}) \le T_+ + T_-$.
We shall refer to this solution as a 
\emph{cylinder of type $(p,q)$}.  An example is shown in
Figure~\ref{fig:cylinderPQ}.\footnote{Thanks to Joel Fish for providing
Figures~\ref{fig:cylinderPQ} and~\ref{fig:cylinderPQa}.}

It is clear from \eqref{eqn:CR_ODE_a} that
$a$ is a proper function with asymptotically linear growth to
$\pm\infty$, as the condition $D(\rho) > 0$ guarantees
that $\lim_{s\rightarrow\pm\infty} a'(s) = qf(\rho_\pm) - 2\pi p g(\rho_\pm)$
cannot be zero.  This expression determines the sign of the puncture at
$s = \pm\infty$ as $\pm \sign( qf(\rho_\pm) - 2\pi p g(\rho_\pm))$.
To put this in a more revealing form, write $f_\pm := f(\rho_\pm)$,
$f'_\pm := f'(\rho_\pm)$ etc., and observe that by assumption there is
a nonzero number 
$$
c_\pm = \frac{2\pi p}{f'_\pm} = \frac{q}{g'_\pm}.
$$
Then the expression above for the sign 
becomes $\pm \sign[ c_\pm (f_\pm g'_\pm - 
f'_\pm g_\pm) ] = \pm \sign(c_\pm)$ since $D(\rho_\pm)$ is positive.  
Now if both tori $L_\pm$ satisfy the 
Morse-Bott condition, then 
$0 \ne f'_\pm g''_\pm - f''_\pm g'_\pm = -\frac{1}{c_\pm} (q f''_\pm - 2\pi p
g''_\pm)$, and our sign convention for $p$ and $q$ implies
$\sign( q f''_\pm - 2\pi p g''_\pm) = \mp 1$, thus 
$\sign(f'_\pm g''_\pm - f''_\pm g'_\pm) 
= - \sign(c_\pm) \sign(q f''_\pm - 2\pi p g''_\pm)
= \pm \sign(c_\pm)$,
and we have
\begin{equation}
\label{eqn:signPuncture}
\text{sign of puncture at $L_\pm$} =
\sign(f'_\pm g''_\pm - f''_\pm g'_\pm).
\end{equation}
This means that in the Morse-Bott case, the sign of a puncture 
approaching $L_r$ is positive if and only if the counterclockwise 
trajectory $\rho \mapsto (f(\rho),g(\rho))$ is
accelerating \emph{inward} at $\rho = r$, and negative if it accelerates
outward.

The equations~\eqref{eqn:CR_ODE} can be thought of as defining a direction 
field in the subset $(\rho_-,\rho_+) \times \RR$ of the $\rho a$--plane,
which integrates to a one-dimensional foliation.  
Since \eqref{eqn:CR_ODE_a} defines $a(s)$ only up to a constant,
this foliation is invariant under the natural $\RR$--action on the 
$a$--coordinate.
Meanwhile the set of trajectories $t\mapsto 
(\theta_0 + qt, \phi_0 - 2\pi pt) \in S^1 \times \RR / 2\pi\ZZ$ for
all choices of $\theta_0$ and $\phi_0$ defines another one-dimensional
foliation.  Putting these together as in \eqref{eqn:ansatz}
creates a two-dimensional foliation of the region
$\{ (a,\theta,\rho,\phi) \in \RR\times M\ | \ \rho\in (\rho_-,\rho_+) \}$
by $\tilde{J}$--holomorphic cylinders with uniformly bounded energy, 
and it projects to a one-dimensional
foliation of $\{ \rho \in (\rho_-,\rho_+) \} \subset M$.
We may assume without loss of generality that $f$ and $g$ are chosen
so that both tori $L_{\rho_\pm}$ are Morse-Bott.  Then
using the frame $(v_1,v_2)$ to trivialize $\xi$ over this region,
it follows from Lemma~\ref{lemma:MBparity} and \eqref{eqn:FredholmIndexBasic}
that each leaf $\tilde{u}$
has $\#\Gamma_0 = 0$ and $\ind(\tilde{u}) = 2$, so Theorem~\ref{thm:IFT}
implies that the foliation is stable.

We can extend this foliation to $\rho=\rho_\pm$
by adding the cylinders over periodic orbits at $L_\pm$.
Moreover, if there exists a radius $\rho_0 \in (0,\rho_-)$ such that
$\rho_0$ and $\rho_-$ satisfy the same 
conditions as $\rho_-$ and $\rho_+$, then we can repeat this construction for 
$\rho\in(\rho_0,\rho_-)$
and thus extend the foliation to the region $\rho \in [\rho_0,\rho_+]$.

It remains to extend the foliation further toward the center in 
the case where there is no $\rho < \rho_-$ with 
$f'(\rho) / 2\pi g'(\rho) = p/q$.  To that end, let us redefine our notation
with $\rho_- = 0$ and $L_+ = \{\rho = \rho_+\}$; choose $\rho_+ > 0$ 
so that
$$
\frac{f'(\rho_+)}{2\pi g'(\rho_+)} = \frac{p}{q} \in \QQ \cup \{\infty\}
\quad\text{ and }\quad
\frac{f'(\rho)}{2\pi g'(\rho)} \ne \frac{p}{q} \text{ for } 
 \rho \in (0,\rho_+).
$$
Choose the signs of $p$ and $q$ so that $qf' - 2\pi p g' > 0$ for
$\rho \in (0,\rho_+)$, and consider once more the family of 
$\tilde{J}$--holomorphic cylinders defined by
$$
\tilde{u} = (a,u):
\RR\times S^1 \to \RR\times M\ :\  (s,t) \mapsto 
(a(s),\theta_0 + qt,\rho(s),\phi_0-2\pi pt),
$$
where $\rho(s)$ and $a(s)$ satisfy the ODEs~\eqref{eqn:CR_ODE} with
$\rho(0) \in (0,\rho_+)$.
Once again $u(s,\cdot)$ converges in $C^\infty$
as $s\rightarrow\infty$ to some parametrization of a simply covered orbit
$P_+ \subset L_+$, and \eqref{eqn:signPuncture} gives 
the sign of this puncture~as
$$
\sigma_+ := \sign( f_+' g_+'' - f_+'' g_+' ).
$$
Define $F(\rho)$ to be the right hand side of \eqref{eqn:CR_ODE_rho}.
The requirement that $J$ be smooth at $\rho=0$ implies
that $\beta(\rho)$ is bounded away from zero as $\rho\to 0$, thus
$\lim_{\rho\to 0} F(\rho) = 0$, and we
conclude that $\rho(s) \to 0$ as $s\to -\infty$.

We must now distinguish between two cases in order to understand fully the
behavior as $s\to -\infty$.  
If $q \ne 0$, $u(s,\cdot)$ converges to the $|q|$--fold cover of $P$,
and the sign of the puncture at $-\infty$ is
$$
\sigma_- := -\sign(q) \cdot \sign[f(0)] =
-\sign(q) \cdot \sign[g''(0)],
$$
where we're using the fact that $f(0)g''(0) > 0$.  We can put this in
a more geometrically revealing form analogous to \eqref{eqn:signPuncture}:
observe first that if $P^{|q|}$ is nondegenerate, Prop.~\ref{prop:symOrbits}
implies $q f_-'' / 2\pi g_-'' \not\in \ZZ$ and thus
$f_-'' / 2\pi g_-'' \ne f_+' / 2\pi g_+' = p/q$.  Meanwhile our sign convention
$q f' - 2\pi p g' > 0$ for $\rho \in (0,\rho_+)$ implies
$$
q f_-'' - 2\pi p g_-'' > 0.
$$
This together with the above expression for $\sigma_-$ yields
$$
\frac{1}{2\pi} \sigma_- \left(\frac{f_-''}{g_-''} - \frac{f_+'}{g_+'}\right) =
\sigma_- \left(\frac{f_-''}{2\pi g_-''} - \frac{p}{q}\right) < 0,
$$
hence
\begin{equation}
\label{eqn:sigmaMinus}
\sigma_- = \sign\left( \frac{f_+'}{g_+'} - \frac{f_-''}{g_-''} \right).
\end{equation}
Thus when $P^{|q|}$ is nondegenerate, $\sigma_-$ depends on whether the 
slope of the trajectory
$\rho \mapsto (f(\rho),g(\rho))$ at $\rho=0$ is greater than or less
than the slope at $\rho = \rho_+$.

\begin{figure}
\begin{center}
\includegraphics{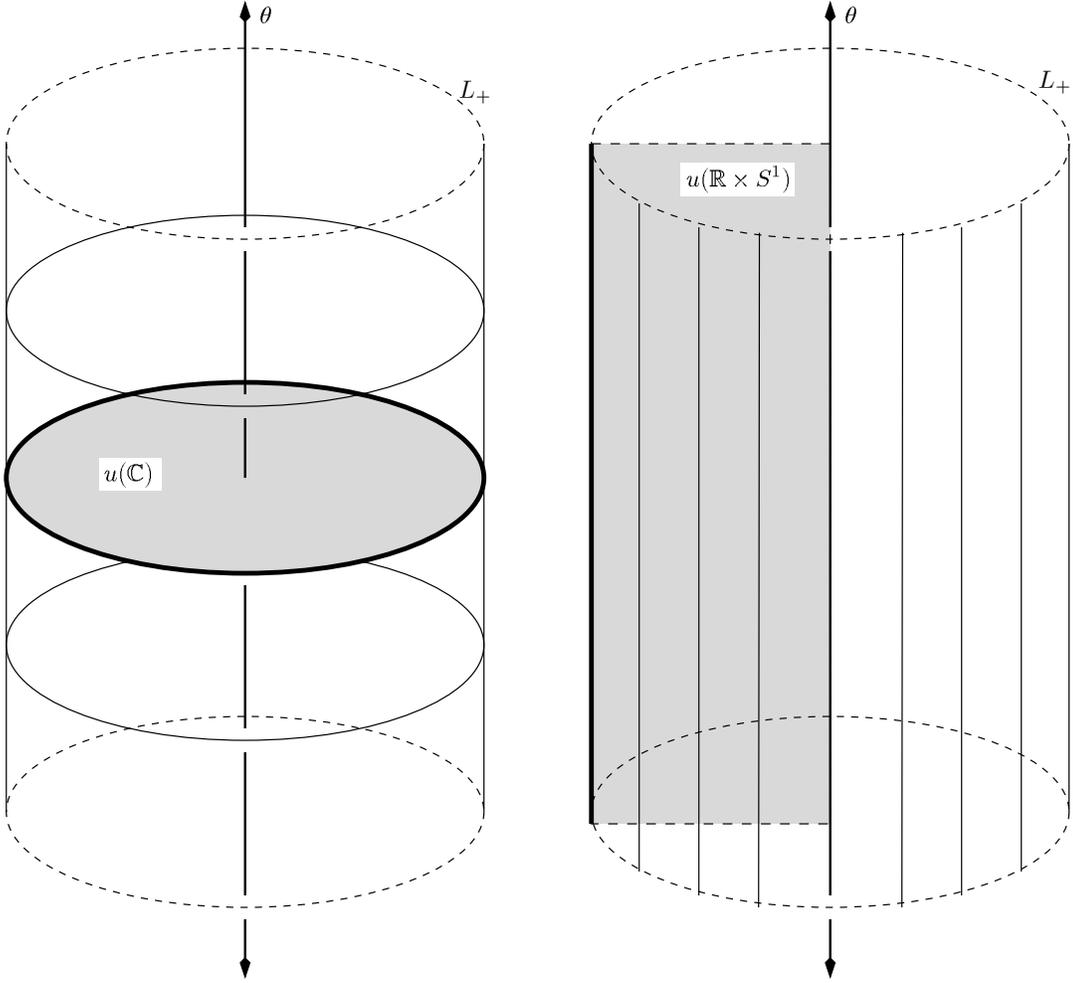}
\end{center}
\caption[Finite energy planes and cylinders inside a Morse-Bott
 torus]{\label{fig:center} Holomorphic curves inside the innermost
torus.  If orbits on $L_+$ have nontrivial $\p_\theta$ component
(right), we get finite energy cylinders with a puncture asymptotic
to the central axis; else that puncture is removable (left) and we
get a finite energy plane.}
\end{figure}

We now have a foliation of the region
$\rho \in (0,\rho_+)$ by an
$\RR$--invariant family of finite energy cylinders, each convergent
to $P^{|q|}$ at one end and a simply covered orbit $P_+ \subset L_+$ 
at the other.  These together with the orbit cylinder over $P$ form a
finite energy foliation in the region $\{\rho < \rho_+\}$.
Figure~\ref{fig:center}, right, shows an example with $(p,q) = (0,1)$.
An example with $p$ and $q$ both nonzero is shown in 
Figure~\ref{fig:cylinderPQa}.

Stability for these cylinders is a somewhat more subtle question than before.
Assume $P^{|q|}$ is nondegenerate and $L_+$
is Morse-Bott.  Then since $P^{|q|}$ has odd Conley-Zehnder
index, Lemma~\ref{lemma:MBparity} and Theorem~\ref{thm:windpi} imply
that each solution $\tilde{u}$ above has $\ind(\tilde{u}) \ge 2$;
in general however, this inequality can be strict.  We claim
that the functions $f$ and $g$ can always be adjusted near~$0$ 
so that $\ind(\tilde{u}) = 2$; in this case Theorem~\ref{thm:IFT} will
imply that the corresponding foliation is stable.
Let $\Phi_0$ be the symplectic trivialization of
$\xi$ along $P$ defined by the Cartesian coordinates, and extend this 
over $\{ \rho \le \rho_+ \}$.  Then accounting for the orientation of
$\xi$ determined by $\lambda$, we have
$$
\wind_{P_+}^{\Phi_0}(v_1) = -\sigma_+ \cdot \sign[f(0)] \cdot p,
$$
and hence
\begin{equation}
\label{eqn:CZ}
\muCZ^{\Phi_0\mp}(P_+) = \sigma_+ \left( 1 - 2\sign[f(0)] \cdot p\right).
\end{equation}
In order to write $\muCZ^{\Phi_0}(P^{|q|})$ in a convenient form, we compute
\begin{multline*}
\left\lfloor -\frac{|q| f_-''}{2\pi g_-''} \right\rfloor =
  \left\lfloor -\sign(q) \left( \frac{q f_-'' - 2\pi p g_-''}{2\pi g_-''} + p
 \right) \right\rfloor \\
= \left\lfloor \sigma_- \frac{q f_-'' - 2\pi p g_-''}{2\pi |g_-''|} 
  \right\rfloor - \sign(q)\cdot p.
\end{multline*}
Then using the index formula in Prop.~\ref{prop:symOrbits},
\begin{equation*}
\begin{split}
\ind(\tilde{u}) &= \sigma_+ \cdot \muCZ^{\Phi_0\mp}(P_+) + 
  \sigma_- \cdot \muCZ^{\Phi_0}(P^{|q|}) \\
&= 1 - 2\sign[f(0)] \cdot p + \sigma_- \left( 2 \left\lfloor
 - \frac{|q| f_-''}{2\pi g_-''} \right\rfloor + 1\right) \\
&= 1 + \sigma_- - (2\sign[f(0)] \cdot p) - (2\sigma_- \cdot \sign(q)\cdot p) \\
&\qquad\qquad + 
 2\sigma_- \left\lfloor \sigma_- \frac{q f_-'' - 2\pi p g_-''}{2\pi |g_-''|}
   \right\rfloor \\
&= 1 + \sigma_- 
  + 2\sigma_- \left\lfloor \sigma_- \frac{q f_-'' - 2\pi p g_-''}{2\pi |g_-''|}
   \right\rfloor.
\end{split}
\end{equation*}
If $\sigma_- = 1$ this gives
$$
\ind(\tilde{u}) = 2 + 2\left\lfloor \frac{q f_-'' - 2\pi p g_-''}{2\pi |g_-''|}
 \right\rfloor,
$$
or if $\sigma_+ = -1$,
$$
\ind(\tilde{u}) = 2\left\lceil \frac{q f_-'' - 2\pi p g_-''}{2\pi |g_-''|}
 \right\rceil,
$$
where $\lceil x \rceil$ denotes the smallest integer $\ge x$.
Recalling $q f_-'' - 2\pi p g_-'' > 0$, we see indeed that in both cases
$\ind(\tilde{u}) \ge 2$, with equality
if and only if $q f_-'' - 2\pi p g_-''$ is sufficiently small.  This can
always be achieved by adjusting the slope of the trajectory
$\rho \mapsto (f(\rho),g(\rho))$ for $\rho$ near~$0$, without creating
any new points at which $q f'(\rho) - 2\pi p g'(\rho) = 0$.  The key point
is to make the slopes of the trajectory $\rho \mapsto (f(\rho),g(\rho))$
at $\rho=0$ and $\rho=\rho_+$ as close as possible.

If on the other hand $q = 0$, we have $p = -\sign(g_-'') =
-\sign[f(0)] = \pm 1$ and
$\lim_{s\to -\infty}u(s,t) = 
(\theta_0,0) \in P \subset S^1\times\RR^2$.  In fact, since
$$
\lim_{\rho\to 0}F'(\rho) = - 
\frac{2\pi p g''(0)}{\lim_{\rho\to 0}\beta(\rho)} \ne 0,
$$
one can easily show
that $\rho(s)$ converges exponentially fast to 0, and plugging this 
behavior into the equation $\rho' = F(\rho)$, so does its derivative.
We now claim that $a(s)$ is bounded at $-\infty$.  For this 
it suffices to prove that the integral
$$
\int_{-\infty}^0 \frac{da}{ds}\ ds = -2\pi p \int_{-\infty}^0 g(\rho(s))\ ds
$$
converges.  We know $\rho'(s)$ satisfies a bound of the form
$|\rho'(s)| \le Me^{\lambda s}$ with $\lambda > 0$.  Since $g'$ is 
continuous and $\rho$ stays within a bounded interval for all $s$, we have
\begin{multline*}
|g(\rho(s))| = \left| \int_{-\infty}^s \frac{d}{d\sigma} 
        g(\rho(\sigma)) \ d\sigma \right|
 \le \int_{-\infty}^s |g'(\rho(\sigma))|\  |\rho'(\sigma)| \ d\sigma \\
 \le M_1 \int_{-\infty}^s e^{\lambda\sigma}\ d\sigma
 = M_2 e^{\lambda s}
\end{multline*}
for some constant $M_2 > 0$.  Then $\int_{-\infty}^0 |g(\rho(s))|\ ds
< \infty$ and the claim follows.  It's clear now that $\tilde{u}$ has
finite area as $s \to -\infty$, thus Gromov's removable singularity
theorem implies that $\tilde{u}$ can be extended smoothly
to a finite energy plane 
$\tilde{v} = (b,v): \CC\to \RR\times M$ with 
$\tilde{v}(e^{2\pi (s+it)}) = \tilde{u}(s,t)$ and
$v(0) = (\theta_0,0) \in S^1\times\RR^2$.  The set of all such planes
forms a finite energy foliation in the region
$\{ \rho < \rho_+ \}$.  Each is positively
asymptotic to a simply covered orbit $P_+ \subset L_+$,
and transverse to the central orbit $P$ (Figure~\ref{fig:center}, left).
From \eqref{eqn:CZ}, we find
$$
\ind(\tilde{v}) = \muCZ^{\Phi_0-}(P_+) - \chi(\CC) = (1 - (-2)) - 1 = 2,
$$
so the foliation is stable.

We now apply these constructions to a contact manifold $(M,\xi_K)$ obtained
from $(S^3,\xi)$ by a twist surgery along a knot $K$.  Let 
$N \subset S^3$ be the corresponding solid torus neighborhood of $K$,
identified with $S^1 \times \overline{B^2_\epsilon(0)}$, and denote by
$N' = S^1 \times \overline{B^2_\epsilon(0)} 
\subset M$ the solid torus that replaces
it after surgery; thus $M \setminus N' = S^3 \setminus N$.
On $N'$, $\xi_K$ is the kernel of $\lambda_K = f_K(\rho')\ d\theta' +
g_K(\rho')\ d\phi'$, which for $\rho \in [\delta,\epsilon]$ is the
pull back of $\lambda_1 = f_1(\rho)\ d\theta + g_1(\rho)\ d\phi$ on $N$
via the gluing map.  Let $\rho_1 \in (\delta,\epsilon)$ be the largest
radius for which $g_1'(\rho_1) = 0$, so the Reeb orbits on the torus
$L_{\rho_1} := \{ \rho = \rho_1 \}$ are negatively oriented meridians on $N$.
They are generally not meridians on $N'$, but will
represent the homology class $q \lambda' - n \mu' \in
H_1(L_{\rho_1})$, where $\lambda'$ and $\mu'$ are the standard longitude and
meridian respectively for $N'$, and $q$ and $n$ are the integers appearing
in the matrix \eqref{eqn:surgeryMatrix}.  There is then a finite set of
radii
$$
\rho_1 > \rho_2 > \ldots > \rho_s > 0
$$
for which the Reeb orbits on $L_{\rho_j}$ represent classes
$\pm (q \lambda' - n \mu') \in H_1(L_{\rho_j})$, and we can foliate the regions
in between each of these concentric tori by cylinders of type $(n,q)$.
For the region $\{ \rho \in (0,\rho_s) \}$, we obtain planes if $q = 0$,
otherwise cylinders asymptotic to the $|q|$--fold cover of the orbit
$S^1 \times \{0\} \subset N'$.

By Prop.~\ref{prop:important}, we're free to assume without loss of
generality that $q$ is either $0$ or $\pm 1$, in which case all the orbits in 
this foliation are simply covered.  By conditions~(3) and~(4) in the
proposition, we can also assume the tori $L_{\rho_j}$ are all Morse-Bott,
the orbit $S^1 \times \{0\}$ is nondegenerate, all punctures
approaching these orbits are positive and the cylinders in the innermost
region have index~$2$.

Before turning to more global considerations, we note another useful local
result, which allows a kind of ``analytic continuation'' for some foliations.

\begin{prop}
\label{prop:continuation}
Let $M = S^1 \times B^2_\epsilon(0)$ with contact form $\lambda =
f(\rho)\ d\theta + g(\rho)\ d\phi$ and $J v_1 = \beta(\rho) v_2$ defining
the almost complex structure $\tilde{J}$ on $\RR\times M$.  Now for some
$\delta \in (0,\epsilon)$, choose new smooth functions 
$f_1$, $g_1$ and $\beta_1$
which match $f$, $g$ and $\beta$ for $\rho \in [\delta,\epsilon)$, such that
$f_1 g_1' - f_1' g_1 > 0$ for $\rho \in (0,\epsilon)$ and $\beta_1(0) > 0$.  
These define new data $\lambda_1$, $J_1$ and $\tilde{J}_1$, which are 
smooth on $M \setminus (S^1 \times \{0\})$ but not necessarily at 
$S^1 \times \{0\}$.
Suppose there is a neighborhood $\uU$ of $S^1 \times \overline{B^2_\delta(0)}$
in $M$ and a family of $\tilde{J}$--holomorphic finite energy
half-cylinders $\tilde{u} = (a,u) : (-\infty,0]\times S^1 \to \RR\times M$,
all asymptotic to a particular cover of the orbit $P := S^1 \times \{0\}$,
and defining a finite energy foliation on $\uU$.
Assume moreover that either of the following is true:
\begin{enumerate}
\setlength{\itemsep}{0in}
\item
$f'g'' - f''g' \equiv f_1'g_1'' - f_1''g_1' \equiv 0$.
\item
$u(s,t) = (\theta(s,t),\rho(s,t),\phi(s,t))$ satisfies
$\theta_s \phi_t - \theta_t \phi_s \equiv 0$.
\end{enumerate}
Then for each $\tilde{u}$ there is a unique $\tilde{J}_1$--holomorphic 
half-cylinder $\tilde{u}_1(s,t)$ which matches $\tilde{u}$ on some annulus
of the form $[-s_0,0] \times S^1$.  If additionally $(f_1,g_1)$ are
$C^1$--close to $(f,g)$ then the new maps $\tilde{u}_1$ are proper and 
asymptotic to the same orbit as $\tilde{u}$, and they form an 
$\RR$--invariant foliation of $\RR \times (\uU \setminus P)$, projecting 
to a foliation of $\uU \setminus P$.
\end{prop}
\begin{proof}
Writing the given curves as 
$\tilde{u}(s,t) = (a(s,t),\theta(s,t),\rho(s,t),\phi(s,t))$,
there are constants $\theta_0$, $\phi_0$ and fixed integers $p$ and $q$
such that
$$
\theta(s,t) \to \theta_0 + qt,
\qquad
\phi(s,t) \to \phi_0 - 2\pi p t
\qquad
\text{ as $s \to -\infty$},
$$
and the functions $(a,\theta,\rho,\phi)$ satisfy \eqref{eqn:CR}.
Combining this with the expressions
$a_{st} - a_{ts} = 0$ and $\rho_{st} - \rho_{ts} = 0$ implies
\begin{equation}
\label{eqn:intega}
f \Delta\theta + g \Delta\phi = 0,
\end{equation}
\begin{equation}
\label{eqn:integrho}
f' \Delta\theta + g' \Delta\phi - \frac{1}{\beta}
 (f'g'' - f''g') (\theta_s\phi_t - \theta_t\phi_s) = 0,
\end{equation}
where $\Delta := \p_{ss} + \p_{tt}$, and $f$, $g$ and $\beta$ are understood 
to depend on $\rho(s,t)$.  We now seek a map of the form
$$
\tilde{u}_1(s,t) = (a_1(s,t), \theta(s,t), \rho_1(s,t), \phi(s,t))
$$
such that $a_1(s,t) = a(s,t)$ and $\rho_1(s,t) = \rho(s,t)$ for
$s \ge -s_0$, and solving the corresponding Cauchy-Riemann equations
with respect to $f_1$, $g_1$ and $\beta_1$.  Since $\theta(s,t)$ and 
$\phi(s,t)$ are now fixed functions, the new equations for $\rho_1$
in \eqref{eqn:CR} can be interpreted as saying that the graph
$\Gamma_{\rho_1} := \{ (s,t,\rho_1(s,t)) \}$ is tangent to a certain
$2$--plane distribution in $(-\infty,0] \times S^1 \times \RR$.  This
distribution turns out to be integrable if and only if
\begin{equation}
\label{eqn:integ1rho}
{f_1}' \Delta\theta + {g_1}' \Delta\phi 
 - \frac{1}{{\beta_1}}
 ({f_1}'{g_1}'' - {f_1}''{g_1}') (\theta_s\phi_t - \theta_t\phi_s) 
 \equiv 0,
\end{equation}
where the expression is to be understood
as a function of three independent variables $(s,t,\rho) \in 
(-\infty,0]\times S^1 \times \RR$.  If it
vanishes identically then solutions $\rho_1(s,t)$ exist locally.
Assume this for the moment: then choosing $s_0 \in (-\infty,0)$ such that 
$\rho(s,t) \ge \delta$ for all $s \ge s_0$, there is a solution
${\rho_1}(s,t)$ on $(-s_1,0] \times S^1$ for some
$-s_1 < -s_0$, with ${\rho_1}(s,t) = \rho(s,t)$ for $s \ge -s_0$.  
For topological reasons, the continued solution is automatically $1$--periodic 
in $t$.  Then for fixed $t$, the function
$s \mapsto \rho_1(s,t)$ satisfies the ODE
$$
\frac{d{\rho_1}}{ds} = \frac{1}{{\beta_1}({\rho_1})} \left( 
 {f_1}'({\rho_1}) \theta_t + {g_1}'({\rho_1}) \phi_t \right),
$$
and we see that the 
solution ${\rho_1}(s,t)$ extends to $(-\infty,0] \times S^1$
with $$
\lim_{s\to -\infty} \rho_1(s,t) = \rho_0,
$$
where $\rho_0 \ge 0$ 
is the largest radius at which $f_1'(\rho_0) / 2\pi g_1'(\rho_0) = p/q$,
or zero if there is no such radius.  The latter is necessarily the case,
in particular, if $(f_1,g_1)$ is $C^1$--close to $(f,g)$, because the
same argument for $\tilde{u}$ shows that $(f,g)$ cannot admit any radius
at which this relation is satisfied.

The remaining two equations in \eqref{eqn:CR} specify the gradient of 
$a_1(s,t)$ in terms of known functions, so solutions exist locally if and 
only if this gradient is curl-free, which in this case means
\begin{equation}
\label{eqn:integ1a}
{f_1} \Delta\theta + {g_1} \Delta\phi \equiv 0
\end{equation}
for all $(s,t) \in (-\infty,0]\times S^1$ and $\rho = \rho_1(s,t)$.
There is then a unique solution on 
$(-\infty,0]\times S^1$ with $a_1(s,t) = a(s,t)$ for all $s \ge s_0$,
and another ODE argument establishes that in the case $\rho_0 = 0$,
$a(s,t)$ blows up linearly as $s \to -\infty$.

We claim that the integrability conditions are satisfied whenever either
of the two additional assumptions in the Proposition are met.  Indeed, if
$f'g'' - f''g' \equiv {f_1}'{g_1}'' - {f_1}''{g_1}' \equiv 0$,
then~\eqref{eqn:intega} and~\eqref{eqn:integrho} give
$f \Delta\theta + g \Delta\phi = f' \Delta\theta + g' \Delta\phi = 0$,
and since the contact condition requires $(f,g)$ and $(f',g')$ to be
linearly independent in $\RR^2$ for all $\rho$, we conclude that both
$\theta(s,t)$ and $\phi(s,t)$ are harmonic.  Thus \eqref{eqn:integ1rho}
and \eqref{eqn:integ1a} are satisfied for all $(s,t,\rho)$.  In the
other case, $\theta_s\phi_t - \theta_t\phi_s \equiv 0$ together with
\eqref{eqn:intega} and \eqref{eqn:integrho} implies again that
$\theta$ and $\phi$ are harmonic, so the same argument applies.
\end{proof}

\newpage

\subsection{Surgery on a holomorphic open book}
\label{subsec:surgery}

\subsubsection{$\fF_0 \to \fF_1$: Stabilizing an open book decomposition}
The global construction begins with a stable foliation of open book type
on the tight $3$--sphere.  Such foliations follow from a general existence
result in \cite{HWZ:tight3sphere},
but for our purposes, we can produce one using much less machinery.

Define $\lambda_0$ on $S^3$ as in 
\eqref{eqn:standardS3}: then the Reeb vector field $X_0$ generates
the Hopf fibration.  At each $z \in S^3$, $\xi_0 = 
\ker\lambda_0$ is the unique complex line in $T_z S^3 \subset \CC^2$,
which therefore admits a natural complex multiplication $i \in 
\Gamma(\End(\xi))$.  Let $\tilde{J}_0$ be the $\RR$--invariant almost complex
structure on $\RR\times S^3$ associated to $\lambda_0$ and $i$.  Then the
diffeomorphism
$$
\Phi : (\RR\times S^3,\tilde{J}_0) \to (\CC^2 \setminus \{0\},i) :
(a,m) \mapsto e^{2a} m
$$
is biholomorphic.  For each
$\zeta \in \CC\setminus \{0\}$, we now define a $\tilde{J}_0$--holomorphic plane
$$
\tilde{u}_\zeta = (a_\zeta,u_\zeta) : \CC \to \RR\times S^3 :
z \mapsto \Phi^{-1}(z,\zeta),
$$
and for $\zeta = 0$, a cylinder (i.e.~punctured plane)
$$
\tilde{u}_0 = (a_0, u_0) : \CC\setminus\{0\} \to \RR\times S^3 :
z \mapsto \Phi^{-1}(z,\zeta).
$$
The latter is in fact the trivial cylinder over the Hopf circle
$$
P_\infty := \{ ( e^{2\pi i \theta},0)\ |\ \theta \in S^1 \},
$$
and the collection of planes $\{\tilde{u}_\zeta\}_{\zeta \in \CC\setminus\{0\}}$ 
is an $\RR$--invariant $2$--parameter family of embedded, pairwise disjoint 
finite energy planes asymptotic to $P_\infty$.  Altogether these define a
finite energy foliation $\fF_0$ on $(S^3,\lambda_0,i)$.  The projection to
$S^3$ is a planar open book decomposition with one binding orbit.

The foliation $\fF_0$ is \emph{not} stable, because the degeneracy of 
$P_\infty$ gives the planes index~$4$.
We can fix this with a small change to $\lambda_0$ near $P_\infty$,
using Prop.~\ref{prop:continuation}.  Indeed, pick any $R \le 1 / \sqrt{2\pi}$
and identify a neighborhood of $P_\infty$ with $S^1 \times B^2_R(0)$ via
the embedding
\begin{equation}
\label{eqn:goodCoords}
\psi : S^1 \times B^2_R(0) \hookrightarrow S^3 : (\theta,\rho,\phi) \mapsto
e^{2\pi i\theta} \left( \sqrt{1 - 2\pi\rho^2},e^{i\phi} \sqrt{2\pi}\rho\right).
\end{equation}
Then $\psi(S^1\times\{0\}) = P_\infty$ and 
$\psi^*\lambda_0 = \pi (d\theta + \rho^2 d\phi) = f(\rho)\ d\theta +
g(\rho)\ d\phi$, where
$f(\rho) = \pi$ and $g(\rho) = \pi\rho^2$.  Defining the vector fields
$v_1$ and $v_2$ as in \eqref{eqn:v1v2}, the complex multiplication is now
specified by $i v_1 = \beta(\rho) v_2$, where
$$
\beta(\rho) = \frac{2\pi}{1 - 2\pi\rho^2}.
$$
For $\zeta = r e^{i\phi_0} \in \CC\setminus\{0\}$, we can express the 
asymptotic behavior of the holomorphic plane 
$\tilde{u}_\zeta(z) = (a_\zeta(z), u_\zeta(z)) = 
 \left( \frac{1}{2} \ln |(z,\zeta)|, \frac{(z,\zeta)}{|(z,\zeta)|} \right)$
in these coordinates by
\begin{multline*}
(a(s,t),\theta(s,t),\rho(s,t),\phi(s,t)) :=
\left( a\left(e^{-2\pi (s+it)}\right), \psi^{-1}\circ 
  u_\zeta\left(e^{-2\pi (s+it)}\right) \right) \\
= \left( \frac{1}{4} \ln(e^{-4\pi s} + r^2), -t, 
 \frac{r}{\sqrt{2\pi (e^{-4\pi s} + r^2)}}, \phi_0 + 2\pi t \right),
\end{multline*}
with $(s,t) \in (-\infty,s_0] \times S^1$ for $s_0$ sufficiently close
to~$-\infty$.  Observe now that $\theta_s \phi_t - \theta_t \phi_s \equiv 0$,
thus by Prop.~\ref{prop:continuation}, any $C^1$--small change in
$f$ and $g$ for $\rho$ near~$0$ admits a new foliation, which is identical
to $\fF$ outside some neighborhood of $P_\infty$.  In particular, 
pick $\delta \in (0,R)$ and define
$$
\lambda_1 = f_1(\rho)\ d\theta + g_1(\rho)\ d\phi =
h(\rho) \cdot ( \pi\ d\theta + \pi\rho^2\ d\phi)
$$
for some function $h$ that satisfies $h(\rho) = 1$ for $\rho \ge \delta$
and is $C^1$--close to this on $[0,R)$,
and such that $h''(0)$ is small but positive.
Then a calculation using Prop.~\ref{prop:symOrbits} shows that for the new
contact form, $P_\infty$ is a nondegenerate orbit with $\muCZ(P_\infty) = 3$.  
The new family of planes asymptotic
to $P_\infty$ then have index~$2$ and form a stable foliation $\fF_1$
on $(S^3,\lambda_1,i)$.

\subsubsection{$\fF_1 \to \fF_2$: Fixing $\lambda$ and $J$ near a link}
\label{subsubsec:simplify}

Next, introduce a positively transverse link $K = K_1 \cup \ldots \cup K_n
\subset S^3$.  By Lemma~\ref{lemma:almostHopf}, there are smooth families
$\gamma_j^\tau : S^1 \to M$ for $\tau\in [0,1]$ such that
$\gamma_j^1(S^1) = K_j$,
$\gamma_j^0(t) = (0,e^{2\pi i k_j t})$ for some $k_j\in\NN$,
and for each fixed $\tau\in (0,1]$, the maps
$\gamma_1^\tau,\ldots,\gamma_n^\tau : S^1\to M$ are mutually non-intersecting
embeddings transverse to $\xi$.  Denote $K^\tau_j = \gamma_j^\tau(S^1)$
and $K^\tau = K^\tau_1 \cup \ldots \cup K^\tau_N$ for $\tau\in(0,1]$.
\begin{lemma}
\label{lemma:stillTransverse}
For $\tau > 0$ sufficiently small,
there is a contact form $\lambda_2$ with $\ker\lambda_2 = \xi_0$
and the following properties:
\begin{enumerate}
\item[(i)] $\lambda_2$ is $C^1$--close to $\lambda$, and differs from
  $\lambda$ only in an arbitrarily small neighborhood of $K^\tau$,
\item[(ii)] Each of the knots $K_j^\tau$ has a tubular neighborhood $N_j 
 \cong S^1 \times B^2_\epsilon(0)$ with coordinates $(\theta,\rho,\phi)$ 
  in which $K_j^\tau = \{ \rho = 0 \}$ and
  $\lambda_2 = c (d\theta + \rho^2\ d\phi)$ for some constant $c > 0$.
\end{enumerate}
\end{lemma}
For a complete proof, we refer to \cite{Wendl:thesis}*{Prop.~5.1.3}.
The main idea is as follows: observe first that a neighborhood of 
$P_0 := \{ (0,e^{2\pi i \theta})\ |\ \theta \in S^1 \}$ admits coordinates
in which $\lambda_1 = c (d\theta + \rho^2\ d\phi)$.  These are defined by an
embedding $\Psi_0 : S^1 \times B^2_\epsilon(0) \hookrightarrow S^3$
quite similar to \eqref{eqn:goodCoords}.  One can then use a parametrized
version of the Moser deformation argument to
construct for each $K_j$ a family of contact immersions
$\psi_j^\tau : S^1 \times B^2_\epsilon(0) \to S^3$, which are embeddings
near $S^1 \times \{0\}$ for $\tau > 0$, taking $S^1\times \{0\}$ to $K_j^\tau$,
and converge as $\tau \to 0$ to a $k_j$--fold cover of $\Psi_0$.
These define coordinate neighborhoods near $K_j^\tau$ in which $\lambda_1$
is $C^1$--close to something of the form $c (d\theta + \rho^2\ d\phi)$.

Let us now redefine notation and call $K^\tau$ (for sufficiently small
$\tau > 0$) simply $K$; we can then assume there is a contact form
$\lambda_2$ as in Lemma~\ref{lemma:stillTransverse}, taking the form
$c (d\theta + \rho^2\ d\phi)$ in coordinates near each component of $K$.
Choose a smooth homotopy of contact forms $\{\lambda_r\}_{r \in [1,2]}$
such that $\ker\lambda_r = \xi_0$ for all $r$, and each $\lambda_r$ is
$C^1$--close to $\lambda_1$, differing from $\lambda_1$ only in a
tubular neighborhood of $K$.  Observe that the corresponding Reeb vector
fields $X_r$ are are all $C^0$--close to $X_1$, and equal to it outside a
compact neighborhood of $K$.  We may therefore assume $X_r$ is always
transverse to the projection of the foliation $\fF_1$ on
$S^3 \setminus P_\infty$.  As a consequence we have, without loss of 
generality:
\begin{prop}
\label{prop:linking}
Every periodic orbit of $X_r$ that's geometrically distinct from
$P_\infty$ is nontrivially linked with $P_\infty$.
\end{prop}

For $r \in [1,2]$, choose also a smooth homotopy of admissible complex 
multiplications $J_r : \xi_0 \to \xi_0$ such that $J_1 \equiv i$,
$J_r$ differs from $i$ only in a neighborhood of $K$, and $J_2$ is
defined in the coordinates $(\theta,\rho,\phi)$ near each component of
$K$ by a relation of the form $J_2 v_1 = \beta(\rho) v_2$, as in
\S\ref{subsec:local}.
These choices define a smooth homotopy of almost complex structures
$\tilde{J}_r$.  Observe that the binding orbit $P_\infty$ remains a closed
orbit with $\muCZ(P_\infty) = 3$ for all $r$.
We can now use the machinery of \S\ref{sec:MBVP} and
\S\ref{sec:compactness} to show that the foliation $\fF_1$ extends to a
continuous family of foliations for $r \in [1,2]$.  

\begin{prop}
\label{prop:folHomotopy}
For each $r \in [1,2]$, there exists a stable finite energy foliation
$\fF_r$ of $(S^3,\lambda_r,J_r)$ which projects to an open book decomposition
of $S^3$, with binding orbit $P_\infty$.
\end{prop}
\begin{proof}
Denote by $\mM_r$ the moduli space of all $\tilde{J}_r$--holomorphic 
finite energy surfaces, and define the space
$$
\mM = \{ (r,\tilde{u}) \ |\ r \in [1,2],\ \tilde{u} \in \mM_r \}.
$$
The latter has a natural topology induced by the same notion of
convergence as in $\mM_r$, and there are natural continuous inclusions
$\mM_r \hookrightarrow \mM$ for each $r$, as well as a natural
$\RR$--action on $\mM$.  Let $\mM_1^*$ denote the
connected component of $\mM_1$ that contains the planes in the foliation
$\fF_1$, let $\mM^*$ be the corresponding connected component of $\mM$ 
containing $\mM_1^*$, and then define $\mM_r^* = \mM^* \cap \mM_r$.  

Combining Theorems~\ref{thm:IFT}, \ref{thm:parametrizedIFT} 
and~\ref{thm:injective}, we see that for every
$(r,\tilde{u}) \in \mM^*$, $\tilde{u} = (a,u)$ is an embedded index~$2$ plane 
asymptotic to $P_\infty$, and $u : \CC \to S^3$ is also embedded.  Moreover
$\tilde{u}$ is regular, and its neighborhood in $\mM_r^*$ foliates 
neighborhoods of the images of $\tilde{u}$ and $u$.  Clearly then,
$\mM^* / \RR$ is a smooth $2$--dimensional manifold, for which
the projection map
$$
\mM^* / \RR \to \RR : (r,[\tilde{u}]) \mapsto r
$$
is always regular.  
In light of the linking condition in Prop.~\ref{prop:linking},
Theorem~\ref{thm:compactness} implies that $\mM^* / \RR$ is compact.
It follows that there is a diffeomorphism
$$
\psi : [1,2] \times \mM_1^* / \RR \to \mM^* / \RR
$$
such that $\psi(1,[\tilde{u}]) = [\tilde{u}]$ and for each $r \in [1,2]$,
$\psi(r,\cdot)$ is a diffeomorphism $\mM_1^* / \RR \to \mM_r^* / \RR$.
Thus $\mM_r^* / \RR \cong \mM_1^* / \RR \cong S^1$.  Applying
Theorems~\ref{thm:IFT} and~\ref{thm:injective} again, we see that any two
elements of $\mM_r^* / \RR$ are either identical or have disjoint images
in $S^3$.  Moreover, if $\uU_r \subset S^3 \setminus P_\infty$ is
the set of points contained in the image of any curve $[\tilde{u}] \in
\mM_r^* / \RR$, then Theorems~\ref{thm:IFT} and~\ref{thm:compactness}
together imply that $\uU_r$ is open and closed, so 
$\uU_r = S^3 \setminus P_\infty$.  The collection of curves $\mM_r^*$,
together with the trivial cylinder over $P_\infty$, therefore form a
stable finite energy foliation of $(S^3,\lambda_r,J_r)$.
\end{proof}

\subsubsection{$\fF_2 \to \fF_3$: Cutting out disks}
\label{subsubsec:cutting}

For the remainder of \S\ref{subsec:surgery}, we impose the following
restrictive assumption:
\begin{assumption}
\label{assumption}
For each component $K_j \subset K$,
$\lk(K_j,P_\infty) = 1$.
\end{assumption}
This is needed for technical reasons in the arguments that follow, but
will be removed in \S\ref{subsec:braids} by a branched covering argument.

By the above results, we have a stable foliation $\fF_2$ of 
$(S^3,\lambda_2,J_2)$, transverse to a link $K = K_1 \cup \ldots \cup K_m$
whose components have disjoint tubular neighborhoods $N_j \cong
S^1 \times B^2_\epsilon(0)$ on which 
$\lambda_2 = c_j (d\theta + \rho^2\ d\phi)$ and $J_2$ has the form
$J_2 v_1 = \beta_j(\rho) v_2$.  The Reeb vector field on $N_j$ is
$X_2 = \frac{1}{c_j} \p_\theta$.  Pick $\delta \in (0,\epsilon)$ and let 
$N_j^\delta = \{ \rho \le \delta \} \subset N_j$, with
$L_j := \p N_j^\delta$.  Observe that since $L_j$ is
foliated by Reeb orbits, which are necessarily transverse to the leaves
of $\fF_2$, $L_j$ is also transverse to these leaves.

We will now replace the planes in $\fF_2$ with solutions to a boundary value
problem, having boundary mapped to the tori $L_j$.  For this it is necessary
to throw out all except one of the curves in $\fF_2$; we will be able to
reconstruct a foliation afterwards.  Therefore, pick any finite energy plane
$\tilde{u} = (a,u) : \CC \to \RR\times S^3$ which parametrizes a leaf of
$\fF_2$, and define the set of $m$ disjoint open disks
$$
\dD_1 \cup \ldots \cup \dD_m \subset \CC
$$
by $\dD_j = u^{-1}(\interior{N_j^\delta})$.  Then
$$
(\Sigma,j) := (S^2 \setminus (\dD_1 \cup \ldots \cup \dD_m), i)
$$
is a compact Riemann surface with boundary
$$
\p\Sigma = \gamma_1 \cup \ldots \cup \gamma_m,
$$
where $\gamma_j := - \p\overline{\dD}_j$, and we have
$u(\gamma_j) \subset L_j$.  Let $\dot{\Sigma} =
\Sigma\setminus \{\infty\}$.  Observe that due to Assumption~\ref{assumption},
each torus $L_j$ meets the image of a \emph{unique} component $\gamma_j$
under $u$.
Then since each Reeb orbit on $L_j$ has a single
transverse intersection with $u(\gamma_j)$, there are unique smooth
functions $g_j : L_j \to \RR$ such that $d g_j(X_2) \equiv 0$ and
$a(z) = g_j(u(z))$ for all $z \in \gamma_j$.  Thus $\tilde{u}_j$ satisfies
the totally real boundary condition
$$
\tilde{u}(\gamma_j) \subset \tilde{L}_j :=
\{ (g_j(x),x) \in \RR\times S^3\ |\ x \in L_j \}.
$$

This boundary condition is not Lagrangian, so it does not naturally
give rise to any obvious energy bounds.\footnote{It is shown in 
\cite{Wendl:thesis} that one can choose
a new definition of energy so that suitable bounds are satisfied and
the compactness argument goes through.  Here we follow an alternative and
somewhat simpler approach.}  However, the fact that
$d g_j(X_2) \equiv 0$ will allow us to identify each $\tilde{L}_j$ with a 
Lagrangian torus in the symplectization of $S^3$ with a \emph{non-contact} 
stable Hamiltonian structure.  
This is why Assumption~\ref{assumption} is necessary.

\begin{prop}
\label{prop:wavy}
Suppose $M$ is an oriented $3$--manifold with positive contact form
$\lambda$ and Reeb vector field $X$, whose flow is globally defined,
and $J$ is an admissible complex multiplication on 
$\xi = \ker\lambda$ which is preserved by the Reeb flow.  Denote by
$\tilde{J}$ the associated almost complex structure on
$\RR\times M$, and define an $\RR$--equivariant diffeomorphism by
\begin{equation}
\label{eqn:Psi}
\Psi : \RR\times M \to \RR\times M : (a,m) \mapsto (a + F(m), m)
\end{equation}
for some smooth function $F : M \to \RR$ that satisfies
$dF(X) \equiv 0$.  Then if $\xi' \subset TM$ is
the unique $2$--plane distribution in $TM = T(\{0\}\times M) \subset
T(\RR\times M)$ which is preserved by $\Psi_*\tilde{J}$, and
$J' := \Psi_*\tilde{J}|_{\xi'} : \xi' \to \xi'$, the data
$$
\hH' := (\xi', X, d\lambda, J')
$$
define a stable Hamiltonian structure on $M$, for which
the associated almost complex structure is precisely $\Psi_* \tilde{J}$.
\end{prop}
\begin{proof}
Denoting by $\pi : TM \to \xi$ the projection along $X$, define a $1$--form
$$
\lambda' = \lambda - dF \circ J \circ \pi
$$
and let $\xi' = \ker\lambda'$ (we're not assuming this is the same $\xi'$
defined in the statement above).  Clearly $\lambda'(X) \equiv 1$, and
we claim that also $d\lambda'(X,\cdot) \equiv 0$.  Since
$L_X \lambda' = d \iota_X \lambda' + \iota_X d\lambda' = \iota_X d\lambda'$,
this is equivalent to the statement that $\xi'$ is preserved by the
flow of $X$.  Denote this flow by $\varphi^t : M \to M$ and observe that
$\varphi^t_* X \equiv X$ for all $t$, and
by assumption similarly $F \circ \varphi^t \equiv F$ and 
$\varphi^t_* J \equiv J$.  For $m \in M$, any $v \in \xi'_m$
can be written as $v = [ dF(m)J\hat{v} ] X(m) + \hat{v}$ where
$\hat{v} := \pi v \in \xi_m$.  Then
\begin{equation*}
\begin{split}
\varphi^t_* v &=
[ dF(m) J\hat{v} ] X(\varphi^t(m)) + \varphi^t_* \hat{v} \\
&= [ d (F \circ \varphi^t)(m) J\hat{v} ] X(\varphi^t(m)) + \varphi^t_* \hat{v} \\
&= [ dF(\varphi^t(m)) \cdot \varphi^t_* (J\hat{v}) ] X(\varphi^t(m)) + \varphi^t_* \hat{v} \\
&= [ dF(\varphi^t(m)) \cdot J (\varphi^t_*\hat{v}) ] X(\varphi^t(m)) + \varphi^t_* \hat{v}
 \in \xi'_{\varphi^t(m)},
\end{split}
\end{equation*}
proving the claim.

Now observe $d\lambda(X,\cdot) \equiv 0$, and since $\xi'$ is transverse to 
$X$, $d\lambda$ is nondegenerate on $\xi'$ and provides a suitable taming
form for any complex multiplication $J' : \xi' \to \xi'$ with the correct
orientation.  We show next that $\xi'$ is in fact
the unique distribution preserved by $\Psi_*\tilde{J}$.  Indeed, for
$v = [dF(m) J\hat{v}] X(m) + \hat{v} \in \xi'_m$, we have
\begin{equation*}
\begin{split}
(\Psi_*\tilde{J}) v &=
T\Psi \circ \tilde{J} \circ T\Psi^{-1} ([dF(m) J\hat{v}] X(m) + \hat{v}) \\
&= T\Psi \circ \tilde{J} \left( [dF(m) J\hat{v}] X(m) - [dF(m)\hat{v}] \p_a
  + \hat{v} \right) \\
&= T\Psi \left( -[dF(m) J\hat{v}]\p_a - [dF(m)\hat{v}] X(m) + J\hat{v} \right)\\
&= -[dF(m) J\hat{v}] \p_a - [dF(m)\hat{v}] X(m) + [dF(m) J\hat{v}] \p_a + 
  J\hat{v} \\
&= -[dF(m)\hat{v}] X(m) + J\hat{v} \\
&= [dF(m) J (J\hat{v})] X(m) + J\hat{v} \in \xi'_m,
\end{split}
\end{equation*}
thus $\Psi_*\tilde{J}$ restricts on $\xi'$ to the unique map
$J' : \xi' \to \xi'$ such that $J \circ \pi|_{\xi'} \equiv \pi \circ J'$.
Finally, observe that $\Psi_*\tilde{J}$ is clearly $\RR$--invariant, and
since $T\Psi^{-1}$ and $T\Psi$ each preserve both $\p_a$ and $X$,
$$
(\Psi_*\tilde{J}) \p_a = T\Psi \circ \tilde{J} \circ T\Psi^{-1}(\p_a) = X.
$$
\end{proof}

\begin{remark}
\label{remark:homotopy}
There is an obvious smooth homotopy between the two stable Hamiltonian 
structures $\hH_0 := (\xi,X,d\lambda,J)$ and $\hH_1 := 
(\xi',X,d\lambda,J')$: just
define $\hH_\tau = (\xi_\tau,X,d\lambda,J_\tau)$ for $\tau \in [0,1]$ by the 
same trick, but using the functions $F_\tau := \tau F : M \to \RR$.
\end{remark}

\begin{remark}
Finding nontrivial examples of the situation in Prop.~\ref{prop:wavy}
requires very precise knowledge of the Reeb dynamics.
One interesting example is the case where $M$ is a principal $S^1$--bundle
over a Riemann surface with compatible symplectic structure, and
the fibers are generated by $X$ (cf.~\cite{SFTcompactness}*{Example~2.2}).
Then the choice of $F$ determines $\xi'$ as a principal connection on this
bundle.
\end{remark}

We apply the above idea as follows.  Pick a smooth function $F :
S^3 \to \RR$ supported in $N_1 \cup \ldots \cup N_m$, such that 
$dF(X_2) \equiv 0$ and $F(x) = -g_j(x)$ for all
$x \in L_j$.  Then the diffeomorphism \eqref{eqn:Psi} satisfies
$$
\Psi(\tilde{L}_j) = \{0\} \times L_j,
$$
and there is a stable Hamiltonian structure 
$\hH_3 = (\xi_2',X_2,d\lambda_2,J_2')$
with associated almost complex structure $\tilde{J}_3$ such that
$$
\tilde{v} = (b,v) := \Psi \circ \tilde{u} : \dot{\Sigma} \to \RR\times S^3
$$
is $\tilde{J}_3$--holomorphic and satisfies the Lagrangian boundary
condition $\tilde{v}(\gamma_j) \subset \{0\}\times L_j$.
Thus writing $\Lambda = (\{0\}\times L_1,\ldots,\{0\}\times L_m)$, we have
$\tilde{v} \in \mM_{\hH_3,\Lambda}$; in fact $\tilde{v}$
satisfies the same assumptions as the sequence in our main compactness
theorem, so Lemma~\ref{lemma:indexuk} implies $\ind(\tilde{v}) = 2$.
Then by Theorem~\ref{thm:IFT}, the connected component $\mM_3^* \subset
\mM_{\hH_3,\Lambda}$ containing $\tilde{v}$ is a smooth $2$--manifold with
free and proper $\RR$--action, so $\mM_3^* / \RR$ is a smooth $1$--manifold,
and Theorem~\ref{thm:compactness} implies it is compact, i.e.~it is
diffeomorphic to $S^1$.  Arguing again as in Prop.~\ref{prop:folHomotopy},
we find that the curves in $\mM_3^*$ form an $\RR$--invariant 
foliation $\fF_3$ of $\RR\times (M\setminus P_\infty)$, where
$$
M := S^3 \setminus \interior( N_1^\delta \cup \ldots \cup N_m^\delta ).
$$
It projects to a smooth foliation of $M \setminus
P_\infty$ by an $S^1$--parametrized family of leaves asymptotic to $P_\infty$
and transverse to $\p M = L_1 \cup \ldots \cup L_m$.  

\begin{defn}
\label{defn:openbookBndry}
A foliation with the properties named above is called a
\emph{stable holomorphic open book decomposition with boundary}.
\end{defn}

\subsubsection{$\fF_3 \to \fF_4$: Twisting}
\label{subsubsec:twisting}

It follows from Remark~\ref{remark:homotopy} that there is a smooth
homotopy of stable Hamiltonian structures $\hH_r$ for $r \in [3,7/2]$,
deforming $\hH_3$ back to the original contact data
$\hH_{7/2} := (\xi_2,X_2,d\lambda_2,J_2)$.  Finally, continue this homotopy
for $r \in [7/2,4]$ by choosing a smooth family of contact forms
$\lambda_r$ on $M$ such that
\begin{enumerate}
\setlength{\itemsep}{0in}
\item $\lambda_r \equiv \lambda_2$ outside coordinate neighborhoods of the
tori $L_j$.
\item In a coordinate neighborhood near $L_j$, $\lambda_r = f_r(\rho)\ d\theta
 + g_r(\rho)\ d\phi$ where $g_r'(\rho) > 0$ for all $r < 4$ and 
 $\rho \ge \delta$, but $g_4'(\delta) = 0$.
\item $f_4$ and $g_4$ extend for $\rho$ in some open neighborhood of $\delta$
such that $g_4''(\delta) > 0$ and $f_4'(\delta) > 0$.
\end{enumerate}
These conditions guarantee that all closed Reeb orbits for $r < 4$ are
nontrivially linked with $P_\infty$, and this remains true at $r = 4$
in $\interior{M}$, but the boundary components $L_j$ then become
Morse-Bott tori with closed orbits forming negatively oriented meridians.
Notice that no such homotopy of contact forms exists globally on $S^3$;
this is why we
introduced the boundary condition, to remove the interiors of
$N_j^\delta$ from the picture.

For any $r_0 \in (3,4)$, we can apply the arguments of 
Prop.~\ref{prop:folHomotopy} and find a continuous family $\fF_r$ of stable
holomorphic open book decompositions with boundary for $r \in [3,r_0]$.
Then taking $r \to 4$, the degeneration theorem~\ref{thm:degeneration}
gives limits in the form of finite energy surfaces without boundary,
having $m+1$ positive punctures asymptotic to $P_\infty$ and the simply
covered Morse-Bott orbits on $N_1^\delta,\ldots,N_m^\delta$.  
In particular, for any $m \in M \setminus (P_\infty \cup \p M)$, 
there exists such a curve
$\tilde{v}_\infty = (b_\infty,v_\infty)$ with $m$ in the image of $v_\infty$: 
it is obtained by taking sequences of
corresponding curves in $\fF_r$ passing through $m$ and letting
$r$ approach~$4$.  By positivity of intersections,
the limit curves are also embedded, and any pair of them has projections
that are either identical or disjoint in $M$.  Finally, a simple
computation using \eqref{eqn:FredholmIndexBasic} and 
Lemma~\ref{lemma:MBparity} shows that these curves have index~$2$.
They therefore constitute a stable finite energy foliation of Morse-Bott
type on the manifold with boundary $M$.  This, together with the
constructions in \S\ref{subsec:local}, proves the main result for any
situation in which Assumption~\ref{assumption} is satisfied.

\subsection{Lifting to general closed braids}
\label{subsec:braids}

We now complete the proof of Theorem~\ref{thm:mainresult} by constructing
foliations in cases where Assumption~\ref{assumption} does not hold.
The idea is to define a branched cover over $S^3$ so that the assumption
does hold on the cover, thus the previous arguments produce a foliation,
which we will then show has a well defined projection.

By way of preparation, define the usual cylindrical coordinates
$(\theta,\rho,\phi)$ on $M := S^1 \times B^2_\epsilon(0)$,
pick $n \in \NN$ and consider the map
\begin{equation}
\label{eqn:localCovering}
p : S^1 \times B^2_\epsilon(0) \to S^1 \times B^2_\epsilon(0) :
(\theta,\rho,\phi) \mapsto (\theta,\rho,n\phi).
\end{equation}
Writing $P := S^1 \times \{0\}$, this map is smooth on 
$M \setminus P$ and continuous everywhere.  Suppose $M$ is endowed with
a smooth contact form of type $\lambda = f(\rho)\ d\theta + g(\rho)\ d\phi$
and admissible complex multiplication defined by $J v_1 = \beta(\rho)v_2$
as in \S\ref{subsec:local}.  Then on $M\setminus P$, $\lambda$ pulls back
to another contact form
\begin{equation}
\label{eqn:lambdan}
\lambda^{(n)} := p^*\lambda = f_n(\rho)\ d\theta + g_n(\rho)\ d\phi
= f(\rho)\ d\theta + n g(\rho)\ d\phi,
\end{equation}
which extends smoothly to $P$.  In fact, let $\Phi_0$ be the trivialization
of $\xi = \ker\lambda$ along $P$ defined by these coordinates, and
suppose $P$ is a nondegenerate Reeb orbit for $\lambda$ with 
$\muCZ^{\Phi_0}(P) = 1$.  From Prop.~\ref{prop:symOrbits}, this means
$-f''(0) / 2\pi g''(0) \in (0,1)$, and thus the same is true of
$-f_n''(0) / 2\pi g_n''(0)$, so $P$ also has $\muCZ^{\Phi_0}(P) = 1$
with respect to the extended contact form $\lambda^{(n)}$.
Writing our standard symplectic frame on $\xi^{(n)} = \ker\lambda^{(n)}$
as $(v_1^{(n)},v_2^{(n)})$, the complex multiplication $J$ also 
pulls back on $M\setminus P$ to $J^{(n)} = p^* J$, satisfying 
$J^{(n)} v_1^{(n)} = \beta_n(\rho) v_2^{(n)}$ where
$\beta_n(\rho) = \beta(\rho)$.  It turns out that $J^{(n)}$ does \emph{not}
have a smooth extension over $P$, but this will be only a minor irritation
in the following.

If Assumption~\ref{assumption} does not hold, choose $n$ to be the
least common multiple of all the linking numbers
$\lk(K_j,P_\infty)$, and define an $n$--fold branched cover of $S^3$
as follows.  By the results of \S\ref{subsubsec:simplify}, there is a
contact form $\lambda_2$ and complex multiplication $J_2$, both of which
take the usual simple forms in the neighborhoods $N_j$ of $K_j$, and
$(S^3,\lambda_2,J_2)$ admits a stable finite energy foliation $\fF_2$
which projects to an open book decomposition of $S^3$ with binding
orbit $P_\infty$.  Denoting $E := S^3 \setminus P_\infty$, this open book
defines a fibration $E \to S^1$, and there is a natural $n$--fold covering
map $p : E^{(n)} \to E$ and smooth fibration $E^{(n)} \to S^1$ such that
$p(E_\tau^{(n)}) = E_{n\tau}$ for each $\tau \in S^1$.  Let $\psi :
E^{(n)} \to E^{(n)}$ be the deck transformation which maps
$$
E_\tau^{(n)} \to E_{\tau + \frac{1}{n}}^{(n)};
$$
then every deck transformation is of the form $\psi^k$ for $k \in \ZZ_n$.

In order to compactify $E^{(n)}$, we shall modify this construction
carefully near $P_\infty$.  Recall from \eqref{eqn:goodCoords} that
a neighborhood of $P_\infty$ admits cylindrical polar coordinates
$(\theta,\rho,\phi)$ such that $P_\infty \cong S^1\times \{0\}$ and
$\lambda_2$ takes the form
$$
\lambda_2 = f(\rho)\ d\theta + g(\rho)\ d\phi,
$$
where $f$ and $g$ are smooth functions with $f'g'' - f''g' \equiv 0$
near~$0$.  Moreover $J_2$ is defined in this neighborhood by a relation
of the form $J_2 v_1 = \beta(\rho)v_2$.  These properties continue to
hold if we change the coordinates by any diffeomorphism of the form 
$(\theta,\rho,\phi) \longleftrightarrow (\theta,\rho,\phi + 2\pi k \theta)$ 
for $k \in \ZZ$, thus we
can assume without loss of generality that the planes in $\fF_2$ have
trivial winding around $P_\infty$ as they approach it in these coordinates.
In this case the coordinates define a trivialization $\Phi_0$ in which
$\muCZ^{\Phi_0}(P_\infty) = 1$.
We can now replace $\fF_2$ with another (homotopic) open book decomposition
whose pages look like $\{ \phi = \text{const} \}$ in a neighborhood of
$P_\infty$.  Then applying the covering construction above,
$E^{(n)}$ admits a compactification
$$
\overline{E}^{(n)} = E^{(n)} \cup P_\infty^{(n)} \cong S^3,
$$
where $P_\infty^{(n)}$ is a circle identified with $S^1\times \{0\}$ in
a certain cylindrical coordinate neighborhood, and
$p$ maps this neighborhood to a
neighborhood of $P_\infty$ via \eqref{eqn:localCovering}.

There is thus a continuous extension $p : \overline{E}^{(n)} \to S^3$,
such that $\lambda^{(n)} = p^*\lambda_2$ has the form \eqref{eqn:lambdan}
and thus extends smoothly over $\overline{E}^{(n)}$; in fact for this
extension, $P_\infty^{(n)}$ is a nondegenerate Reeb orbit with
$\muCZ^{\Phi_0}(P_\infty^{(n)}) = 1$.  The lift $J^{(n)} = p^*J_2$ is
uniquely defined over $E^{(n)}$ but singular at $P_\infty^{(n)}$.
These are both preserved by the deck transformation $\psi$, and together
they define an almost complex structure $\tilde{J}^{(n)}$ on
$\RR\times E^{(n)} \subset \RR\times \overline{E}^{(n)}$, such that
the diffeomorphism
$$
\tilde{\psi} : \RR\times E^{(n)} \to \RR\times E^{(n)} :
(a,m) \mapsto (a,\psi(m))
$$
is $\tilde{J}^{(n)}$--holomorphic.  Each leaf $\tilde{u} \in \fF_2$
lifts to an embedded $\tilde{J}^{(n)}$--holomorphic plane
$\CC \to \RR\times E^{(n)}$, giving a foliation of $E^{(n)}$ by planes
asymptotic to $P_\infty^{(n)}$.  These can be made into an honest
stable foliation by changing $\beta_n(\rho)$ for $\rho$ near~$0$ so that
$J^{(n)}$ becomes smooth.  This is possible by Prop.~\ref{prop:continuation},
because $f_n' g_n'' - f_n'' g_n' \equiv 0$ in this neighborhood:
thus for a suitable smooth choice of $J^{(n)}$, we obtain a stable
foliation $\fF_2^{(n)}$ which matches the original lift of $\fF_2$
outside some neighborhood of $P_\infty^{(n)}$.

\begin{figure}
\begin{center}
\includegraphics{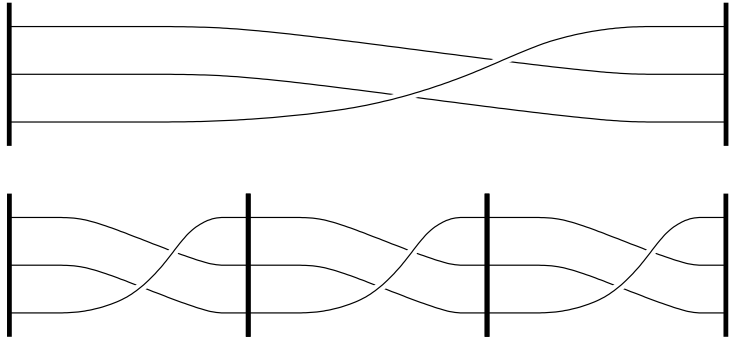}
\end{center}
\caption[A closed braid with one component and its $3$--fold 
cover]{\label{fig:braid}
The top is a transverse knot $K$ with $\lk(K,P_\infty) = 3$, represented as
a closed braid.  The bottom is the $3$--fold cover $K^{(3)} \subset
E^{(3)}$, with three components cyclically permuted by $\psi$.}
\end{figure}

Here is the main point: the link $K \subset E$ lifts to
another transverse link $K^{(n)} = p^{-1}(K) \subset E^{(n)}$,
and our choice of $n$ guarantees that every connected component
$K^{(n)}_j \subset K^{(n)}$ satisfy $\lk(K^{(n)}_j, P_\infty^{(n)}) = 1$
(see Fig.~\ref{fig:braid}).
Thus the arguments of \S\ref{subsubsec:cutting}
and \S\ref{subsubsec:twisting} produce a stable finite energy
foliation $\fF_4$ of Morse-Bott type on $(M^{(n)},\lambda^{(n)}_4,
J^{(n)}_4)$, where $M^{(n)}$ is the complement of a neighborhood
$N^{(n)}$ of $K^{(n)}$ in $\overline{E}^{(n)} \cong S^3$.
We can easily arrange moreover that $N^{(n)}$, 
$\lambda^{(n)}_4$ and $J^{(n)}_4$ 
be invariant under the deck transformation $\psi$, so the diffeomorphisms
$\tilde{\psi}^k$ for $k \in \ZZ_n$ are $\tilde{J}^{(n)}_4$--holomorphic.
This gives rise to a set of $n$ stable foliations
$\fF^k_4 := \tilde{\psi}^k(\fF_4)$ for $k \in \ZZ_n$.

\begin{prop}
\label{prop:lift}
The foliations $\fF^k_4$ are all identical.
\end{prop}
\begin{proof}
It suffices to show that for any leaf $\tilde{u} = (a,u) \in \fF_4$, the curve 
$$
\tilde{\psi} \circ \tilde{u} = (a, \psi\circ u ) : 
\dot{\Sigma} \to \RR\times M^{(n)}
$$
is also a leaf of the foliation.
This follows from positivity of intersections.
Indeed, if $\tilde{\psi} \circ \tilde{u}$ is not a leaf, it must have
finitely many isolated intersections with some other leaf
$\tilde{v} = (b,v) \in \fF_4$, and these
cannot be eliminated under homotopies.  Thus $\tilde{\psi}
\circ \tilde{u}$ also has isolated intersections with $\tilde{u}$.
Now applying an $\RR$--translation, $\tilde{\psi}\circ \tilde{u}$ also
intersects $\tilde{u}^\sigma := (a + \sigma, u)$ for all $\sigma \in \RR$.
Since $a : \dot{\Sigma} \to \RR$ is a proper map, choosing $\sigma$ large
forces these intersections toward the asymptotic limits.  But $\tilde{u}$
and $\tilde{\psi}\circ\tilde{u}$ clearly have distinct asymptotic limits,
and neither curve intersects the asymptotic limits of the other, thus we
have a contradiction.
\end{proof}

It follows that $\fF_4$ has a well defined projection under $p$ to an
$\RR$--invariant foliation of $\RR \times S^3$, inducing also a 
foliation of $S^3\setminus P_\infty$ by planes asymptotic to $P_\infty$.
The projection of $J_4$ is singular at $P_\infty$, but this can again
be fixed by an application of Prop.~\ref{prop:continuation}.
The proof of Theorem~\ref{thm:mainresult} is now complete.

\end{document}